\documentclass[12pt,diff_eng]{article} 
\usepackage[english]{babel}
\usepackage{amssymb,amscd,amsmath}     

\PassOptionsToPackage{linktocpage}{hyperref}

\begin{document}

\large  

%%%%%%%%%%%%%%%%%%%%%%%%%%%%%%%%%%%%%%%%%%%%%%%%%%%%%%%%%%%%%
%%% insert the name of article file here
%%%%

\bigskip

$$
{}
$$
$$
{}
$$

\vspace{2mm}

\begin{center}
\Large{\bf A new proof of the expansion 
of iterated It\^{o} stochastic integrals with
respect to the components of a multidimensional Wiener process
based on generalized multiple Fourier series
and Hermite polynomials}
\end{center}

$$
{}
$$

\centerline{\large {Dmitriy F. Kuznetsov}}
$$
{}
$$
\begin{center}
Peter the Great Saint-Petersburg Polytechnic University\\
Polytechnicheskaya ul., 29,\\
195251, Saint-Petersburg, Russia\\
\hspace{-2mm}e-mail: sde\_\hspace{0.5mm}kuznetsov@inbox.ru
\end{center}

$$
{}
$$
$$
{}
$$

\noindent
{\bf Abstract.} The article is devoted to a new proof of the expansion 
for iterated It\^{o} stochastic integrals with
respect to the components of a multidimensional Wiener process.
The above expansion is
based on Hermite polynomials and generalized multiple Fourier series in arbitrary 
complete orthonormal systems of functions in a Hilbert space.
In 2006, the author obtained a similar expansion, but with a lesser degree
of generality, namely, for the case of continuous or piecewise continuous 
complete orthonornal systems of functions in a Hilbert space.
In this article, the author generalizes the expansion of iterated
It\^{o} stochastic integrals obtained by him in 2006 to the case
of an arbitrary 
complete orthonormal systems of functions in a Hilbert space
using a new approach based on the It\^{o} formula.
The obtained expansion of iterated It\^{o} stochastic integrals is useful
for constructing of high-order strong numerical methods
for systems of It\^{o} stochastic differential equations with
multidimensional non-commutative noise.

\renewcommand{\baselinestretch}{1.3}

\vspace{7mm}

{\normalsize 

\linespread{1.0}

\tableofcontents

\vspace{8mm}

\linespread{1.0}

}

%{\normalsize 

%\setlength{\baselineskip}{1.5em}

%\tableofcontents

%\setlength{\baselineskip}{1.2em}

\vspace{5mm}

\section{Introduction}

\vspace{3mm}

Let $(\Omega,$ ${\rm F},$ ${\sf P})$ be a complete probability space, let 
$\{{\rm F}_{\tau}, \tau\in[0,\bar T]\}$ be a nondecreasing 
right-continous family of $\sigma$-algebras of ${\rm F},$
and let ${\bf w}_{\tau}$ be a standard $m$-dimensional Wiener 
stochastic process, which is
${\rm F}_{\tau}$-measurable for any $\tau\in[0, \bar T].$ We assume that the components
${\bf w}_{\tau}^{(i)}$ $(i=1,\ldots,m)$ of this process 
are independent. 
Consider
an It\^{o} stochastic differential equation (SDE) in the integral form
\begin{equation}
\label{1.5.2}
{\bf x}_s={\bf x}_0+\int\limits_0^s {\bf a}({\bf x}_{\tau},\tau)d\tau+
\sum\limits_{j=1}^m \int\limits_0^s B_j({\bf x}_{\tau},\tau)d{\bf w}_{\tau}^{(j)},\ \ \
{\bf x}_0={\bf x}(0,\omega),\ \ \ \omega\in\Omega.
\end{equation}
Here ${\bf x}_s$ is $n$-dimensional stochastic process 
satisfying the equation (\ref{1.5.2}). 
The nonrandom functions ${\bf a}({\bf x},\tau),$ $B_j({\bf x},\tau):\ 
{\bf R}^n\times[0, \bar T]\to{\bf R}^n$ $(j=1,\dots,m)$
guarantee the existence and uniqueness up to stochastic equivalence 
of the strong solution
of equation (\ref{1.5.2}) \cite{1}. The second integral on the 
right-hand side of (\ref{1.5.2}) is 
the It\^{o} stochastic integral.
Let ${\bf x}_0$ be an $n$-dimensional random variable, which is 
${\rm F}_0$-measurable and 
${\sf M}\bigl\{\left|{\bf x}_0\right|^2\bigr\}<\infty$ 
(${\sf M}$ denotes a mathematical expectation).
We assume that
${\bf x}_0$ and ${\bf w}_{\tau}-{\bf w}_0$ are independent when $\tau>0.$
In addition to the above conditions, we will assume that the functions
${\bf a}({\bf x},\tau),$ $B_j({\bf x},\tau)$ $(j=1,\ldots,m)$ are sufficiently 
smooth functions in both arguments.

It is well known \cite{KlPl2}-\cite{KPS}
that It\^{o} SDEs are 
adequate mathematical models of dynamic systems of various 
physical nature under 
the influence of random disturbances. One of the effective approaches 
to the numerical integration of 
It\^{o} SDEs is an approach based on 
the Taylor--It\^{o} and 
Taylor--Stratonovich expansions
\cite{KlPl2}-\cite{2018a}. 
The most important feature of such 
expansions is a presence in them of the so-called iterated
It\^{o} and Stratonovich stochastic integrals, which play the key 
role for solving the 
problem of numerical integration of It\^{o} SDEs
and have the 
following form 
\begin{equation}
\label{ito}
J[\psi^{(k)}]_{T,t}^{(i_1 \ldots i_k)}=\int\limits_t^T\psi_k(t_k) \ldots 
\int\limits_t^{t_{2}}
\psi_1(t_1) d{\bf w}_{t_1}^{(i_1)}\ldots
d{\bf w}_{t_k}^{(i_k)},
\end{equation}
\begin{equation}
\label{str}
J^{*}[\psi^{(k)}]_{T,t}^{(i_1 \ldots i_k)}=
{\int\limits_t^{*}}^T
\psi_k(t_k) \ldots 
{\int\limits_t^{*}}^{t_2}
\psi_1(t_1) d{\bf w}_{t_1}^{(i_1)}\ldots
d{\bf w}_{t_k}^{(i_k)},
\end{equation}
where $\psi_1(\tau),\ldots,\psi_k(\tau)$ are nonrandom
functions 
on $[t,T],$ ${\bf w}_{\tau}^{(i)}$ ($i=1,\ldots,m$) are
independent standard Wiener processes and
${\bf w}_{\tau}^{(0)}=\tau,$
$i_1,\ldots,i_k = 0,$ $1,$ $\ldots,$ $m,$
$$
\int\limits\ \hbox{and}\ \int\limits^{*}
$$ 
denote It\^{o} and 
Stratonovich stochastic integrals,
respectively.

Generalization of the method 
of expansion of iterated It\^{o} stochastic integrals (\ref{ito})
based on generalized multiple Fourier series
(see Theorem~5.1 (\cite{2006}, p.~236) and
Sect.~5.1 (\cite{2006}, pp.~235-245))
composes the subject of the article.

Note that another approaches to the expansion and mean-square approximation
of iterated It\^{o} and Stratonovich
stochastic integrals 
(\ref{ito}) and (\ref{str}) can be found in 
\cite{KlPl2}-\cite{KPS},
\cite{allen}-\cite{Zapad-9}.

Suppose that $\psi_1(\tau),\ldots,\psi_k(\tau)
\in L_2([t,T])$.
Define the following function (the so-called factorized Volterra-type kernel)
on the hypercube $[t, T]^k$
\begin{equation}
\label{ppp}
K(t_1,\ldots,t_k)=
\left\{\begin{matrix}
\psi_1(t_1)\ldots \psi_k(t_k),\ &t_1<\ldots<t_k\cr\cr
0,\ &\hbox{\rm otherwise}
\end{matrix}
\right.,
\end{equation}

\noindent
where $t_1,\ldots,t_k\in [t, T]$ $(k\ge 2)$ and 
$K(t_1)\equiv\psi_1(t_1)$ for $t_1\in[t, T].$

Suppose that $\{\phi_j(x)\}_{j=0}^{\infty}$
is an arbitrary  complete orthonormal system of functions in 
the space $L_2([t, T])$.

It is well known that the generalized 
multiple Fourier series 
of $K(t_1,\ldots,t_k)$ $\in L_2([t, T]^k)$ is converging 
to $K(t_1,\ldots,t_k)$ in the hypercube $[t, T]^k$ in 
the mean-square sense, i.e.

\vspace{-2mm}
\begin{equation}
\label{sos1z}
\hbox{\vtop{\offinterlineskip\halign{
\hfil#\hfil\cr
{\rm lim}\cr
$\stackrel{}{{}_{p_1,\ldots,p_k\to \infty}}$\cr
}} }\biggl\Vert
K-K_{p_1\ldots p_k}\biggr\Vert_{L_2([t, T]^k)}=0,
\end{equation}

\vspace{1mm}
\noindent
where
\begin{equation}
\label{chain30001x}
K_{p_1\ldots p_k}(t_1,\ldots,t_k)=
\sum_{j_1=0}^{p_1}\ldots
\sum_{j_k=0}^{p_k}
C_{j_k\ldots j_1}
\prod_{l=1}^k\phi_{j_l}(t_l),
\end{equation}

\vspace{1mm}
\begin{equation}
\label{ppppa}
C_{j_k\ldots j_1}=\int\limits_{[t,T]^k}
K(t_1,\ldots,t_k)\prod_{l=1}^{k}\phi_{j_l}(t_l)dt_1\ldots dt_k
\end{equation}

\vspace{1mm}
\noindent
is the Fourier coefficient, and
$$
\left\Vert f\right\Vert_{L_2([t, T]^k)}=\left(~\int\limits_{[t,T]^k}
f^2(t_1,\ldots,t_k)dt_1\ldots dt_k\right)^{1/2}.
$$

\vspace{1mm}

Consider the partition $\{\tau_j\}_{j=0}^N$ of $[t,T]$ such that
\begin{equation}
\label{1111}
t=\tau_0<\ldots <\tau_N=T,\ \ \
\Delta_N=
\hbox{\vtop{\offinterlineskip\halign{
\hfil#\hfil\cr
{\rm max}\cr
$\stackrel{}{{}_{0\le j\le N-1}}$\cr
}} }\Delta\tau_j\to 0\ \ \hbox{if}\ \ N\to \infty,\ \ \ 
\Delta\tau_j=\tau_{j+1}-\tau_j.
\end{equation}

{\bf Theorem 1}\ \cite{2006} (2006) (also see \cite{2007}-\cite{29a}).\
{\it Suppose that
$\psi_1(\tau),\ldots,\psi_k(\tau)$ are continuous 
nonrandom functions on 
$[t, T]$ and
$\{\phi_j(x)\}_{j=0}^{\infty}$ is an arbitrary complete orthonormal system  
of continuous or piecewise continuous functions in the space $L_2([t,T]).$ 
Then
$$
J[\psi^{(k)}]_{T,t}^{(i_1\ldots i_k)} =
\hbox{\vtop{\offinterlineskip\halign{
\hfil#\hfil\cr
{\rm l.i.m.}\cr
$\stackrel{}{{}_{p_1,\ldots,p_k\to \infty}}$\cr
}} }\sum_{j_1=0}^{p_1}\ldots\sum_{j_k=0}^{p_k}
C_{j_k\ldots j_1}\Biggl(
\prod_{l=1}^k\zeta_{j_l}^{(i_l)} -
\Biggr.
$$
\begin{equation}
\label{tyyy}
-\Biggl.
\hbox{\vtop{\offinterlineskip\halign{
\hfil#\hfil\cr
{\rm l.i.m.}\cr
$\stackrel{}{{}_{N\to \infty}}$\cr
}} }\sum_{(l_1,\ldots,l_k)\in {\rm G}_k}
\phi_{j_{1}}(\tau_{l_1})
\Delta{\bf w}_{\tau_{l_1}}^{(i_1)}\ldots
\phi_{j_{k}}(\tau_{l_k})
\Delta{\bf w}_{\tau_{l_k}}^{(i_k)}\Biggr),
\end{equation}
where
$$
{\rm G}_k={\rm H}_k\backslash{\rm L}_k,\ \ \
{\rm H}_k=\bigl\{(l_1,\ldots,l_k):\ l_1,\ldots,l_k=0,\ 1,\ldots,N-1\bigr\},
$$
$$
{\rm L}_k=\bigl\{(l_1,\ldots,l_k):\ l_1,\ldots,l_k=0,\ 1,\ldots,N-1;\
l_g\ne l_r\ (g\ne r);\ g, r=1,\ldots,k\bigr\},
$$

\noindent
${\rm l.i.m.}$ is a limit in the mean-square sense$,$
$i_1,\ldots,i_k=0,1,\ldots,m,$ 
\begin{equation}
\label{rr23}
\zeta_{j}^{(i)}=
\int\limits_t^T \phi_{j}(\tau) d{\bf w}_{\tau}^{(i)}
\end{equation} 
are independent standard Gaussian random variables
for various
$i$ or $j$ {\rm(}in the case when $i\ne 0${\rm),}
$C_{j_k\ldots j_1}$ is the Fourier coefficient {\rm(\ref{ppppa}),}
$\Delta{\bf w}_{\tau_{j}}^{(i)}=
{\bf w}_{\tau_{j+1}}^{(i)}-{\bf w}_{\tau_{j}}^{(i)}$
$(i=0, 1,\ldots,m),$\
$\left\{\tau_{j}\right\}_{j=0}^{N}$ is a partition of
$[t,T],$ which satisfies the condition {\rm (\ref{1111})}.}

A number of generalizations and modifications of Theorem 1 
can be found in \cite{2018a}, Chapter~1 (see also bibliography therein).

Let us consider corollaries from Theorem 1 
(see (\ref{tyyy})) for
$k=1,\ldots,5$ \cite{2006} 

\begin{equation}
\label{a1}
J[\psi^{(1)}]_{T,t}^{(i_1)}
=\hbox{\vtop{\offinterlineskip\halign{
\hfil#\hfil\cr
{\rm l.i.m.}\cr
$\stackrel{}{{}_{p_1\to \infty}}$\cr
}} }\sum_{j_1=0}^{p_1}
C_{j_1}\zeta_{j_1}^{(i_1)},
\end{equation}

\begin{equation}
\label{a2}
J[\psi^{(2)}]_{T,t}^{(i_1i_2)}
=\hbox{\vtop{\offinterlineskip\halign{
\hfil#\hfil\cr
{\rm l.i.m.}\cr
$\stackrel{}{{}_{p_1,p_2\to \infty}}$\cr
}} }\sum_{j_1=0}^{p_1}\sum_{j_2=0}^{p_2}
C_{j_2j_1}\Biggl(\zeta_{j_1}^{(i_1)}\zeta_{j_2}^{(i_2)}
-{\bf 1}_{\{i_1=i_2\ne 0\}}
{\bf 1}_{\{j_1=j_2\}}\Biggr),
\end{equation}
$$
J[\psi^{(3)}]_{T,t}^{(i_1 i_2 i_3)}=
\hbox{\vtop{\offinterlineskip\halign{
\hfil#\hfil\cr
{\rm l.i.m.}\cr
$\stackrel{}{{}_{p_1,p_2,p_3\to \infty}}$\cr
}} }\sum_{j_1=0}^{p_1}\sum_{j_2=0}^{p_2}\sum_{j_3=0}^{p_3}
C_{j_3j_2j_1}\Biggl(
\zeta_{j_1}^{(i_1)}\zeta_{j_2}^{(i_2)}\zeta_{j_3}^{(i_3)}
-\Biggr.
$$
\begin{equation}
\label{a3}
\Biggl.-{\bf 1}_{\{i_1=i_2\ne 0\}}
{\bf 1}_{\{j_1=j_2\}}
\zeta_{j_3}^{(i_3)}
-{\bf 1}_{\{i_2=i_3\ne 0\}}
{\bf 1}_{\{j_2=j_3\}}
\zeta_{j_1}^{(i_1)}-
{\bf 1}_{\{i_1=i_3\ne 0\}}
{\bf 1}_{\{j_1=j_3\}}
\zeta_{j_2}^{(i_2)}\Biggr),
\end{equation}

\vspace{3mm}
$$
J[\psi^{(4)}]_{T,t}^{(i_1 \ldots i_4)}
=
\hbox{\vtop{\offinterlineskip\halign{
\hfil#\hfil\cr
{\rm l.i.m.}\cr
$\stackrel{}{{}_{p_1,\ldots,p_4\to \infty}}$\cr
}} }\sum_{j_1=0}^{p_1}\ldots\sum_{j_4=0}^{p_4}
C_{j_4\ldots j_1}\Biggl(
\prod_{l=1}^4\zeta_{j_l}^{(i_l)}
\Biggr.
-
$$
$$
-
{\bf 1}_{\{i_1=i_2\ne 0\}}
{\bf 1}_{\{j_1=j_2\}}
\zeta_{j_3}^{(i_3)}
\zeta_{j_4}^{(i_4)}
-
{\bf 1}_{\{i_1=i_3\ne 0\}}
{\bf 1}_{\{j_1=j_3\}}
\zeta_{j_2}^{(i_2)}
\zeta_{j_4}^{(i_4)}-
$$
$$
-
{\bf 1}_{\{i_1=i_4\ne 0\}}
{\bf 1}_{\{j_1=j_4\}}
\zeta_{j_2}^{(i_2)}
\zeta_{j_3}^{(i_3)}
-
{\bf 1}_{\{i_2=i_3\ne 0\}}
{\bf 1}_{\{j_2=j_3\}}
\zeta_{j_1}^{(i_1)}
\zeta_{j_4}^{(i_4)}-
$$
$$
-
{\bf 1}_{\{i_2=i_4\ne 0\}}
{\bf 1}_{\{j_2=j_4\}}
\zeta_{j_1}^{(i_1)}
\zeta_{j_3}^{(i_3)}
-
{\bf 1}_{\{i_3=i_4\ne 0\}}
{\bf 1}_{\{j_3=j_4\}}
\zeta_{j_1}^{(i_1)}
\zeta_{j_2}^{(i_2)}+
$$
$$
+
{\bf 1}_{\{i_1=i_2\ne 0\}}
{\bf 1}_{\{j_1=j_2\}}
{\bf 1}_{\{i_3=i_4\ne 0\}}
{\bf 1}_{\{j_3=j_4\}}
+
{\bf 1}_{\{i_1=i_3\ne 0\}}
{\bf 1}_{\{j_1=j_3\}}
{\bf 1}_{\{i_2=i_4\ne 0\}}
{\bf 1}_{\{j_2=j_4\}}+
$$
\begin{equation}
\label{a4}
+\Biggl.
{\bf 1}_{\{i_1=i_4\ne 0\}}
{\bf 1}_{\{j_1=j_4\}}
{\bf 1}_{\{i_2=i_3\ne 0\}}
{\bf 1}_{\{j_2=j_3\}}\Biggr),
\end{equation}

\vspace{3mm}
$$
J[\psi^{(5)}]_{T,t}^{(i_1 \ldots i_5)}
=\hbox{\vtop{\offinterlineskip\halign{
\hfil#\hfil\cr
{\rm l.i.m.}\cr
$\stackrel{}{{}_{p_1,\ldots,p_5\to \infty}}$\cr
}} }\sum_{j_1=0}^{p_1}\ldots\sum_{j_5=0}^{p_5}
C_{j_5\ldots j_1}\Biggl(
\prod_{l=1}^5\zeta_{j_l}^{(i_l)}
-\Biggr.
$$
$$
-
{\bf 1}_{\{i_1=i_2\ne 0\}}
{\bf 1}_{\{j_1=j_2\}}
\zeta_{j_3}^{(i_3)}
\zeta_{j_4}^{(i_4)}
\zeta_{j_5}^{(i_5)}-
{\bf 1}_{\{i_1=i_3\ne 0\}}
{\bf 1}_{\{j_1=j_3\}}
\zeta_{j_2}^{(i_2)}
\zeta_{j_4}^{(i_4)}
\zeta_{j_5}^{(i_5)}-
$$
$$
-
{\bf 1}_{\{i_1=i_4\ne 0\}}
{\bf 1}_{\{j_1=j_4\}}
\zeta_{j_2}^{(i_2)}
\zeta_{j_3}^{(i_3)}
\zeta_{j_5}^{(i_5)}-
{\bf 1}_{\{i_1=i_5\ne 0\}}
{\bf 1}_{\{j_1=j_5\}}
\zeta_{j_2}^{(i_2)}
\zeta_{j_3}^{(i_3)}
\zeta_{j_4}^{(i_4)}-
$$
$$
-
{\bf 1}_{\{i_2=i_3\ne 0\}}
{\bf 1}_{\{j_2=j_3\}}
\zeta_{j_1}^{(i_1)}
\zeta_{j_4}^{(i_4)}
\zeta_{j_5}^{(i_5)}-
{\bf 1}_{\{i_2=i_4\ne 0\}}
{\bf 1}_{\{j_2=j_4\}}
\zeta_{j_1}^{(i_1)}
\zeta_{j_3}^{(i_3)}
\zeta_{j_5}^{(i_5)}-
$$
$$
-
{\bf 1}_{\{i_2=i_5\ne 0\}}
{\bf 1}_{\{j_2=j_5\}}
\zeta_{j_1}^{(i_1)}
\zeta_{j_3}^{(i_3)}
\zeta_{j_4}^{(i_4)}
-{\bf 1}_{\{i_3=i_4\ne 0\}}
{\bf 1}_{\{j_3=j_4\}}
\zeta_{j_1}^{(i_1)}
\zeta_{j_2}^{(i_2)}
\zeta_{j_5}^{(i_5)}-
$$
$$
-
{\bf 1}_{\{i_3=i_5\ne 0\}}
{\bf 1}_{\{j_3=j_5\}}
\zeta_{j_1}^{(i_1)}
\zeta_{j_2}^{(i_2)}
\zeta_{j_4}^{(i_4)}
-{\bf 1}_{\{i_4=i_5\ne 0\}}
{\bf 1}_{\{j_4=j_5\}}
\zeta_{j_1}^{(i_1)}
\zeta_{j_2}^{(i_2)}
\zeta_{j_3}^{(i_3)}+
$$
$$
+
{\bf 1}_{\{i_1=i_2\ne 0\}}
{\bf 1}_{\{j_1=j_2\}}
{\bf 1}_{\{i_3=i_4\ne 0\}}
{\bf 1}_{\{j_3=j_4\}}\zeta_{j_5}^{(i_5)}+
{\bf 1}_{\{i_1=i_2\ne 0\}}
{\bf 1}_{\{j_1=j_2\}}
{\bf 1}_{\{i_3=i_5\ne 0\}}
{\bf 1}_{\{j_3=j_5\}}\zeta_{j_4}^{(i_4)}+
$$
$$
+
{\bf 1}_{\{i_1=i_2\ne 0\}}
{\bf 1}_{\{j_1=j_2\}}
{\bf 1}_{\{i_4=i_5\ne 0\}}
{\bf 1}_{\{j_4=j_5\}}\zeta_{j_3}^{(i_3)}+
{\bf 1}_{\{i_1=i_3\ne 0\}}
{\bf 1}_{\{j_1=j_3\}}
{\bf 1}_{\{i_2=i_4\ne 0\}}
{\bf 1}_{\{j_2=j_4\}}\zeta_{j_5}^{(i_5)}+
$$
$$
+
{\bf 1}_{\{i_1=i_3\ne 0\}}
{\bf 1}_{\{j_1=j_3\}}
{\bf 1}_{\{i_2=i_5\ne 0\}}
{\bf 1}_{\{j_2=j_5\}}\zeta_{j_4}^{(i_4)}+
{\bf 1}_{\{i_1=i_3\ne 0\}}
{\bf 1}_{\{j_1=j_3\}}
{\bf 1}_{\{i_4=i_5\ne 0\}}
{\bf 1}_{\{j_4=j_5\}}\zeta_{j_2}^{(i_2)}+
$$
$$
+
{\bf 1}_{\{i_1=i_4\ne 0\}}
{\bf 1}_{\{j_1=j_4\}}
{\bf 1}_{\{i_2=i_3\ne 0\}}
{\bf 1}_{\{j_2=j_3\}}\zeta_{j_5}^{(i_5)}+
{\bf 1}_{\{i_1=i_4\ne 0\}}
{\bf 1}_{\{j_1=j_4\}}
{\bf 1}_{\{i_2=i_5\ne 0\}}
{\bf 1}_{\{j_2=j_5\}}\zeta_{j_3}^{(i_3)}+
$$
$$
+
{\bf 1}_{\{i_1=i_4\ne 0\}}
{\bf 1}_{\{j_1=j_4\}}
{\bf 1}_{\{i_3=i_5\ne 0\}}
{\bf 1}_{\{j_3=j_5\}}\zeta_{j_2}^{(i_2)}+
{\bf 1}_{\{i_1=i_5\ne 0\}}
{\bf 1}_{\{j_1=j_5\}}
{\bf 1}_{\{i_2=i_3\ne 0\}}
{\bf 1}_{\{j_2=j_3\}}\zeta_{j_4}^{(i_4)}+
$$
$$
+
{\bf 1}_{\{i_1=i_5\ne 0\}}
{\bf 1}_{\{j_1=j_5\}}
{\bf 1}_{\{i_2=i_4\ne 0\}}
{\bf 1}_{\{j_2=j_4\}}\zeta_{j_3}^{(i_3)}+
{\bf 1}_{\{i_1=i_5\ne 0\}}
{\bf 1}_{\{j_1=j_5\}}
{\bf 1}_{\{i_3=i_4\ne 0\}}
{\bf 1}_{\{j_3=j_4\}}\zeta_{j_2}^{(i_2)}+
$$
$$
+
{\bf 1}_{\{i_2=i_3\ne 0\}}
{\bf 1}_{\{j_2=j_3\}}
{\bf 1}_{\{i_4=i_5\ne 0\}}
{\bf 1}_{\{j_4=j_5\}}\zeta_{j_1}^{(i_1)}+
{\bf 1}_{\{i_2=i_4\ne 0\}}
{\bf 1}_{\{j_2=j_4\}}
{\bf 1}_{\{i_3=i_5\ne 0\}}
{\bf 1}_{\{j_3=j_5\}}\zeta_{j_1}^{(i_1)}+
$$
\begin{equation}
\label{a5}
+\Biggl.
{\bf 1}_{\{i_2=i_5\ne 0\}}
{\bf 1}_{\{j_2=j_5\}}
{\bf 1}_{\{i_3=i_4\ne 0\}}
{\bf 1}_{\{j_3=j_4\}}\zeta_{j_1}^{(i_1)}\Biggr),
\end{equation}

\vspace{3mm}
\noindent
where ${\bf 1}_A$ is the indicator of the set $A$.

Consider a generalization of the formulas 
(\ref{a1})--(\ref{a5}) for the case of arbitrary multiplicity
$k$ ($k\in {\bf N})$ of the iterated It\^{o} stochastic integral (\ref{ito}).

In order to do this, let us
consider the unordered
set $\{1, 2, \ldots, k\}$ 
and separate it into two parts:
the first part consists of $r$ unordered 
pairs (sequence order of these pairs is also unimportant) and the 
second one consists of the 
remaining $k-2r$ numbers.
So, we have
\begin{equation}
\label{leto5007}
(\{
\underbrace{\{g_1, g_2\}, \ldots, 
\{g_{2r-1}, g_{2r}\}}_{\small{\hbox{part 1}}}
\},
\{\underbrace{q_1, \ldots, q_{k-2r}}_{\small{\hbox{part 2}}}
\}),
\end{equation}
where 
$\{g_1, g_2, \ldots, 
g_{2r-1}, g_{2r}, q_1, \ldots, q_{k-2r}\}=\{1, 2, \ldots, k\},$
braces   
mean an unordered 
set, and pa\-ren\-the\-ses mean an ordered set.

We will say that (\ref{leto5007}) is a partition 
and consider the sum with respect to all possible
partitions
\begin{equation}
\label{leto5008}
\sum_{\stackrel{(\{\{g_1, g_2\}, \ldots, 
\{g_{2r-1}, g_{2r}\}\}, \{q_1, \ldots, q_{k-2r}\})}
{{}_{\{g_1, g_2, \ldots, 
g_{2r-1}, g_{2r}, q_1, \ldots, q_{k-2r}\}=\{1, 2, \ldots, k\}}}}
a_{g_1 g_2, \ldots, 
g_{2r-1} g_{2r}, q_1 \ldots q_{k-2r}},
\end{equation}

\noindent
where $a_{g_1 g_2, \ldots, g_{2r-1} g_{2r}, q_1 \ldots q_{k-2r}}\in {\bf R}$.

Below there are several examples of sums in the form (\ref{leto5008})
$$
\sum_{\stackrel{(\{g_1, g_2\})}{{}_{\{g_1, g_2\}=\{1, 2\}}}}
a_{g_1 g_2}=a_{12},
$$

\vspace{-2mm}
$$
\sum_{\stackrel{(\{\{g_1, g_2\}, \{g_3, g_4\}\})}
{{}_{\{g_1, g_2, g_3, g_4\}=\{1, 2, 3, 4\}}}}
a_{g_1 g_2 g_3 g_4}=a_{12,34} + a_{13,24} + a_{23,14},
$$

\vspace{-2mm}
$$
\sum_{\stackrel{(\{g_1, g_2\}, \{q_1, q_{2}\})}
{{}_{\{g_1, g_2, q_1, q_{2}\}=\{1, 2, 3, 4\}}}}
a_{g_1 g_2, q_1 q_{2}}=
a_{12,34}+a_{13,24}+a_{14,23}
+a_{23,14}+a_{24,13}+a_{34,12},
$$
$$
\sum_{\stackrel{(\{g_1, g_2\}, \{q_1, q_{2}, q_3\})}
{{}_{\{g_1, g_2, q_1, q_{2}, q_3\}=\{1, 2, 3, 4, 5\}}}}
a_{g_1 g_2, q_1 q_{2}q_3}
=
$$

\vspace{-2mm}
$$
=a_{12,345}+a_{13,245}+a_{14,235}
+a_{15,234}+a_{23,145}+a_{24,135}+
$$
$$
+a_{25,134}+a_{34,125}+a_{35,124}+a_{45,123},
$$

$$
\sum_{\stackrel{(\{\{g_1, g_2\}, \{g_3, g_{4}\}\}, \{q_1\})}
{{}_{\{g_1, g_2, g_3, g_{4}, q_1\}=\{1, 2, 3, 4, 5\}}}}
a_{g_1 g_2, g_3 g_{4},q_1}
=
$$

\vspace{-5mm}
$$
=
a_{12,34,5}+a_{13,24,5}+a_{14,23,5}+
a_{12,35,4}+a_{13,25,4}+a_{15,23,4}+
a_{12,54,3}+a_{15,24,3}+
$$
$$
+a_{14,25,3}+a_{15,34,2}+a_{13,54,2}+a_{14,53,2}+
a_{52,34,1}+a_{53,24,1}+a_{54,23,1}.
$$

\vspace{5mm}

Now we can formulate Theorem 1
(see (\ref{tyyy})) 
in another form.

{\bf Theorem 2}\ \cite{2009} (2009) (also see 
\cite{12aa-afterxxx}--\cite{arxiv-1}).\
{\it Suppose that
$\psi_1(\tau),$ $\ldots,\psi_k(\tau)$ are continuous 
nonrandom functions on 
$[t, T]$ and
$\{\phi_j(x)\}_{j=0}^{\infty}$ is an arbitrary complete orthonormal system  
of continuous or piecewise continuous functions in the space $L_2([t,T]).$ 
Then the following expansion
$$
J[\psi^{(k)}]_{T,t}^{(i_1 \ldots i_k)}=
\hbox{\vtop{\offinterlineskip\halign{
\hfil#\hfil\cr
{\rm l.i.m.}\cr
$\stackrel{}{{}_{p_1,\ldots,p_k\to \infty}}$\cr
}} }
\sum\limits_{j_1=0}^{p_1}\ldots
\sum\limits_{j_k=0}^{p_k}
C_{j_k\ldots j_1}\Biggl(
\prod_{l=1}^k\zeta_{j_l}^{(i_l)}+\sum\limits_{r=1}^{[k/2]}
(-1)^r \times \Biggr.
$$

\begin{equation}
\label{leto6000}
\times
\sum_{\stackrel{(\{\{g_1, g_2\}, \ldots, 
\{g_{2r-1}, g_{2r}\}\}, \{q_1, \ldots, q_{k-2r}\})}
{{}_{\{g_1, g_2, \ldots, 
g_{2r-1}, g_{2r}, q_1, \ldots, q_{k-2r}\}=\{1,2, \ldots, k\}}}}
\prod\limits_{s=1}^r
{\bf 1}_{\{i_{g_{{}_{2s-1}}}=~i_{g_{{}_{2s}}}\ne 0\}}
\Biggl.{\bf 1}_{\{j_{g_{{}_{2s-1}}}=~j_{g_{{}_{2s}}}\}}
\prod_{l=1}^{k-2r}\zeta_{j_{q_l}}^{(i_{q_l})}\Biggr)
\end{equation}

\vspace{1mm}
\noindent 
that converges in the mean-square sense is valid$,$ where $i_1,\ldots,i_k=0,1,\ldots,m,$
$[x]$ is an integer part of a real number $x,$
$\prod\limits_{\emptyset}
\stackrel{\sf def}{=}1,$ 
$\sum\limits_{\emptyset}
\stackrel{\sf def}{=}0;$ 
another notations are the same as in Theorem~{\rm 1.}}

Further in this article, we will consider a generalization
of the expansion (\ref{leto6000}) to the case
of an arbitrary complete orthonormal 
systems of functions in the space $L_2([t,T])$ and
$\psi_1(\tau),\ldots,\psi_k(\tau)\in L_2([t, T]).$
Moreover, we will consider a modification 
of (\ref{leto6000}) based on the Hermite polynomials.

It should be  noted that there is a work \cite{Rybakov1000} in which
an expansion similar to (\ref{new9999x}) was obtained (see Sect.~4 for details).
A comparison of our results with the results from \cite{Rybakov1000}
and with other publications will be given in Sect.~4.

\section{Preliminary Results}

\subsection{Expansion of Ite\-ra\-ted It\^{o} Sto\-chas\-tic Integrals 
based on 
Generalized Multiple Fo\-u\-ri\-er Series}

Suppose that $\Phi(t_1,\ldots,t_k)\in L_2([t, T]^k),$ $i_1,\ldots,i_k=0,1,\ldots,m,$\ $d{\bf w}_{\tau}^{(0)}
\stackrel{\sf def}{=}d\tau.$

Let us introduce the following notation for the sum of iterated It\^{o} stochastic integrals
\begin{equation}
\label{chain10100}
J''[\Phi]_{T,t}^{(i_1\ldots i_k)}\stackrel{\sf def}{=}\sum_{(t_1,\ldots,t_k)}
\int\limits_{t}^{T}
\ldots
\int\limits_{t}^{t_2}
\Phi(t_1,\ldots,t_k)d{\bf w}_{t_1}^{(i_1)}
\ldots
d{\bf w}_{t_k}^{(i_k)},
\end{equation}

\noindent
where all permutations $(t_1,\ldots,t_k)$ when summing are 
performed only in the values
$d{\bf w}_{t_1}^{(i_1)}
\ldots $
$d{\bf w}_{t_k}^{(i_k)}.$ At the same time the indices near 
upper 
limits of integration in the iterated stochastic integrals are changed 
correspondently and if $t_r$ swapped with $t_q$ in the  
permutation $(t_1,\ldots,t_k)$, then $i_r$ swapped with $i_q$ in 
the permutation $(i_1,\ldots,i_k).$ 
In addition, 
$$
\int\limits_{t}^{T}
\ldots
\int\limits_{t}^{t_2}
\Phi(t_1,\ldots,t_k)d{\bf w}_{t_1}^{(i_1)}
\ldots
d{\bf w}_{t_k}^{(i_k)}
$$
is the iterated It\^{o} stochastic integral.

Let us give an exumple of the sum (\ref{chain10100}) for $k=3$
$$
J''[\Phi]_{T,t}^{(i_1 i_2 i_3)}\stackrel{\sf def}{=}\sum_{(t_1,t_2,t_3)}
\int\limits_{t}^{T}
\int\limits_{t}^{t_3}
\int\limits_{t}^{t_2}
\Phi(t_1,t_2,t_3)d{\bf w}_{t_1}^{(i_1)}
d{\bf w}_{t_2}^{(i_2)}
d{\bf w}_{t_3}^{(i_3)}=
$$
$$
=\int\limits_{t}^{T}
\int\limits_{t}^{t_3}
\int\limits_{t}^{t_2}
\Phi(t_1,t_2,t_3)d{\bf w}_{t_1}^{(i_1)}
d{\bf w}_{t_2}^{(i_2)}
d{\bf w}_{t_3}^{(i_3)}+
\int\limits_{t}^{T}
\int\limits_{t}^{t_2}
\int\limits_{t}^{t_3}
\Phi(t_1,t_2,t_3)d{\bf w}_{t_1}^{(i_1)}
d{\bf w}_{t_3}^{(i_3)}
d{\bf w}_{t_2}^{(i_2)}+
$$
$$
+\int\limits_{t}^{T}
\int\limits_{t}^{t_3}
\int\limits_{t}^{t_1}
\Phi(t_1,t_2,t_3)d{\bf w}_{t_2}^{(i_2)}
d{\bf w}_{t_1}^{(i_1)}
d{\bf w}_{t_3}^{(i_3)}+
\int\limits_{t}^{T}
\int\limits_{t}^{t_1}
\int\limits_{t}^{t_3}
\Phi(t_1,t_2,t_3)d{\bf w}_{t_2}^{(i_2)}
d{\bf w}_{t_3}^{(i_3)}
d{\bf w}_{t_1}^{(i_1)}+
$$
$$
+\int\limits_{t}^{T}
\int\limits_{t}^{t_2}
\int\limits_{t}^{t_1}
\Phi(t_1,t_2,t_3)d{\bf w}_{t_3}^{(i_3)}
d{\bf w}_{t_1}^{(i_1)}
d{\bf w}_{t_2}^{(i_2)}+
\int\limits_{t}^{T}
\int\limits_{t}^{t_1}
\int\limits_{t}^{t_2}
\Phi(t_1,t_2,t_3)d{\bf w}_{t_3}^{(i_3)}
d{\bf w}_{t_2}^{(i_2)}
d{\bf w}_{t_1}^{(i_1)}.
$$

{\bf Theorem~3}\ \cite{2018a}, \cite{arxiv-1}.\ 
{\it Suppose that 
$\psi_1(\tau),\ldots,\psi_k(\tau)\in L_2([t, T])$ and
$\{\phi_j(x)\}_{j=0}^{\infty}$ is an arbitrary complete orthonormal system  
of functions in the space $L_2([t,T]).$
Then the following expansion
$$
J[\psi^{(k)}]_{T,t}^{(i_1\ldots i_k)} =
\hbox{\vtop{\offinterlineskip\halign{
\hfil#\hfil\cr
{\rm l.i.m.}\cr
$\stackrel{}{{}_{p_1,\ldots,p_k\to \infty}}$\cr
}} }\sum_{j_1=0}^{p_1}\ldots\sum_{j_k=0}^{p_k}
C_{j_k\ldots j_1}J''[\phi_{j_1}\ldots \phi_{j_k}]_{T,t}^{(i_1\ldots i_k)}
$$

\noindent
converging in the mean-square sense is valid, where $J[\psi^{(k)}]_{T,t}^{(i_1\ldots i_k)}$
is the iterated It\^{o} stochastic integral {\rm (\ref{ito}),}
$J''[\phi_{j_1}\ldots \phi_{j_k}]_{T,t}^{(i_1\ldots i_k)}$ is defined by
{\rm (\ref{chain10100})} or has the form
$$
J''[\phi_{j_1}\ldots \phi_{j_k}]_{T,t}^{(i_1\ldots i_k)}=
\sum\limits_{(j_1,\ldots,j_k)}
\int\limits_t^T \phi_{j_k}(t_k)
\ldots
\int\limits_t^{t_{2}}\phi_{j_{1}}(t_{1})
d{\bf w}_{t_1}^{(i_1)}\ldots
d{\bf w}_{t_k}^{(i_k)},
$$

\noindent
where 
$$
\sum\limits_{(j_1,\ldots,j_k)}
$$ 

\noindent
means the sum with respect to all
possible permutations 
$(j_1,\ldots,j_k).$ At the same time if 
$j_r$ swapped with $j_q$ in the permutation $(j_1,\ldots,j_k)$,
then $i_r$ swapped with $i_q$ in the permutation $(i_1,\ldots,i_k).$
Another notations are the same as in Theorems~{\rm 1} and {\rm 2}.
}

{\bf Proof.}\ Using (\ref{chain10100}), we have
\begin{equation}
\label{wi1001}
J[\psi^{(k)}]_{T,t}^{(i_1\ldots i_k)}=\int\limits_t^T\psi_k(t_k) \ldots \int\limits_t^{t_{2}}
\psi_1(t_1) d{\bf w}_{t_1}^{(i_1)}\ldots
d{\bf w}_{t_k}^{(i_k)}=J''[K]_{T,t}^{(i_1\ldots i_k)}\ \ \ \hbox{w.\ p.\ 1},
\end{equation}
where 
$K=K(t_1,\ldots,t_k)$ is defined by (\ref{ppp}).

Applying the linearity property of the It\^{o} stochastic integral
and (\ref{wi1001}), we obtain w.~p.~1

\vspace{-2mm}
$$
J[\psi^{(k)}]_{T,t}^{(i_1\ldots i_k)}=J''[K]_{T,t}^{(i_1\ldots i_k)}=
J''[K_{p_1\ldots p_k}]_{T,t}^{(i_1\ldots i_k)}+J''[K-K_{p_1\ldots p_k}]_{T,t}^{(i_1\ldots i_k)}=
$$

\vspace{-2mm}
\begin{equation}
\label{chain102}
=\sum_{j_1=0}^{p_1}\ldots
\sum_{j_k=0}^{p_k}
C_{j_k\ldots j_1}
J''[\phi_{j_1}\ldots \phi_{j_k}]_{T,t}^{(i_1\ldots i_k)}+
J''[R_{p_1\ldots p_k}]_{T,t}^{(i_1\ldots i_k)},
\end{equation}

\vspace{1mm}
\noindent
where

\vspace{-3.5mm}
$$
R_{p_1\ldots p_k}(t_1,\ldots,t_k)
=K(t_1,\ldots,t_k)-K_{p_1\ldots p_k}(t_1,\ldots,t_k),
$$

\vspace{4mm}
\noindent
$K(t_1,\ldots,t_k)$ and $K_{p_1\ldots p_k}(t_1,\ldots,t_k)$ are defined by (\ref{ppp}) and
(\ref{chain30001x}), respectively;
the Fourier coefficient $C_{j_k\ldots j_1}$ has the form (\ref{ppppa}).

Note that (see (\ref{chain10100}))

\vspace{-2.5mm}
$$
J''[R_{p_1\ldots p_k}]_{T,t}^{(i_1\ldots i_k)}=
$$

\vspace{-4.5mm}
$$
=
\sum_{(t_1,\ldots,t_k)}
\int\limits_{t}^{T}
\ldots
\int\limits_{t}^{t_2}
\Biggl(K(t_1,\ldots,t_k)-\Biggr.
\Biggl.
\sum_{j_1=0}^{p_1}\ldots
\sum_{j_k=0}^{p_k}
C_{j_k\ldots j_1}
\prod_{l=1}^k\phi_{j_l}(t_l)\Biggr)\times
$$

\vspace{-1.5mm}
$$
\times d{\bf w}_{t_1}^{(i_1)}
\ldots
d{\bf w}_{t_k}^{(i_k)},
$$

\vspace{3mm}
\noindent
where notations are the same as in (\ref{chain10100}).

According to  
the standard moment properties of the It\^{o}
stochastic integral \cite{1} and the properties of the Lebesgue
integral, we get the following estimate

$$
{\sf M}\left\{\left(J''[R_{p_1\ldots p_k}]_{T,t}^{(i_1\ldots i_k)}\right)^2\right\}
\le 
$$

\vspace{-3mm}
$$
\le C_k
\sum_{(t_1,\ldots,t_k)}
\int\limits_{t}^{T}
\ldots
\int\limits_{t}^{t_2}
\left(K(t_1,\ldots,t_k)-
\sum_{j_1=0}^{p_1}\ldots
\sum_{j_k=0}^{p_k}
C_{j_k\ldots j_1}
\prod_{l=1}^k\phi_{j_l}(t_l)\right)^2\times
$$

\vspace{-1mm}
\begin{equation}
\label{new321}
\times
dt_1
\ldots
dt_k=
\end{equation}

\vspace{-5mm}
$$
=C_k\int\limits_{[t,T]^k}
\left(K(t_1,\ldots,t_k)-
\sum_{j_1=0}^{p_1}\ldots
\sum_{j_k=0}^{p_k}
C_{j_k\ldots j_1}
\prod_{l=1}^k\phi_{j_l}(t_l)\right)^2
dt_1
\ldots
dt_k=
$$

\vspace{1mm}
\begin{equation}
\label{chain7771}
=C_k\biggl\Vert K-K_{p_1\ldots p_k}\biggr\Vert_{L_2([t,T]^k)}^2,
\end{equation}

\vspace{4mm}
\noindent
where constant $C_k$ 
depends only
on the multiplicity $k$ of the
iterated  It\^{o} stochastic integral
$J[\psi^{(k)}]_{T,t}^{(i_1\ldots i_k)},$ and permutations $(t_1,\ldots,t_k)$ when summing
in (\ref{new321}) are performed 
in the expression $dt_1\ldots dt_k$. At the same time the indices near upper 
limits of integration in the
iterated integrals from (\ref{new321})
are changed correspondently.

Combining (\ref{sos1z}) and 
(\ref{chain7771}), we get
\begin{equation}
\label{new3211}
\lim\limits_{p_1,\ldots,p_k\to\infty}
{\sf M}\left\{\left(J''[R_{p_1\ldots p_k}]_{T,t}^{(i_1\ldots i_k)}\right)^2\right\}=0.
\end{equation}

\vspace{1mm}

From (\ref{chain102}) and (\ref{new3211}) we obtain the following expansion
for the iterated It\^{o} stochastic integral (\ref{ito})
\begin{equation}
\label{new3212}
J[\psi^{(k)}]_{T,t}^{(i_1\ldots i_k)} =
\hbox{\vtop{\offinterlineskip\halign{
\hfil#\hfil\cr
{\rm l.i.m.}\cr
$\stackrel{}{{}_{p_1,\ldots,p_k\to \infty}}$\cr
}} }\sum_{j_1=0}^{p_1}\ldots\sum_{j_k=0}^{p_k}
C_{j_k\ldots j_1}J''[\phi_{j_1}\ldots \phi_{j_k}]_{T,t}^{(i_1\ldots i_k)},
\end{equation}

\vspace{1mm}
\noindent
where $J''[\phi_{j_1}\ldots \phi_{j_k}]_{T,t}^{(i_1\ldots i_k)}$ is defined by
(\ref{chain10100}).

It is easy to see that 
$J''[\phi_{j_1}\ldots \phi_{j_k}]_{T,t}^{(i_1\ldots i_k)}$ can be written in the 
form
\begin{equation}
\label{new3213}
J''[\phi_{j_1}\ldots \phi_{j_k}]_{T,t}^{(i_1\ldots i_k)}=
\sum\limits_{(j_1,\ldots,j_k)}
\int\limits_t^T \phi_{j_k}(t_k)
\ldots
\int\limits_t^{t_{2}}\phi_{j_{1}}(t_{1})
d{\bf w}_{t_1}^{(i_1)}\ldots
d{\bf w}_{t_k}^{(i_k)},
\end{equation}

\noindent
where 
$$
\sum\limits_{(j_1,\ldots,j_k)}
$$ 

\noindent
means the sum with respect to all
possible permutations 
$(j_1,\ldots,j_k).$ At the same time if 
$j_r$ swapped with $j_q$ in the permutation $(j_1,\ldots,j_k)$,
then $i_r$ swapped with $i_q$ in the permutation $(i_1,\ldots,i_k).$

The relations (\ref{new3212}) and (\ref{new3213}) complete the 
proof of Theorem~3. Theorem~3 is proved.

\subsection{Modification and Generalization 
of It\^{o}'s Theorem.
Proof on the Base of the It\^{o} Formula and Without Explicit Use of the Multiple Wiener 
Stochastic Integral}

In this section, we generalize Theorem~3.1 from \cite{ito1951} (1951)
which gives the relaionship between the multiple Wiener stochastic integral
and the Hermite polynomials. Recall that in \cite{ito1951} the case
$i_1=\ldots=i_k\ne 0$ (the case of a scalar standard Wiener process)
has been considered. In the main result of this section,
we will consider the case $i_1,\ldots,i_k=0,1,\ldots,m$
(the case of a multidimensional Wiener process). Moreover, 
our proof diffes from that given in \cite{ito1951} and is based
on the It\^{o} formula. Also, we do not explicitly use
the multiple Wiener stochastic integral in the proof of Theorem~4. 
Although it should be noted that the sum (\ref{chain10100}),
which plays a central role in the proof of Theorem~4, is equal w.~p.~1
to the multiple Wiener stochastic integral with respect
to the components of a multidimensional Wiener process
(see the proof in \cite{2018a}, Sect.~1.11 for details).

Let us introduce some notations.

{\it We will say that the condition {\rm ($\star$)} is fulfilled
for the multi-index $(i_1\ldots i_k)$ $(i_1,\ldots,i_k=0, 1,\ldots, m)$ if
$m_1,\ldots,m_k$ are multiplicities of the elements $i_1,\ldots,i_k,$ respectively$,$ i.e.
$$
\{i_1,\ldots, i_k\}\hspace{-0.4mm}=\hspace{-0.4mm}\{\overbrace{{i_1, \ldots, i_1}}^{m_1},
\overbrace{{i_2, \ldots, i_2}}^{m_2},
\ldots, \overbrace{{i_r, \ldots, i_r}}^{m_r}\},
$$
where $r=1,\ldots, k,$ braces   
mean an unordered 
set, and pa\-ren\-the\-ses mean an ordered set. At that, 
$m_1+\ldots+m_k=k,$\ $m_1,\ldots, m_k=0,1,\ldots,k,$\ 
and all elements with nonzero multiplicities are pairwise different.}

It is not difficult to see that 

\vspace{-1mm}
$$
J''\left[\phi_{j_1}\ldots \phi_{j_k}\right]_{T,t}^{(i_1\ldots i_k)}=
J''\biggl[\underbrace{\phi_{j_{g_1}}
\ldots \phi_{j_{g_{{}_{m_1}}}}}_{m_1}
\underbrace{\phi_{j_{g_{m_1+1}}}
\ldots \phi_{j_{g_{m_1+m_2}}}}_{m_2}\ldots \biggr.
$$

\vspace{1mm}
$$
\biggl.\ldots
\underbrace{\phi_{j_{g_{m_1+m_2+\ldots+m_{k-1}+1}}}\ldots
\phi_{j_{g_{m_1+m_2+\ldots+m_k}}}}_{m_k}\biggr]_{T,t}^
{(\overbrace{{}_{i_1 \ldots i_1}}^{m_1}
\overbrace{{}_{i_2 \ldots i_2}}^{m_2}
\ldots \overbrace{{}_{i_k \ldots i_k}}^{m_k})}
$$

\noindent
w.~p.~1, where we suppose that 
the condition {\rm ($\star$)} is fulfilled
for the multi-index $(i_1 \ldots i_k)$ and
$\{j_{g_1},\ldots,j_{g_{m_1+m_2+\ldots+m_k}}\}=
\{j_{g_1},\ldots,j_{g_{k}}\}=
\{j_1,\ldots,j_k\}.$

Suppose that 
$$
\left\{j_{g_{m_1+m_2+\ldots+m_{l-1}+1}}, \ldots, j_{g_{m_1+m_2+\ldots+m_{l}}}
\right\}=
$$
\begin{equation}
\label{ziko999}
=\biggl\{\underbrace{j_{h_{1,l}}, \ldots, j_{h_{1,l}}}_{n_{1,l}}\ \hspace{-1mm},
\underbrace{j_{h_{2,l}}, \ldots, j_{h_{2,l}}}_{n_{2,l}}\ \hspace{-1mm}, \ldots,
\underbrace{j_{h_{d_l,l}}, \ldots, j_{h_{d_l,l}}}_{n_{d_l,l}}\biggr\},
\end{equation}

\noindent
where
$n_{1,l}+n_{2,l}+\ldots+n_{d_l,l}=m_l;$\ \ $n_{1,l}, n_{2,l}, \ldots, n_{d_l,l}=1,\ldots, m_l;$\ \ 
$d_l=1,\ldots,m_l;$\ \ $l=1,\ldots,k.$ Note that the numbers $m_1,\ldots,m_k,$\ $g_1,\ldots,g_k$
depend on $(i_1,\ldots ,i_k)$ and the numbers
$n_{1,l},\ldots,n_{d_l,l},$\ $h_{1,l},\ldots,h_{d_l,l},$\ $d_l$
depend on $\{j_1,\ldots,j_k\}$. Moreover, 
$\left\{j_{g_1},\ldots,j_{g_k}\right\}
=\{j_1,\ldots,j_k\}.$

Let $H_n(x)$ be the Hermite polynomial of degree $n$

\vspace{-2mm}
$$
H_n(x)=(-1)^n e^{x^2/2} \frac{d^n}{dx^n}\left(e^{-x^2/2}\right)
$$

\vspace{2mm}
\noindent
or
\begin{equation}
\label{ziko500}
H_n(x)=n!\sum\limits_{m=0}^{[n/2]}\frac{(-1)^m x^{n-2m}}{m!(n-2m)! 2^m}\ \ \ (n\in{\bf N}).
\end{equation}

\vspace{2mm}

For example,
$$
H_0(x)=1,
$$
$$
H_1(x)=x,
$$
$$
H_2(x)=x^2-1,
$$
$$
H_3(x)=x^3-3x,
$$
$$
H_4(x)=x^4-6x^2 + 3,
$$
$$
H_5(x)=x^5-10x^3 + 15x.
$$

\vspace{1mm}

Let us formulate the following modification and generalization 
of Theorem~3.1 from \cite{ito1951} for the case $i_1,\ldots,i_k=0, 1,\ldots,m$.

{\bf Theorem~4}\ \cite{2018a}, \cite{arxiv-1}.\ {\it Suppose that
the condition {\rm ($\star$)} is fulfilled
for the multi-index $(i_1 \ldots i_k)$ 
and the condition {\rm (\ref{ziko999})} is also 
fulfilled.
Furthermore$,$ let 
$\{\phi_j(x)\}_{j=0}^{\infty}$ is an arbitrary complete orthonormal system  
of functions in the space $L_2([t,T]).$
Then 

\vspace{-9mm}
$$
J''[\phi_{j_1}\ldots \phi_{j_k}]_{T,t}^{(i_1\ldots i_k)}=
$$

\vspace{-2mm}
$$
=\prod_{l=1}^k\left({\bf 1}_{\{m_l=0\}}+{\bf 1}_{\{m_l>0\}}\left\{
\begin{matrix}
H_{n_{1,l}}\left(\zeta_{j_{h_{1,l}}}^{(i_l)}\right)\ldots 
H_{n_{d_l,l}}\left(\zeta_{j_{h_{d_l,l}}}^{(i_l)}\right),\ 
&\hbox{\rm if}\ \ \ 
i_l\ne 0\cr\cr
\left(\zeta_{j_{h_{1,l}}}^{(0)}\right)^{n_{1,l}}\ldots
\left(\zeta_{j_{h_{d_l,l}}}^{(0)}\right)^{n_{d_l,l}},\  &\hbox{\rm if}\ \ \ 
i_l=0
\end{matrix}\right.\ \right)
$$

\vspace{3mm}
\noindent
w.~p.~{\rm 1,} where 
$H_n(x)$ is the Hermite polynomial {\rm (\ref{ziko500}),}
${\bf 1}_A$ is the indicator of the set $A,$
$i_1,\ldots,i_k=0, 1,\ldots,m;$\ \
$n_{1,l}+n_{2,l}+\ldots+n_{d_l,l}=m_l;$\ \ $n_{1,l}, n_{2,l}, \ldots, n_{d_l,l}=1,\ldots, m_l;$\ \ 
$d_l=1,\ldots,m_l;$\ \ $l=1,\ldots,k;$\ \ $m_1+\ldots+m_k=k;$\ \  
the numbers $m_1,\ldots,m_k,$\ $g_1,\ldots,g_k$
depend on $(i_1,\ldots,i_k)$ and 
the numbers $n_{1,l},\ldots,n_{d_l,l},$\ $h_{1,l},\ldots,h_{d_l,l},$\ $d_l$
depend on $\{j_1,\ldots,j_k\};$ moreover$,$ $\left\{j_{g_1},\ldots,j_{g_k}\right\}
=\{j_1,\ldots,j_k\};$
$$
\zeta_{j}^{(i)}=
\int\limits_t^T \phi_{j}(\tau) d{\bf w}_{\tau}^{(i)}\ \ \ (i=0,1,\ldots,m;\ \ j=0,1,2,\ldots)
$$
are independent standard Gaussian random variables
for various
$i$ or $j$ {\rm(}in the case when $i\ne 0${\rm)}
and $d{\bf w}_{\tau}^{(0)}=d\tau.$}

{\bf Proof.}\ First, consider the case $i_1=\ldots=i_k=1,\ldots, m$~ and~ $j_1,\ldots,j_k\in $
$\{0\}\cup {\bf N}$.
This case has been considered in \cite{ito1951}, but we give a different proof here.
By induction, we prove the following equality
$$
p! \int\limits_t^T \phi_l(t_p)\ldots \int\limits_t^{t_2}
\phi_l(t_1)d{\bf w}_{t_1}^{(1)}\ldots d{\bf w}_{t_p}^{(1)}\times
$$
$$
\times \sum\limits_{(j_1,\ldots,j_q)}
\int\limits_t^T \phi_{j_q}(t_q)\ldots \int\limits_t^{t_2}
\phi_{j_1}(t_1)d{\bf w}_{t_1}^{(1)}\ldots d{\bf w}_{t_q}^{(1)}=
$$
$$
=\sum\limits_{(j_1,\ldots,j_q, \underbrace{{}_{l, \ldots ,l}}_{p})}
\int\limits_t^T \phi_{j_q}(t_q)\ldots \int\limits_t^{t_2}
\phi_{j_1}(t_1)
\int\limits_t^{t_1} \phi_{l}(t_p')\ldots \int\limits_t^{t_2'}
\phi_{l}(t_1')\times
$$
\begin{equation}
\label{new1010}
\times
d{\bf w}_{t_1'}^{(1)}\ldots d{\bf w}_{t_p'}^{(1)}
d{\bf w}_{t_1}^{(1)}\ldots d{\bf w}_{t_q}^{(1)}
\end{equation}

\vspace{5mm}
\noindent
w.~p.~1, where $p\in{\bf N},$\ $l\ne j_1,\ldots,j_q,$ and
$$
\sum\limits_{(q_1,\ldots, q_n)}
$$

\noindent
means the sum with respect to all possible permutations
$(q_1,\ldots, q_n)$.

Consider the case $p=1.$ Using the It\^{o} formula, we get w.~p.~1 for
$s\in[t, T]$
$$
\int\limits_t^s \phi_l(\tau)
d{\bf w}_{\tau}^{(1)}
\int\limits_t^s \phi_{j_q}(t_q)\ldots \int\limits_t^{t_2}
\phi_{j_1}(t_1)d{\bf w}_{t_1}^{(1)}\ldots d{\bf w}_{t_q}^{(1)}=
$$
$$
=\int\limits_t^s \phi_l(\tau)\phi_{j_q}(\tau)
\int\limits_t^{\tau} \phi_{j_{q-1}}(t_{q-1})\ldots \int\limits_t^{t_2}
\phi_{j_1}(t_1)d{\bf w}_{t_1}^{(1)}\ldots d{\bf w}_{t_{q-1}}^{(1)}d\tau+
$$
$$
+\int\limits_t^s \phi_l(\tau)\
\int\limits_t^{\tau} \phi_{j_q}(t_q)\ldots \int\limits_t^{t_2}
\phi_{j_1}(t_1)d{\bf w}_{t_1}^{(1)}\ldots d{\bf w}_{t_q}^{(1)}d{\bf w}_{\tau}^{(1)}+
$$
\begin{equation}
\label{new1011}
+\int\limits_t^s \phi_{j_q}(\tau)\hspace{-0.5mm}
\left(\int\limits_t^{\tau} \phi_{l}(\theta)
d{\bf w}_{\theta}^{(1)}
\int\limits_t^{\tau} \phi_{j_{q-1}}(t_{q-1})\ldots \int\limits_t^{t_2}
\phi_{j_1}(t_1)d{\bf w}_{t_1}^{(1)}\ldots d{\bf w}_{t_{q-1}}^{(1)}\right)
\hspace{-0.5mm}d{\bf w}_{\tau}^{(1)}\hspace{-0.2mm}.
\end{equation}

Hereinafter in this section always $s\in [t, T].$
Differentiating by the It\^{o} formula the expression in parentheses
on the right-hand side of equality (\ref{new1011}) and combining the result 
of differentiation with (\ref{new1011}), we obtain w.~p.~1

\vspace{-4mm}
$$
J_{(l)s,t} J_{(j_q\ldots j_1)s,t}=
$$

\vspace{-4mm}
$$
=\int\limits_t^s \phi_l(\tau)\phi_{j_q}(\tau)
\int\limits_t^{\tau} \phi_{j_{q-1}}(t_{q-1})\ldots \int\limits_t^{t_2}
\phi_{j_1}(t_1)d{\bf w}_{t_1}^{(1)}\ldots d{\bf w}_{t_{q-1}}^{(1)}d\tau+
$$

\vspace{-3mm}
$$
+
J_{(l j_q\ldots j_1)s,t}+
$$

\vspace{-4mm}
$$
+
\int\limits_t^s \hspace{-0.2mm}\phi_{j_q}(\tau)
\hspace{-0.2mm}
\int\limits_t^{\tau}\hspace{-0.2mm}
\phi_{l}(\theta)\phi_{j_{q-1}}(\theta)\hspace{-0.1mm}\int\limits_t^{\theta}
\hspace{-0.2mm}\phi_{j_{q-2}}(t_{q-2})
\ldots \int\limits_t^{t_2}\hspace{-0.2mm}
\phi_{j_1}(t_1)d{\bf w}_{t_1}^{(1)}\ldots d{\bf w}_{t_{q-2}}^{(1)}d\theta
d{\bf w}_{\tau}^{(1)}+
$$

\vspace{-3mm}
$$
+
J_{(j_q l j_{q-1}\ldots j_1)s,t}+
$$

$$
+\int\limits_t^s \phi_{j_q}(\tau)
\int\limits_t^{\tau} \phi_{j_{q-1}}(\theta)\times
$$
$$
\times\left(\int\limits_t^{\theta} \phi_l(u)\
d{\bf w}_{u}^{(1)}\int\limits_t^{\theta} 
\phi_{j_{q-2}}(t_{q-2})\ldots \int\limits_t^{t_2}
\phi_{j_1}(t_1)d{\bf w}_{t_1}^{(1)}\ldots d{\bf w}_{t_{q-2}}^{(1)}
\right)
d{\bf w}_{\theta}^{(1)}d{\bf w}_{\tau}^{(1)},
$$

\vspace{2mm}
\noindent
where
$$
\int\limits_t^s \phi_{j_q}(t_q)\ldots \int\limits_t^{t_2}
\phi_{j_1}(t_1)d{\bf w}_{t_1}^{(1)}\ldots d{\bf w}_{t_q}^{(1)}
\stackrel{\sf def}{=}J_{(j_q \ldots j_1)s,t}.
$$

\vspace{1mm}

Continuing the process of iterative application of the It\^{o} formula, we have w.~p.~1

\vspace{-5mm}
$$
J_{(l)s,t} J_{(j_q\ldots j_1)s,t}=
$$

\vspace{1mm}
$$
=J_{(l j_q\ldots j_1)s,t}+ J_{(j_q l j_{q-1}\ldots j_1)s,t}+\ldots + J_{(j_q\ldots j_1 l)s,t}+
$$

\vspace{-1mm}
$$
+\int\limits_t^s \phi_l(\tau)\phi_{j_q}(\tau)
\int\limits_t^{\tau} \phi_{j_{q-1}}(t_{q-1})\ldots \int\limits_t^{t_2}
\phi_{j_1}(t_1)d{\bf w}_{t_1}^{(1)}\ldots d{\bf w}_{t_{q-1}}^{(1)}d\tau+\ldots
$$
\begin{equation}
\label{new1040ss}
\ldots +
\int\limits_t^{s} \phi_{j_{q}}(t_{q})\ldots \int\limits_t^{t_3}
\phi_{j_2}(t_2)\int\limits_t^{t_2} \phi_{l}(\tau)\phi_{j_1}(\tau)      
d\tau d{\bf w}_{t_2}^{(1)}\ldots d{\bf w}_{t_{q}}^{(1)}.
\end{equation}

\vspace{3mm}

Summing the equality (\ref{new1040ss}) over permutations $(j_1,\ldots, j_q)$, we get 

\vspace{-1mm}
\begin{equation}
\label{new1025}
\sum\limits_{(j_1,\ldots, j_q)}J_{(l)s,t} J_{(j_q\ldots j_1)s,t}=
\sum\limits_{(j_1,\ldots, j_q,l)}J_{(l j_q\ldots j_1)s,t}+ S(s)
\end{equation}

\vspace{3mm}
\noindent
w.~p.~1, where
$$
S(s)=
$$

\vspace{-6mm}
$$
=\sum\limits_{(j_1,\ldots, j_q)}\left(\int\limits_t^s \phi_l(\tau)\phi_{j_q}(\tau)
\int\limits_t^{\tau} \phi_{j_{q-1}}(t_{q-1})\ldots \int\limits_t^{t_2}
\phi_{j_1}(t_1)d{\bf w}_{t_1}^{(1)}\ldots d{\bf w}_{t_{q-1}}^{(1)}d\tau+\ldots\right.
$$

\vspace{-2mm}
\begin{equation}
\label{new1040}
\left.\ldots +
\int\limits_t^{s} \phi_{j_{q}}(t_{q})\ldots \int\limits_t^{t_3}
\phi_{j_2}(t_2)\int\limits_t^{t_2} \phi_{l}(\tau)\phi_{j_1}(\tau)      
d\tau d{\bf w}_{t_2}^{(1)}\ldots d{\bf w}_{t_{q}}^{(1)}\right).
\end{equation}

\vspace{2mm}

Consider 
$$
\int\limits_t^s \phi_l(\tau)\phi_{j_q}(\tau)d\tau
\int\limits_t^{s} \phi_{j_{q-1}}(t_{q-1})\ldots \int\limits_t^{t_2}
\phi_{j_1}(t_1)d{\bf w}_{t_1}^{(1)}\ldots d{\bf w}_{t_{q-1}}^{(1)}.
$$

\vspace{2mm}

Applying the It\^{o} formula, we get w.~p.~1
$$
\int\limits_t^s \phi_l(\tau)\phi_{j_q}(\tau)d\tau
\int\limits_t^{s} \phi_{j_{q-1}}(t_{q-1})\ldots \int\limits_t^{t_2}
\phi_{j_1}(t_1)d{\bf w}_{t_1}^{(1)}\ldots d{\bf w}_{t_{q-1}}^{(1)}=
$$
$$
=\int\limits_t^s \phi_l(\tau)\phi_{j_q}(\tau)
\int\limits_t^{\tau} \phi_{j_{q-1}}(t_{q-1})\ldots \int\limits_t^{t_2}
\phi_{j_1}(t_1)d{\bf w}_{t_1}^{(1)}\ldots d{\bf w}_{t_{q-1}}^{(1)}d\tau+
$$

\vspace{-2mm}
$$
+\int\limits_t^s \phi_{j_{q-1}}(t_{q-1})\times
$$
$$
\times\left(\int\limits_t^{t_{q-1}}\phi_l(\tau)\phi_{j_q}(\tau)d\tau
\int\limits_t^{t_{q-1}} \phi_{j_{q-2}}(t_{q-2})\ldots \int\limits_t^{t_2}
\phi_{j_1}(t_1)d{\bf w}_{t_1}^{(1)}\ldots d{\bf w}_{t_{q-2}}^{(1)}\right)d{\bf w}_{t_{q-1}}^{(1)}.
$$

\vspace{3mm}

By iterative application of the It\^{o} formula (as above), we obtain w.~p.~1

\vspace{-2mm}
$$
\int\limits_t^s \phi_l(\tau)\phi_{j_q}(\tau)d\tau
\int\limits_t^{s} \phi_{j_{q-1}}(t_{q-1})\ldots \int\limits_t^{t_2}
\phi_{j_1}(t_1)d{\bf w}_{t_1}^{(1)}\ldots d{\bf w}_{t_{q-1}}^{(1)}=
$$

\vspace{-3mm}
$$
=\int\limits_t^s \phi_l(\tau)\phi_{j_q}(\tau)
\int\limits_t^{\tau} \phi_{j_{q-1}}(t_{q-1})\ldots \int\limits_t^{t_2}
\phi_{j_1}(t_1)d{\bf w}_{t_1}^{(1)}\ldots d{\bf w}_{t_{q-1}}^{(1)}d\tau+\ldots
$$

\vspace{-3mm}
\begin{equation}
\label{newx1020ss}
\ldots +
\int\limits_t^{s} \phi_{j_{q-1}}(t_{q-1})\ldots \int\limits_t^{t_2}
\phi_{j_1}(t_1)\int\limits_t^{t_1} \phi_{l}(\tau)\phi_{j_q}(\tau)      
d\tau d{\bf w}_{t_1}^{(1)}\ldots d{\bf w}_{t_{q-1}}^{(1)}.
\end{equation}

\vspace{3mm}

Summing the equality (\ref{newx1020ss}) over permutations $(j_1,\ldots, j_q)$, we get 

\vspace{-4mm}
\begin{equation}
\label{new1020}
\sum\limits_{(j_1,\ldots, j_q)}\int\limits_t^s \phi_l(\tau)\phi_{j_q}(\tau)d\tau
\int\limits_t^{s} \phi_{j_{q-1}}(t_{q-1})\ldots \int\limits_t^{t_2}
\phi_{j_1}(t_1)d{\bf w}_{t_1}^{(1)}\ldots d{\bf w}_{t_{q-1}}^{(1)}=S_1(s),
\end{equation}

\vspace{2mm}
\noindent
w.~p.~1, where
$$
S_1(s)=
$$

\vspace{-4mm}
$$
=\sum\limits_{(j_1,\ldots, j_q)}\left(\int\limits_t^s \phi_l(\tau)\phi_{j_q}(\tau)
\int\limits_t^{\tau} \phi_{j_{q-1}}(t_{q-1})\ldots \int\limits_t^{t_2}
\phi_{j_1}(t_1)d{\bf w}_{t_1}^{(1)}\ldots d{\bf w}_{t_{q-1}}^{(1)}d\tau+\ldots\right.
$$
\begin{equation}
\label{newx1020}
\left.\ldots +
\int\limits_t^{s} \phi_{j_{q-1}}(t_{q-1})\ldots \int\limits_t^{t_2}
\phi_{j_1}(t_1)\int\limits_t^{t_1} \phi_{l}(\tau)\phi_{j_q}(\tau)      
d\tau d{\bf w}_{t_1}^{(1)}\ldots d{\bf w}_{t_{q-1}}^{(1)}\right).
\end{equation}

\vspace{3mm}

It is not difficult to see that
\begin{equation}
\label{new1021}
S(s)=S_1(s)\ \ \ \hbox{w.~p.~1.}
\end{equation}

\vspace{2mm}

Moreover, due to the orthogonality of $\{\phi_j(x)\}_{j=0}^{\infty}$
and (\ref{new1020}), (\ref{new1021}), we have
\begin{equation}
\label{new1026}
S(T)=S_1(T)=0\ \ \ \hbox{w.~p.~1.}
\end{equation}

\vspace{3mm}

Thus (see (\ref{new1025}), (\ref{new1026})), the equality (\ref{new1010}) is proved for the case
$p=1.$
Let us assume that the equality (\ref{new1010}) is true for $p=2, 3, \ldots, k-1$, and prove
its validity for $p=k.$

From (\ref{new1025}) for the case $q=k-1,$ $j_1=\ldots=j_{k-1}=l$ we obtain

\vspace{-1mm}
\begin{equation}
\label{new1042}
\left(J_1\right)_{s,t} (k-1)! \left(J_{k-1}\right)_{s,t}=k! \left(J_{k}\right)_{s,t} + S_2(s)
\end{equation}

\vspace{3mm}
\noindent
w.~p.~1, where 

\vspace{-3mm}
$$
S_2(s)=S(s)\biggl|_{j_1=\ldots=j_q=l,\ q=k-1}\biggr.\ \ (k\ge 2)\ \ \ \hbox{and}\ \ \
S_2(s)\stackrel{\sf def}{=}0\ \ (q=k-1,\ k=1),
$$

\vspace{-2mm}
$$
\int\limits_t^s \phi_{l}(t_r)\ldots \int\limits_t^{t_2}
\phi_{l}(t_1)d{\bf w}_{t_1}^{(1)}\ldots d{\bf w}_{t_r}^{(1)}
\stackrel{\sf def}{=}\left(J_{r}\right)_{s,t}\ \ (r\in {\bf N})\ \ \ \hbox{and}\ \ \
\left(J_{0}\right)_{s,t}\stackrel{\sf def}{=}1.
$$

\vspace{3mm}

Taking into account (\ref{new1040}), (\ref{new1020})--(\ref{new1021})
and the orthonormality of $\{\phi_j(x)\}_{j=0}^{\infty}$, we have
\begin{equation}
\label{new1043}
S_2(T)=(k-1)!\left(J_{k-2}\right)_{T,t}.
\end{equation}

\vspace{3mm}

Combining (\ref{new1042}) and (\ref{new1043}), we obtain the following recurrence 
relation

\vspace{-1mm}
\begin{equation}
\label{new1044}
k! \left(J_{k}\right)_{T,t}=\left(J_1\right)_{T,t} (k-1)! \left(J_{k-1}\right)_{T,t}-
(k-1)!\left(J_{k-2}\right)_{T,t}
\end{equation}

\vspace{1mm}
\noindent
w.~p.~1.

Using (\ref{new1044}) and the induction hypothesis, we get w.~p.~1
$$
k! \int\limits_t^T \phi_l(t_k)\ldots \int\limits_t^{t_2}
\phi_l(t_1)d{\bf w}_{t_1}^{(1)}\ldots d{\bf w}_{t_k}^{(1)}\times
$$
$$
\times \sum\limits_{(j_1,\ldots,j_q)}
\int\limits_t^T \phi_{j_q}(t_q)\ldots \int\limits_t^{t_2}
\phi_{j_1}(t_1)d{\bf w}_{t_1}^{(1)}\ldots d{\bf w}_{t_q}^{(1)}=
$$

\vspace{-2mm}
$$
=\int\limits_t^T \phi_l(\tau)\
d{\bf w}_{\tau}^{(1)}\Biggl(
(k-1)!\int\limits_t^T \phi_l(t_{k-1})\ldots \int\limits_t^{t_2}
\phi_l(t_1)d{\bf w}_{t_1}^{(1)}\ldots d{\bf w}_{t_{k-1}}^{(1)}\times\Biggr.
$$
$$
\Biggl.\times \sum\limits_{(j_1,\ldots,j_q)}
\int\limits_t^T \phi_{j_q}(t_q)\ldots \int\limits_t^{t_2}
\phi_{j_1}(t_1)d{\bf w}_{t_1}^{(1)}\ldots d{\bf w}_{t_q}^{(1)}\Biggr)-
$$
$$
-
(k-1)!\int\limits_t^T \phi_l(t_{k-2})\ldots \int\limits_t^{t_2}
\phi_l(t_1)d{\bf w}_{t_1}^{(1)}\ldots d{\bf w}_{t_{k-2}}^{(1)}\times
$$
$$
\times \sum\limits_{(j_1,\ldots,j_q)}
\int\limits_t^T \phi_{j_q}(t_q)\ldots \int\limits_t^{t_2}
\phi_{j_1}(t_1)d{\bf w}_{t_1}^{(1)}\ldots d{\bf w}_{t_q}^{(1)}=
$$
$$
=\int\limits_t^T \phi_l(\tau)\
d{\bf w}_{\tau}^{(1)}
\sum\limits_{(j_1,\ldots,j_q, \underbrace{{}_{l, \ldots ,l}}_{k-1})}
\int\limits_t^T \phi_{j_q}(t_q)\ldots \int\limits_t^{t_2}
\phi_{j_1}(t_1)
\int\limits_t^{t_1} \phi_{l}(t_{k-1}')\ldots \int\limits_t^{t_2'}
\phi_{l}(t_1')\times
$$
$$
\times
d{\bf w}_{t_1'}^{(1)}\ldots d{\bf w}_{t_{k-1}'}^{(1)}
d{\bf w}_{t_1}^{(1)}\ldots d{\bf w}_{t_{q}}^{(1)}-
$$

\vspace{-3mm}
$$
-(k-1)\sum\limits_{(j_1,\ldots,j_q, \underbrace{{}_{l, \ldots ,l}}_{k-2})}
\int\limits_t^T \phi_{j_q}(t_q)\ldots \int\limits_t^{t_2}
\phi_{j_1}(t_1)
\int\limits_t^{t_1} \phi_{l}(t_{k-2}')\ldots \int\limits_t^{t_2'}
\phi_{l}(t_1')\times
$$
\begin{equation}
\label{new1050a}
\times
d{\bf w}_{t_1'}^{(1)}\ldots d{\bf w}_{t_{k-2}'}^{(1)}
d{\bf w}_{t_1}^{(1)}\ldots d{\bf w}_{t_{q}}^{(1)}.
\end{equation}

\vspace{4mm}

Let $\fbox{\it l}$ be the symbol $l$ which does not participate
in the following sum with respect to permutations

\newpage
\noindent
$$
\sum\limits_{(j_1,\ldots,j_q, \underbrace{{}_{l, \ldots ,l}}_{k-1})}.
$$
\vspace{2mm}

Using (\ref{new1025}), we have w.~p.~1

\vspace{-2mm}
$$
\int\limits_t^s \phi_{l}(\tau)\
d{\bf w}_{\tau}^{(1)}
\sum\limits_{(j_1,\ldots,j_q, \underbrace{{}_{l, \ldots ,l}}_{k-1})}
\int\limits_t^s \phi_{j_q}(t_q)\ldots \int\limits_t^{t_2}
\phi_{j_1}(t_1)
\int\limits_t^{t_1} \phi_{l}(t_{k-1}')\ldots \int\limits_t^{t_2'}
\phi_{l}(t_1')\times
$$
$$
\times
d{\bf w}_{t_1'}^{(1)}\ldots d{\bf w}_{t_{k-1}'}^{(1)}
d{\bf w}_{t_1}^{(1)}\ldots d{\bf w}_{t_{q}}^{(1)}=
$$

\vspace{-2mm}
$$
=
\sum\limits_{(j_1,\ldots,j_q, \underbrace{{}_{l, \ldots ,l}}_{k-1})}
\int\limits_t^s \phi_{\small{\fbox{\it l}}}(\tau)\
d{\bf w}_{\tau}^{(1)}\int\limits_t^s \phi_{j_q}(t_q)\ldots \int\limits_t^{t_2}
\phi_{j_1}(t_1)
\int\limits_t^{t_1} \phi_{l}(t_{k-1}')\ldots \int\limits_t^{t_2'}
\phi_{l}(t_1')\times
$$

\vspace{-2mm}
$$
\times
d{\bf w}_{t_1'}^{(1)}\ldots d{\bf w}_{t_{k-1}'}^{(1)}
d{\bf w}_{t_1}^{(1)}\ldots d{\bf w}_{t_{q}}^{(1)}=
$$

$$
=\sum\limits_{(j_1,\ldots,j_q, \underbrace{{}_{l, \ldots ,l}}_{k-1})}
\left(J_{(\small{\fbox{\it l}}j_q\ldots j_1 \underbrace{l \ldots l}_{k-1})s,t}+
J_{(\small{j_q\fbox{\it l}}j_{q-1}\ldots j_1 \underbrace{l \ldots l}_{k-1})s,t}+\ldots\right.
$$
$$
\left.\ldots +J_{(j_q\ldots j_1 \small{\fbox{\it l}}\underbrace{l \ldots l}_{k-1})s,t}+
J_{(j_q\ldots j_1 l\small{\fbox{\it l}}\underbrace{l \ldots l}_{k-2})s,t}+
\ldots + J_{(j_q\ldots j_1 \small{\underbrace{l \ldots l}_{k-1}}\small{\fbox{\it l}})s,t}\right)+
S_3(s)=
$$

\vspace{3mm}
\begin{equation}
\label{new1050b}
=\sum\limits_{(j_1,\ldots,j_q, \underbrace{{}_{l, \ldots ,l}}_{k})}
J_{(j_q\ldots j_1 \small{\underbrace{l \ldots l}_{k}})s,t}
+S_3(s),
\end{equation}

\vspace{3mm}
\noindent
where
$$
S_3(s)=
$$

\vspace{-3mm}
$$
=
\sum\limits_{(j_1,\ldots,j_q, \underbrace{{}_{l, \ldots ,l}}_{k-1})}
\Biggl(\int\limits_t^s \phi_{\small{\fbox{\it l}}}(\tau)\phi_{j_q}(\tau)
\int\limits_t^{\tau} \phi_{j_{q-1}}(t_{q-1})\ldots \int\limits_t^{t_2}
\phi_{j_1}(t_1)
\times\Biggr.
$$
$$
\times\int\limits_t^{t_1} \phi_{l}(t_{k-1}')\ldots \int\limits_t^{t_2'}
\phi_{l}(t_1')
d{\bf w}_{t_1'}^{(1)}\ldots d{\bf w}_{t_{k-1}'}^{(1)}
d{\bf w}_{t_1}^{(1)}\ldots d{\bf w}_{t_{q-1}}^{(1)}d\tau +  \ldots
$$

\vspace{-2mm}
$$
+\ldots 
\int\limits_t^{s} \phi_{j_{q}}(t_{q})\ldots \int\limits_t^{t_3}
\phi_{j_2}(t_2)\int\limits_t^{t_2} \phi_{\small{\fbox{\it l}}}(\tau)\phi_{j_1}(\tau)
\times\Biggr.
$$

\vspace{-2mm}
$$
\times\int\limits_t^{\tau} \phi_{l}(t_{k-1}')\ldots \int\limits_t^{t_2'}
\phi_{l}(t_1')
d{\bf w}_{t_1'}^{(1)}\ldots d{\bf w}_{t_{k-1}'}^{(1)}
d\tau d{\bf w}_{t_2}^{(1)}\ldots d{\bf w}_{t_{q}}^{(1)}+
$$

\vspace{-2mm}
$$
+
\int\limits_t^{s} \phi_{j_{q}}(t_{q})\ldots \int\limits_t^{t_2}
\phi_{j_1}(t_1)\int\limits_t^{t_1} \phi_{\small{\fbox{\it l}}}(\tau)\phi_{l}(\tau)
\times\Biggr.
$$

\vspace{-2mm}
$$
\times\int\limits_t^{\tau} \phi_{l}(t_{k-2}')\ldots \int\limits_t^{t_2'}
\phi_{l}(t_1')
d{\bf w}_{t_1'}^{(1)}\ldots d{\bf w}_{t_{k-2}'}^{(1)}d\tau
d{\bf w}_{t_1}^{(1)}\ldots d{\bf w}_{t_{q}}^{(1)}+ \ldots
$$

\vspace{-2mm}
$$
\ldots +
\int\limits_t^{s} \phi_{j_{q}}(t_{q})\ldots \int\limits_t^{t_2}
\phi_{j_1}(t_1)
\times\Biggr.
$$

\vspace{-2mm}
$$
\Biggl.\times\int\limits_t^{t_1} \phi_{l}(t_{k-1}')\ldots \int\limits_t^{t_3'}
\phi_{l}(t_2')\int\limits_t^{t_2'} \phi_{\small{\fbox{\it l}}}(\tau)\phi_{l}(\tau)
d\tau d{\bf w}_{t_2'}^{(1)}\ldots d{\bf w}_{t_{k-1}'}^{(1)}
d{\bf w}_{t_1}^{(1)}\ldots d{\bf w}_{t_{q}}^{(1)}\Biggr).
$$

\vspace{6mm}

Using (\ref{new1040}), (\ref{new1020})--(\ref{new1021}), we get w.~p.~1

$$
S_3(s)=
$$

\vspace{-2mm}
$$
=\sum\limits_{(j_1,\ldots,j_q, \underbrace{{}_{l, \ldots ,l}}_{k-1})}
\int\limits_t^{s} \phi_{\small{\fbox{\it l}}}(\tau)\phi_{l}(\tau)d\tau
\int\limits_t^{s} \phi_{j_{q}}(t_{q})\ldots \int\limits_t^{t_2}
\phi_{j_1}(t_1)\times
$$
$$
\times
\int\limits_t^{t_1}\phi_{l}(t_{k-2}')\ldots \int\limits_t^{t_2'}
\phi_{l}(t_1')
d{\bf w}_{t_1'}^{(1)}\ldots d{\bf w}_{t_{k-2}'}^{(1)}
d{\bf w}_{t_1}^{(1)}\ldots d{\bf w}_{t_{q}}^{(1)} =
$$
$$
=(k-1)\sum\limits_{(j_1,\ldots,j_q, \underbrace{{}_{l, \ldots ,l}}_{k-2})}
\int\limits_t^{s} \phi_{\small{\fbox{\it l}}}(\tau)\phi_{l}(\tau)d\tau
\int\limits_t^{s} \phi_{j_{q}}(t_{q})\ldots \int\limits_t^{t_2}
\phi_{j_1}(t_1)\times
$$
$$
\times
\int\limits_t^{t_1}\phi_{l}(t_{k-2}')\ldots \int\limits_t^{t_2'}
\phi_{l}(t_1')
d{\bf w}_{t_1'}^{(1)}\ldots d{\bf w}_{t_{k-2}'}^{(1)}
d{\bf w}_{t_1}^{(1)}\ldots d{\bf w}_{t_{q}}^{(1)} +
$$

$$
+\sum\limits_{(j_1,\ldots,j_{q-1}, \underbrace{{}_{l, \ldots ,l}}_{k-1})}
\int\limits_t^{s} \phi_{\small{\fbox{\it l}}}(\tau)\phi_{j_q}(\tau)d\tau
\int\limits_t^{s} \phi_{j_{q-1}}(t_{q-1})\ldots \int\limits_t^{t_2}
\phi_{j_1}(t_1)\times
$$
$$
\times
\int\limits_t^{t_1}\phi_{l}(t_{k-1}')\ldots \int\limits_t^{t_2'}
\phi_{l}(t_1')
d{\bf w}_{t_1'}^{(1)}\ldots d{\bf w}_{t_{k-1}'}^{(1)}
d{\bf w}_{t_1}^{(1)}\ldots d{\bf w}_{t_{q-1}}^{(1)} + 
$$
$$
+\sum\limits_{(j_1,\ldots,j_{q-2}, j_q \underbrace{{}_{l, \ldots ,l}}_{k-1})}
\int\limits_t^{s} \phi_{\small{\fbox{\it l}}}(\tau)\phi_{j_{q-1}}(\tau)d\tau
\int\limits_t^{s} \phi_{j_{q}}(t_{q}) \int\limits_t^{t_q} \phi_{j_{q-2}}(t_{q-2})\ldots \int\limits_t^{t_2}
\phi_{j_1}(t_1)\times
$$

$$
\times
\int\limits_t^{t_1}\phi_{l}(t_{k-1}')\ldots \int\limits_t^{t_2'}
\phi_{l}(t_1')
d{\bf w}_{t_1'}^{(1)}\ldots d{\bf w}_{t_{k-1}'}^{(1)}
d{\bf w}_{t_1}^{(1)}\ldots d{\bf w}_{t_{q-2}}^{(1)}d{\bf w}_{t_{q}}^{(1)} + 
$$

\vspace{-6mm}
$$
\ldots
$$

\vspace{-6mm}
$$
+\sum\limits_{(j_2,\ldots,j_q \underbrace{{}_{l, \ldots ,l}}_{k-1})}
\int\limits_t^{s} \phi_{\small{\fbox{\it l}}}(\tau)\phi_{j_{1}}(\tau)d\tau
\int\limits_t^{s} \phi_{j_{q}}(t_{q}) \ldots \int\limits_t^{t_3}
\phi_{j_2}(t_2)\times
$$
\begin{equation}
\label{new1047}
\times
\int\limits_t^{t_2}\phi_{l}(t_{k-1}')\ldots \int\limits_t^{t_2'}
\phi_{l}(t_1')
d{\bf w}_{t_1'}^{(1)}\ldots d{\bf w}_{t_{k-1}'}^{(1)}
d{\bf w}_{t_2}^{(1)}\ldots d{\bf w}_{t_{q}}^{(1)}.
\end{equation}

\vspace{3mm}

Applying (\ref{new1047}) and the orthonormality of $\{\phi_j(x)\}_{j=0}^{\infty}$, we finally have

\vspace{-2mm}
$$
S_3(T)=(k-1)\sum\limits_{(j_1,\ldots,j_q, \underbrace{{}_{l, \ldots ,l}}_{k-2})}
\int\limits_t^{T} \phi_{j_{q}}(t_{q})\ldots \int\limits_t^{t_2}
\phi_{j_1}(t_1)
\times
$$
\begin{equation}
\label{new1048}
\times
\int\limits_t^{t_1}\phi_{l}(t_{k-2}')\ldots \int\limits_t^{t_2'}
\phi_{l}(t_1')
d{\bf w}_{t_1'}^{(1)}\ldots d{\bf w}_{t_{k-2}'}^{(1)}
d{\bf w}_{t_1}^{(1)}\ldots d{\bf w}_{t_{q}}^{(1)}.
\end{equation}

\vspace{3mm}

Combining (\ref{new1050a}), (\ref{new1050b}), (\ref{new1048}), we obtain w.~p.~1

\vspace{1mm}
$$
k! \int\limits_t^T \phi_l(t_k)\ldots \int\limits_t^{t_2}
\phi_l(t_1)d{\bf w}_{t_1}^{(1)}\ldots d{\bf w}_{t_k}^{(1)}\times
$$
$$
\times \sum\limits_{(j_1,\ldots,j_q)}
\int\limits_t^T \phi_{j_q}(t_q)\ldots \int\limits_t^{t_2}
\phi_{j_1}(t_1)d{\bf w}_{t_1}^{(1)}\ldots d{\bf w}_{t_q}^{(1)}=
$$
$$
=\sum\limits_{(\underbrace{{}_{l, \ldots ,l}}_{k})}
\int\limits_t^T \phi_l(t_k)\ldots \int\limits_t^{t_2}
\phi_l(t_1)d{\bf w}_{t_1}^{(1)}\ldots d{\bf w}_{t_k}^{(1)}\times
$$
$$
\times \sum\limits_{(j_1,\ldots,j_q)}
\int\limits_t^T \phi_{j_q}(t_q)\ldots \int\limits_t^{t_2}
\phi_{j_1}(t_1)d{\bf w}_{t_1}^{(1)}\ldots d{\bf w}_{t_q}^{(1)}=
$$
$$
=\sum\limits_{(j_1,\ldots,j_q, \underbrace{{}_{l, \ldots ,l}}_{k})}
\int\limits_t^T \phi_{j_q}(t_q)\ldots \int\limits_t^{t_2}
\phi_{j_1}(t_1)
\int\limits_t^{t_1} \phi_{l}(t_k')\ldots \int\limits_t^{t_2'}
\phi_{l}(t_1')\times
$$
\begin{equation}
\label{new1060}
\times
d{\bf w}_{t_1'}^{(1)}\ldots d{\bf w}_{t_k'}^{(1)}
d{\bf w}_{t_1}^{(1)}\ldots d{\bf w}_{t_q}^{(1)},
\end{equation}

\vspace{3mm}
\noindent
where $l\ne j_1,\ldots, j_q.$

The equality (\ref{new1010}) is proved. From the other hand, (\ref{new1060}) means that
\begin{equation}
\label{new1061}
J''[\phi_{j_1}\ldots \phi_{j_q}\underbrace{\phi_{l}\ldots \phi_{l}}_{n}]^
{(\hspace{0.5mm}\small{\overbrace{1 \ldots 1}^{q+n}}\hspace{0.5mm})}_{T,t}=
J''[\underbrace{\phi_{l}\ldots \phi_{l}}_{n}]^
{(\hspace{0.5mm}\small{\overbrace{1 \ldots 1}^{n}}\hspace{0.5mm})}_{T,t}
\cdot J''[\phi_{j_1}\ldots \phi_{j_q}]^
{(\hspace{0.5mm}\small{\overbrace{1 \ldots 1}^{q}}\hspace{0.5mm})}_{T,t}
\end{equation}
w.~p.~1, where $n, q=0,1,2\ldots;$\ $l\ne j_1,\ldots, j_q$ and
$$
J''[\phi_{j_1}\ldots \phi_{j_q}]^
{(\hspace{0.5mm}\small{\overbrace{1 \ldots 1}^{q}}\hspace{0.5mm})}_{T,t}\stackrel{\sf def}{=}1
$$
for $q=0.$

Consider 
polynomials 
$H_n(x,y),$ $n=0, 1,\ldots$ defined by \cite{Ch}
\begin{equation}
\label{new1090}
\Biggl.H_n(x,y)=\left(\frac{d^n}{d\alpha^n} 
e^{\alpha x-\alpha^2 y/2}\right)
\Biggr|_{\alpha=0}\ \ \ (H_0(x,y)\stackrel{\sf def}{=}1).
\end{equation}

\vspace{1mm}

It is well known that polynomials $H_n(x,y)$ are connected with 
the Hermite polynomials (\ref{ziko500}) by the formula \cite{Ch}
\begin{equation}
\label{ziko1000}
H_n(x,y)=y^{n/2}
H_n\left(\frac{x}{\sqrt{y}}\right)=
n!\sum\limits_{i=0}^{[n/2]}\frac{(-1)^i x^{n-2i} y^i}{i!(n-2i)! 2^i}.
\end{equation}

For example,
$$
H_1(x,y)
=x,\ \ \ 
H_2(x,y)
=x^2-y,\ \ \
H_3(x,y)
=x^3-3xy,
$$
$$
H_4(x,y)
=x^4-6x^2 y
+3y^2,\ \ \ 
H_5(x,y)=x^5-10x^3 y+15xy^2.
$$

\vspace{2mm}

From (\ref{ziko500}) and (\ref{ziko1000}) we get
\begin{equation}
\label{ziko1001}
H_n(x,1)=H_n(x).
\end{equation}

Note that \cite{Ch} (also see \cite{2018a} (Chapter 6, Sect.~6.6) for details)
$$
\int\limits_t^T \phi_l(t_n)\ldots \int\limits_t^{t_2}
\phi_l(t_1)d{\bf w}_{t_1}^{(1)}\ldots d{\bf w}_{t_n}^{(1)}=
\frac{1}{n!} H_n\left(\int\limits_t^T \phi_l(\tau)d{\bf w}_{\tau}^{(1)},
\int\limits_t^T \phi_l^2(\tau)d\tau\right)=
$$
\begin{equation}
\label{new1100}
=
\frac{1}{n!} H_n\left(\int\limits_t^T \phi_l(\tau)d{\bf w}_{\tau}^{(1)},1\right)=
\frac{1}{n!} H_n\left(\int\limits_t^T \phi_l(\tau)d{\bf w}_{\tau}^{(1)}\right)
\end{equation}

\vspace{2mm}
\noindent
w.~p.~1, where $n\in {\bf N},$ $H_n(x, y)$ is defined by (\ref{new1090})
(also see (\ref{ziko1000})), and 
$H_n(x)$ is the Hermite polynomial (\ref{ziko500}).

From (\ref{new1100}) we have w.~p.~1
$$
J''[\underbrace{\phi_{l}\ldots \phi_{l}}_{n}]^
{(\hspace{0.5mm}\small{\overbrace{1 \ldots 1}^{n}}\hspace{0.5mm})}_{T,t}
=n! \int\limits_t^T \phi_l(t_n)\ldots \int\limits_t^{t_2}
\phi_l(t_1)d{\bf w}_{t_1}^{(1)}\ldots d{\bf w}_{t_n}^{(1)}=
$$
\begin{equation}
\label{new1101}
=n! \frac{1}{n!} H_n\left(\int\limits_t^T \phi_l(\tau)d{\bf w}_{\tau}^{(1)}\right)=
H_n\left(\int\limits_t^T \phi_l(\tau)d{\bf w}_{\tau}^{(1)}\right),
\end{equation}

\vspace{2mm}
\noindent
where $n\in{\bf N}.$

Combining (\ref{new1061}) and (\ref{new1101}), we obtain
\begin{equation}
\label{new1102}
J''[\phi_{j_1}\ldots \phi_{j_q}\underbrace{\phi_{l}\ldots \phi_{l}}_{n}]^
{(\hspace{0.5mm}\small{\overbrace{1 \ldots 1}^{q+n}}\hspace{0.5mm})}_{T,t}=
H_n\left(\int\limits_t^T \phi_l(\tau)d{\bf w}_{\tau}^{(1)}\right)
\cdot J''[\phi_{j_1}\ldots \phi_{j_q}]^
{(\hspace{0.5mm}\small{\overbrace{1 \ldots 1}^{q}}\hspace{0.5mm})}_{T,t}
\end{equation}
w.~p.~1, where $n, q=0,1,2\ldots;$\ $l\ne j_1,\ldots, j_q.$ 

The iterated application of the formula (\ref{new1102})
completes the proof of Theorem~4 for the case 
$i_1=\ldots=i_k=1,\ldots, m$ and $j_1,\ldots,j_k\in \{0\}\cup {\bf N}$.

To prove Theorem~4 for the case $i_1=\ldots=i_k=0, 1,\ldots, m$ and 
$j_1,\ldots,j_k\in \{0\}\cup {\bf N}$, we need to prove
the following formula in addition to the previous proof
$$
p! \int\limits_t^T \phi_l(t_p)\ldots \int\limits_t^{t_2}
\phi_l(t_1)dt_1\ldots dt_p
\sum\limits_{(j_1,\ldots,j_q)}
\int\limits_t^T \phi_{j_q}(t_q)\ldots \int\limits_t^{t_2}
\phi_{j_1}(t_1)dt_1\ldots dt_q=
$$
\begin{equation}
\label{new1200}
=\sum\limits_{(j_1,\ldots,j_q, \underbrace{{}_{l, \ldots ,l}}_{p})}
\int\limits_t^T \phi_{j_q}(t_q)\ldots \int\limits_t^{t_2}
\phi_{j_1}(t_1)
\int\limits_t^{t_1} \phi_{l}(t_p')\ldots \int\limits_t^{t_2'}
\phi_{l}(t_1')
dt_1'\ldots dt_p'
dt_1\ldots dt_q,
\end{equation}

\vspace{3mm}
\noindent
where $p\in{\bf N}$,
$$
\sum\limits_{(j_1,\ldots, j_{d})}
$$

\noindent
means the sum with respect to all possible permutations $(j_1,\ldots,j_{d}).$ 

First, consider the case $p=1.$ We have 
$$
d\left(\int\limits_t^s \phi_l(\theta)d\theta
\int\limits_t^s \phi_{j_q}(t_q)\ldots \int\limits_t^{t_2}
\phi_{j_1}(t_1)dt_1\ldots dt_q\right)=
$$
$$
=\phi_l(s)
\int\limits_t^s \phi_{j_q}(t_q)\ldots \int\limits_t^{t_2}
\phi_{j_1}(t_1)dt_1\ldots dt_q ds+
$$
$$
+\phi_{j_q}(s)\left(\int\limits_t^{s}\phi_{j_{q-1}}(t_{q-1})\ldots \int\limits_t^{t_2}
\phi_{j_1}(t_1)dt_{1}\ldots dt_{q-1}
\cdot
\int\limits_t^s \phi_l(\theta)d\theta\right)ds.
$$

\vspace{3mm}

Then
$$
\int\limits_t^s \phi_l(\theta)d\theta
\int\limits_t^s \phi_{j_q}(t_q)\ldots \int\limits_t^{t_2}
\phi_{j_1}(t_1)dt_1\ldots dt_q=
$$
$$
=I_{(l j_q \ldots j_1)s,t}+
$$
$$
+\int\limits_t^s
\phi_{j_q}(\tau)\left(\int\limits_t^{\tau}\phi_{j_{q-1}}(t_{q-1})\ldots \int\limits_t^{t_2}
\phi_{j_1}(t_1)dt_{1}\ldots dt_{q-1}
\cdot
\int\limits_t^{\tau} \phi_l(\theta)d\theta\right)d\tau,
$$

\vspace{2mm}
\noindent
where
\begin{equation}
\label{new1701}
\int\limits_t^s \phi_{j_r}(t_r)\ldots \int\limits_t^{t_2}
\phi_{j_1}(t_1)dt_1\ldots dt_r
\stackrel{\sf def}{=}I_{(j_r \ldots j_1)s,t}.
\end{equation}

\vspace{2mm}

Continuing this process, we get
\begin{equation}
\label{new1301}
\int\limits_t^s \phi_l(\theta)d\theta\sum\limits_{(j_1,\ldots,j_q)}
I_{(j_q \ldots j_1)s,t}=
\sum\limits_{(j_1,\ldots,j_q, l)}
I_{(l j_q \ldots j_1)s,t}.
\end{equation}

\vspace{2mm}

The equality (\ref{new1200}) is proved for the case $p=1.$
Let us assume that the equality (\ref{new1200}) is true for $p=2, 3, \ldots, k-1$, and prove
its validity for $p=k.$

From (\ref{new1301}) for $j_1=\ldots=j_q=l,$ $q=k-1$ we have

\vspace{-1mm}
\begin{equation}
\label{new1400}
\left(I_1\right)_{s,t} (k-1)! \left(I_{k-1}\right)_{s,t}=k!\left(I_{k}\right)_{s,t},
\end{equation}

\vspace{1mm}
\noindent
where $k\in {\bf N}$ and
$$
\int\limits_t^s \phi_{l}(t_k)\ldots \int\limits_t^{t_2}
\phi_{l}(t_1)dt_1\ldots dt_k
\stackrel{\sf def}{=}\left(I_k\right)_{s,t},\ \ \ \left(I_0\right)_{s,t}\stackrel{\sf def}{=}1.
$$

\vspace{2mm}

Using (\ref{new1400}) and the induction hypothesis, we obtain 

\vspace{-2mm}
$$
k! \left(I_k\right)_{s,t}
\sum\limits_{(j_1,\ldots,j_q)}
I_{(j_q\ldots j_1)s,t}=
\left(I_1\right)_{s,t} (k-1)! \left(I_{k-1}\right)_{s,t}
\sum\limits_{(j_1,\ldots,j_q)}I_{(j_q\ldots j_1)s,t}=
$$

\vspace{-2mm}
\begin{equation}
\label{new1499}
~~~=I_{(l)s,t} 
\sum\limits_{(j_1,\ldots,j_q, \underbrace{{}_{l, \ldots ,l}}_{k-1})}
I_{(j_q \ldots j_1 \underbrace{{}_{l, \ldots ,l}}_{k-1})s,t}=
\sum\limits_{(j_1,\ldots,j_q, \underbrace{{}_{l, \ldots ,l}}_{k-1})}
I_{(\small{\fbox{\it l}})s,t} 
I_{(j_q \ldots j_1 \underbrace{{}_{l, \ldots ,l}}_{k-1})s,t},
\end{equation}

\vspace{3mm}
\noindent
where $I_{(j_r \ldots j_1)s,t}$ is defined by (\ref{new1701})
and $\fbox{\it l}$ is the symbol $l$ which does not participate
in the following sum with respect to permutations
$$
\sum\limits_{(j_1,\ldots,j_q, \underbrace{{}_{l, \ldots ,l}}_{k-1})}.
$$

\vspace{1mm}

By analogy with (\ref{new1301}) we obtain

\vspace{-2mm}
$$
\sum\limits_{(j_1,\ldots,j_q, \underbrace{{}_{l, \ldots ,l}}_{k-1})}
I_{(\small{\fbox{\it l}})s,t}
I_{(j_q \ldots j_1 \underbrace{{}_{l, \ldots ,l}}_{k-1})s,t}=
$$
$$
=\sum\limits_{(j_1,\ldots,j_q, \underbrace{{}_{l, \ldots ,l}}_{k-1})}
\left(I_{(\small{\fbox{\it l}}j_q\ldots j_1 \underbrace{l \ldots l}_{k-1})s,t}+
I_{(\small{j_q\fbox{\it l}}j_{q-1}\ldots j_1 \underbrace{l \ldots l}_{k-1})s,t}
+\ldots\right.
$$
$$
\left.\ldots +I_{(j_q\ldots j_1 \small{\fbox{\it l}}\underbrace{l \ldots l}_{k-1})s,t}
+
I_{(j_q\ldots j_1 l\small{\fbox{\it l}}\underbrace{l \ldots l}_{k-2})s,t}
+
\ldots + I_{(j_q\ldots j_1 \small{\underbrace{l \ldots l}_{k-1}}\small{\fbox{\it l}})s,t}
\right)=
$$

\vspace{3mm}
\begin{equation}
\label{new1500}
=\sum\limits_{(j_1,\ldots,j_q, \underbrace{{}_{l, \ldots ,l}}_{k})}
I_{(j_q\ldots j_1 \small{\underbrace{l \ldots l}_{k}})s,t}.
\end{equation}

\vspace{1mm}

Substituting $s=T$ into (\ref{new1499}), (\ref{new1500})
and combining (\ref{new1499}), (\ref{new1500}), we conlude that the equality (\ref{new1200}) 
is proved for $p=k.$ The equality (\ref{new1200}) 
is proved.

Note that
\begin{equation}
\label{new1505}
n! \int\limits_t^T \phi_l(t_n)\ldots \int\limits_t^{t_2}
\phi_l(t_1)dt_1\ldots dt_n = n! \frac{1}{n!}
\left(\int\limits_t^T \phi_l(\tau)d\tau\right)^n=
\left(\int\limits_t^T \phi_l(\tau)d\tau\right)^n,
\end{equation}

\noindent
where $n\in{\bf N}.$

After substituting (\ref{new1505}) into (\ref{new1200}), we have for $p=n$
\begin{equation}
\label{new1506}
\left(\int\limits_t^T \phi_l(\tau)d\tau\right)^n
\sum\limits_{(j_1,\ldots,j_q)}
J_{(j_q\ldots j_1)T,t}
=\sum\limits_{(j_1,\ldots,j_q, \underbrace{{}_{l, \ldots ,l}}_{n})}
J_{(j_q\ldots j_1 \small{\underbrace{l \ldots l}_{n}})T,t}.
\end{equation}

\vspace{2mm}

The equality (\ref{new1506}) means that
\begin{equation}
\label{new1507b}
J''[\phi_{j_1}\ldots \phi_{j_q}\underbrace{\phi_{l}\ldots \phi_{l}}_{n}]^
{(\hspace{0.5mm}\small{\overbrace{0 \ldots 0}^{q+n}}\hspace{0.5mm})}_{T,t}=
\left(\int\limits_t^T \phi_l(\tau)d\tau\right)^n
\cdot J''[\phi_{j_1}\ldots \phi_{j_q}]^
{(\hspace{0.5mm}\small{\overbrace{0 \ldots 0}^{q}}\hspace{0.5mm})}_{T,t},
\end{equation}
where $n, q=0,1,2\ldots $ 
and $J''[\phi_{j_1}\ldots \phi_{j_q}]^
{(0 \ldots 0)}_{T,t}\stackrel{\sf def}{=}1$
for $q=0.$

The relations (\ref{new1102}) and (\ref{new1507b}) prove 
Theorem~4 for the case $i_1=\ldots=i_k=0, 1,\ldots, m$ and $j_1,\ldots,j_k\in \{0\}\cup {\bf N}$.

\vspace{2mm}

{\bf Remark~1.}\ {\it Note that the equality
{\rm (\ref{new1200})} can be obtained in another way.
Let $D_q=\left\{(t_1,\ldots,t_q)\in [t, T]^q:\
\exists\ i\ne j\ \hbox{such that}\ t_i=t_j\right\}$ be the "diagonal set" of $[t,T]^q$
$(q=2,3,\ldots)$
{\rm \cite{Kuo}}. Since the Lebesgue meashure of the set $D_q$ is equal to zero {\rm \cite{Kuo}}$,$ 
then {\rm(}see {\rm (\ref{chain10100}))}
\begin{equation}
\label{new9000a}
J''[\phi_{j_1}\ldots \phi_{j_q}]^
{(\hspace{0.5mm}\small{\overbrace{0 \ldots 0}^{q}}\hspace{0.5mm})}_{T,t}=
\int\limits_{[t,T]^q}\phi_{j_1}(t_1)\ldots \phi_{j_q}(t_q)dt_1\ldots dt_q.
\end{equation}

From {\rm (\ref{new9000a})} we have
$$
J''[\phi_{l}\ldots \phi_{l}]^
{(\hspace{0.5mm}\small{\overbrace{0 \ldots 0}^{p}}\hspace{0.5mm})}_{T,t}\cdot
J''[\phi_{j_1}\ldots \phi_{j_q}]^
{(\hspace{0.5mm}\small{\overbrace{0 \ldots 0}^{q}}\hspace{0.5mm})}_{T,t}=
$$

\vspace{-6mm}
$$
=\int\limits_{[t,T]^q}\phi_{j_1}(t_1)\ldots \phi_{j_q}(t_q)dt_1\ldots dt_q
\int\limits_{[t,T]^p}\phi_{l}(t_1)\ldots \phi_{l}(t_p)dt_1\ldots dt_p=
$$

\vspace{-1mm}
$$
=\int\limits_{[t,T]^{p+q}}\phi_{j_1}(t_1)\ldots \phi_{j_q}(t_q)
\phi_l(t_1')\ldots \phi_l(t_p')dt_1'\ldots dt_p'dt_1\ldots dt_q=
$$
\begin{equation}
\label{new100000}
=J''[\phi_{j_1}\ldots \phi_{j_q} \phi_{l}\ldots \phi_{l}]^
{(\hspace{0.5mm}\small{\overbrace{0 \ldots 0}^{p+q}}\hspace{0.5mm})}_{T,t}.
\end{equation}

\vspace{2mm}

It is not difficult to see that the equality {\rm (\ref{new100000})} is nothing but
the equality {\rm (\ref{new1200})} written in another form.}

To complete the proof of Theorem~4, we need to consider the case
$i_1,\ldots,i_k=0, 1,\ldots, m$ and $j_1,\ldots,j_k\in \{0\}\cup {\bf N}$.

Obviously, the proof of Theorem~4 will be completed if we prove the following equalities
\newpage
\noindent
$$
\sum\limits_{(j_1,\ldots,j_q)}\int\limits_t^T \phi_{j_q}(t_q)\ldots
\int\limits_t^{t_2}\phi_{j_1}(t_1)
d{\bf w}_{t_1}^{(i_1)}\ldots d{\bf w}_{t_q}^{(i_q)}\times
$$

\vspace{-2mm}
$$
\times \sum\limits_{(j_1',\ldots,j_n')}
\int\limits_t^T 
\phi_{j_n'}(t_n')\ldots \int\limits_t^{t_2'}\phi_{j_1'}(t_1')d{\bf w}_{t_1'}^{(1)}\ldots 
d{\bf w}_{t_n'}^{(1)}=
$$
$$
=\sum\limits_{(j_1,\ldots,j_q,j_1',\ldots,j_n')}
\int\limits_t^T \phi_{j_q}(t_q)\ldots \int\limits_t^{t_2}\phi_{j_1}(t_1)
\int\limits_t^{t_1}\phi_{j_n'}(t_n')\ldots \int\limits_t^{t_2'}
\phi_{j_1'}(t_1')\times
$$

\vspace{1mm}
\begin{equation}
\label{new1600}
\times d{\bf w}_{t_1'}^{(1)}\ldots d{\bf w}_{t_n'}^{(1)}d{\bf w}_{t_1}^{(i_1)}\ldots d{\bf w}_{t_q}^{(i_q)},
\end{equation}

\vspace{3mm}
$$
\sum\limits_{(j_1,\ldots,j_q)}\int\limits_t^T \phi_{j_q}(t_q)\ldots
\int\limits_t^{t_2}\phi_{j_1}(t_1)
d{\bf w}_{t_1}^{(i_1)}\ldots d{\bf w}_{t_q}^{(i_q)}\times
$$

\vspace{-2mm}
$$
\times \sum\limits_{(j_1',\ldots,j_n')}
\int\limits_t^T \phi_{j_n'}(t_n')\ldots \int\limits_t^{t_2'}\phi_{j_1'}(t_1')
d{\bf w}_{t_1'}^{(0)}\ldots d{\bf w}_{t_n'}^{(0)}=
$$

\vspace{-2mm}
$$
=\sum\limits_{(j_1,\ldots,j_q,j_1',\ldots,j_n')}
\int\limits_t^T \phi_{j_q}(t_q)\ldots \int\limits_t^{t_2}\phi_{j_1}(t_1)
\int\limits_t^{t_1}\phi_{j_n'}(t_n')\ldots \int\limits_t^{t_2'}
\phi_{j_1'}(t_1')\times
$$

\vspace{2mm}
\begin{equation}
\label{new1600a}
\times d{\bf w}_{t_1'}^{(0)}\ldots d{\bf w}_{t_n'}^{(0)}d{\bf w}_{t_1}^{(i_1)}\ldots d{\bf w}_{t_q}^{(i_q)}
\end{equation}

\vspace{4mm}
\noindent
w.~p.~1, where $n, q\in {\bf N},$\ $d{\bf w}_{\tau}^{(0)}
=d\tau,$\ $i_1,\ldots,i_q\ne 1$ in (\ref{new1600})
and $i_1,\ldots,i_q\ne 0$ in (\ref{new1600a}),
$$
\sum\limits_{(j_1,\ldots,j_g)}
$$

\vspace{1mm}
\noindent
means the sum with respect to all possible permutations $(j_1,\ldots,j_g)$.
At the same time if $j_r$ swapped with $j_d$ in the permutation $(j_1,\ldots,j_g)$, then
$i_r$ swapped with $i_d$ in the permutation $(i_1,\ldots,i_g).$

The equalities (\ref{new1600}) and (\ref{new1600a}) mean that

\vspace{-1mm}
\begin{equation}
\label{new1601}
J''[\phi_{j_1}\ldots\phi_{j_q}\phi_{j_1'}\ldots\phi_{j_n'}]_{T,t}^{(i_1\ldots i_q 1\ldots 1)}=
J''[\phi_{j_1}\ldots\phi_{j_q}]_{T,t}^{(i_1\ldots i_q)}
\cdot J''[\phi_{j_1'}\ldots\phi_{j_n'}]^{(1\ldots 1)}_{T,t},
\end{equation}

\newpage
\noindent
\begin{equation}
\label{new1601a}
J''[\phi_{j_1}\ldots\phi_{j_q}\phi_{j_1'}\ldots\phi_{j_n'}]_{T,t}^{(i_1\ldots i_q 0\ldots 0)}=
J''[\phi_{j_1}\ldots\phi_{j_q}]_{T,t}^{(i_1\ldots i_q)}
\cdot J''[\phi_{j_1'}\ldots\phi_{j_n'}]^{(0\ldots 0)}_{T,t}
\end{equation}

\vspace{3mm}
\noindent
w.~p.~1, where $i_1,\ldots,i_q\ne 1$ in (\ref{new1601}) and
$i_1,\ldots,i_q\ne 0$ in (\ref{new1601a}).

First, we prove the equality (\ref{new1600}). Consider the case $n=1.$
Using the It\^{o} formula, we get w.~p.~1
$$
\int\limits_t^{s}\phi_{j_1'}(\theta)d{\bf w}_{\theta}^{(1)}
\int\limits_t^s \phi_{j_q}(t_q)\ldots
\int\limits_t^{t_2}\phi_{j_1}(t_1)
d{\bf w}_{t_1}^{(i_1)}\ldots d{\bf w}_{t_q}^{(i_q)}=
$$
$$
=J_{(j_1' j_q\ldots j_1)s,t}^{(1 i_q\ldots i_1)}+
$$

\vspace{-5mm}
$$
+\hspace{-0.6mm}
\int\limits_t^s\hspace{-1mm} \phi_{j_q}(\tau)\hspace{-0.6mm}\left(\int\limits_t^{\tau}
\hspace{-0.2mm}\phi_{j_{q-1}}(t_{q-1})\ldots
\int\limits_t^{t_2}\hspace{-0.2mm}\phi_{j_1}(t_1)
d{\bf w}_{t_1}^{(i_1)}\ldots d{\bf w}_{t_{q-1}}^{(i_{q-1})}
\hspace{-0.6mm}\int\limits_t^{\tau}\hspace{-0.2mm}
\phi_{j_1'}(\theta)d{\bf w}_{\theta}^{(1)}\right)\hspace{-0.6mm}
d{\bf w}_{\tau}^{(i_{q})}\hspace{-0.5mm}=
$$
$$
\ldots 
$$

\vspace{-3mm}
\begin{equation}
\label{new1700}
=J_{(j_1' j_q\ldots j_1)s,t}^{(1 i_q\ldots i_1)}+J_{(j_q j_1' j_{q-1}\ldots j_1)s,t}
^{(i_q 1 i_{q-1}\ldots i_1)}+\ldots +J_{(j_q\ldots j_1 j_1')s,t}^{(i_q\ldots i_1 1)},
\end{equation}

\vspace{5mm}
\noindent
where
\begin{equation}
\label{new100001}
\int\limits_t^s \phi_{j_r}(t_r)\ldots \int\limits_t^{t_2}
\phi_{j_1}(t_1)d{\bf w}_{t_1}^{(i_1)}\ldots d{\bf w}_{t_r}^{(i_r)}
\stackrel{\sf def}{=}J_{(j_r \ldots j_1)s,t}^{(i_r\ldots i_1)},
\end{equation}

\vspace{3mm}
\noindent
$i_1,\ldots,i_r=0,1,\ldots,m.$

From (\ref{new1700}) we obtain

\vspace{-1mm}
$$
\int\limits_t^{s}\phi_{j_1'}(\theta)d{\bf w}_{\theta}^{(1)}
\sum\limits_{(j_1,\ldots,j_q)}\int\limits_t^s \phi_{j_q}(t_q)\ldots
\int\limits_t^{t_2}\phi_{j_1}(t_1)
d{\bf w}_{t_1}^{(i_1)}\ldots d{\bf w}_{t_q}^{(i_q)}=
$$

\vspace{-3mm}
$$
=\sum\limits_{(j_1,\ldots,j_q)}\int\limits_t^{s}\phi_{j_1'}(\theta)d{\bf w}_{\theta}^{(1)}
\int\limits_t^s \phi_{j_q}(t_q)\ldots
\int\limits_t^{t_2}\phi_{j_1}(t_1)
d{\bf w}_{t_1}^{(i_1)}\ldots d{\bf w}_{t_q}^{(i_q)}=
$$

\vspace{3mm}
$$
=
\sum\limits_{(j_1,\ldots,j_q)}\left(
J_{(j_1' j_q\ldots j_1)s,t}^{(1 i_q\ldots i_1)}+J_{(j_q j_1' j_{q-1}\ldots j_1)s,t}
^{(i_q 1 i_{q-1}\ldots i_1)}+\ldots +J_{(j_q\ldots j_1 j_1')s,t}^{(i_q\ldots i_1 1)}\right)=
$$
\begin{equation}
\label{new1700a}
=
\sum\limits_{(j_1,\ldots,j_q, j_1')}
J_{(j_q\ldots j_1 j_1')s,t}^{(i_q\ldots i_1 1)}
\end{equation}

\vspace{3mm}
\noindent
w.~p.~1, where $J_{(j_r\ldots j_1)s,t}^{(i_r\ldots i_1)}$
is defined by (\ref{new100001}).
The equality (\ref{new1600}) is proved for the case $n=1.$

Let us assume that the equality (\ref{new1600}) is true for $n=2, 3, \ldots, k-1$, and prove
its validity for $n=k.$

Applying (\ref{new1025}), (\ref{new1040}), (\ref{new1020})--(\ref{new1021}), we obtain w.~p.~1
$$
\sum\limits_{(j_1',\ldots,j_k')}
\int\limits_t^s 
\phi_{j_k'}(t_k')\ldots \int\limits_t^{t_2'}\phi_{j_1'}(t_1')d{\bf w}_{t_1'}^{(1)}\ldots 
d{\bf w}_{t_k'}^{(1)}=
$$

\vspace{-2mm}
$$
=\int\limits_t^{s}\phi_{j_k'}(\theta)d{\bf w}_{\theta}^{(1)}
\sum\limits_{(j_1',\ldots,j_{k-1}')}
\int\limits_t^s 
\phi_{j_{k-1}'}(t_{k-1})\ldots \int\limits_t^{t_2}\phi_{j_1'}(t_1)d{\bf w}_{t_1}^{(1)}\ldots 
d{\bf w}_{t_{k-1}}^{(1)}-
$$

\vspace{-2mm}
\begin{equation}
\label{new1800}
-\sum\limits_{(j_1',\ldots,j_{k-1}')}
\int\limits_t^s \phi_{j_k'}(\theta)\phi_{j_{k-1}'}(\theta)d\theta
\int\limits_t^s 
\phi_{j_{k-2}'}(t_{k-2})\ldots \int\limits_t^{t_2}\phi_{j_1'}(t_1)d{\bf w}_{t_1}^{(1)}\ldots 
d{\bf w}_{t_{k-2}}^{(1)}.
\end{equation}

\vspace{3mm}

Substituting $s=T$ in (\ref{new1800})
and applying the orthonormality of 
$\{\phi_j(x)\}_{j=0}^{\infty}$, we get w.~p.~1
$$
\sum\limits_{(j_1',\ldots,j_k')}
\int\limits_t^T
\phi_{j_k'}(t_k')\ldots \int\limits_t^{t_2'}\phi_{j_1'}(t_1')d{\bf w}_{t_1'}^{(1)}\ldots 
d{\bf w}_{t_k'}^{(1)}=
$$

\vspace{-4mm}
$$
=\int\limits_t^{T}\phi_{j_k'}(\theta)d{\bf w}_{\theta}^{(1)}
\sum\limits_{(j_1',\ldots,j_{k-1}')}
\int\limits_t^T 
\phi_{j_{k-1}'}(t_{k-1})\ldots \int\limits_t^{t_2}\phi_{j_1'}(t_1)d{\bf w}_{t_1}^{(1)}\ldots 
d{\bf w}_{t_{k-1}}^{(1)}-
$$

\begin{equation}
\label{new1801}
-\sum\limits_{(j_1',\ldots,j_{k-1}')}
{\bf 1}_{\{j_k'=j_{k-1}'\}}
\int\limits_t^T
\phi_{j_{k-2}'}(t_{k-2})\ldots \int\limits_t^{t_2}\phi_{j_1'}(t_1)d{\bf w}_{t_1}^{(1)}\ldots 
d{\bf w}_{t_{k-2}}^{(1)},
\end{equation}

\vspace{3mm}
\noindent
where ${\bf 1}_{A}$ is the indicator of the set $A$.

Using (\ref{new1801}) and the induction hypothesis, we obtain w.~p.~1
$$
\sum\limits_{(j_1',\ldots,j_k')}
\int\limits_t^T 
\phi_{j_k'}(t_k)\ldots \int\limits_t^{t_2}\phi_{j_1'}(t_1)d{\bf w}_{t_1}^{(1)}\ldots 
d{\bf w}_{t_k}^{(1)}\times
$$

\vspace{-2mm}
$$
\times\sum\limits_{(j_1,\ldots,j_q)}\int\limits_t^T \phi_{j_q}(t_q)\ldots
\int\limits_t^{t_2}\phi_{j_1}(t_1)
d{\bf w}_{t_1}^{(i_1)}\ldots d{\bf w}_{t_q}^{(i_q)}=
$$

\vspace{-2mm}
$$
=\int\limits_t^{T}\phi_{j_k'}(\theta)d{\bf w}_{\theta}^{(1)}
\sum\limits_{(j_1',\ldots,j_{k-1}')}
\int\limits_t^T 
\phi_{j_{k-1}'}(t_{k-1})\ldots \int\limits_t^{t_2}\phi_{j_1'}(t_1)d{\bf w}_{t_1}^{(1)}\ldots 
d{\bf w}_{t_{k-1}}^{(1)}\times
$$

$$
\times\sum\limits_{(j_1,\ldots,j_q)}\int\limits_t^T \phi_{j_q}(t_q)\ldots
\int\limits_t^{t_2}\phi_{j_1}(t_1)
d{\bf w}_{t_1}^{(i_1)}\ldots d{\bf w}_{t_q}^{(i_q)}-
$$

\vspace{-2mm}
$$
-\sum\limits_{(j_1',\ldots,j_{k-1}')}
{\bf 1}_{\{j_k'=j_{k-1}'\}}
\int\limits_t^T
\phi_{j_{k-2}'}(t_{k-2})\ldots \int\limits_t^{t_2}\phi_{j_1'}(t_1)d{\bf w}_{t_1}^{(1)}\ldots 
d{\bf w}_{t_{k-2}}^{(1)}\times
$$

\vspace{-2mm}
$$
\times\sum\limits_{(j_1,\ldots,j_q)}\int\limits_t^T \phi_{j_q}(t_q)\ldots
\int\limits_t^{t_2}\phi_{j_1}(t_1)
d{\bf w}_{t_1}^{(i_1)}\ldots d{\bf w}_{t_q}^{(i_q)}=
$$

\vspace{1mm}
$$
=\int\limits_t^{T}\phi_{j_k'}(\theta)d{\bf w}_{\theta}^{(1)}\times
$$

\vspace{-3mm}
$$
\times
\sum\limits_{(j_1,\ldots,j_q,j_1',\ldots,j_{k-1}')}
\int\limits_t^T \phi_{j_q}(t_q)\ldots \int\limits_t^{t_2}\phi_{j_1}(t_1)
\int\limits_t^{t_1}\phi_{j_{k-1}'}(t_{k-1}')\ldots \int\limits_t^{t_2'}
\phi_{j_1'}(t_1')\times
$$

\vspace{1mm}
$$
\times d{\bf w}_{t_1'}^{(1)}\ldots 
d{\bf w}_{t_{k-1}'}^{(1)}d{\bf w}_{t_1}^{(i_1)}\ldots d{\bf w}_{t_q}^{(i_q)}-
$$

\vspace{-2mm}
$$
-\sum\limits_{(j_1',\ldots,j_{k-1}')}
{\bf 1}_{\{j_k'=j_{k-1}'\}}
\int\limits_t^T
\phi_{j_{k-2}'}(t_{k-2})\ldots \int\limits_t^{t_2}\phi_{j_1'}(t_1)d{\bf w}_{t_1}^{(1)}\ldots 
d{\bf w}_{t_{k-2}}^{(1)}\times
$$
\begin{equation}
\label{newi9000}
\times\sum\limits_{(j_1,\ldots,j_q)}\int\limits_t^T \phi_{j_q}(t_q)\ldots
\int\limits_t^{t_2}\phi_{j_1}(t_1)
d{\bf w}_{t_1}^{(i_1)}\ldots d{\bf w}_{t_q}^{(i_q)}.
\end{equation}

\vspace{3mm}

Further, applying the induction hypothesis, we have w.~p.~1
$$
\sum\limits_{(j_1',\ldots,j_{k-1}')}
{\bf 1}_{\{j_k'=j_{k-1}'\}}
\int\limits_t^T
\phi_{j_{k-2}'}(t_{k-2})\ldots \int\limits_t^{t_2}\phi_{j_1'}(t_1)d{\bf w}_{t_1}^{(1)}\ldots 
d{\bf w}_{t_{k-2}}^{(1)}\times
$$

\vspace{-2mm}
$$
\times\sum\limits_{(j_1,\ldots,j_q)}\int\limits_t^T \phi_{j_q}(t_q)\ldots
\int\limits_t^{t_2}\phi_{j_1}(t_1)
d{\bf w}_{t_1}^{(i_1)}\ldots d{\bf w}_{t_q}^{(i_q)}=
$$

\vspace{-2mm}
$$
=\Biggl(
\sum\limits_{(j_1',\ldots,j_{k-2}')}
{\bf 1}_{\{j_k'=j_{k-1}'\}}
\int\limits_t^T
\phi_{j_{k-2}'}(t_{k-2})\ldots \int\limits_t^{t_2}\phi_{j_1'}(t_1)d{\bf w}_{t_1}^{(1)}\ldots 
d{\bf w}_{t_{k-2}}^{(1)}+
\Biggr.
$$

\vspace{-2mm}
$$
+\sum\limits_{(j_1',\ldots,j_{k-3}', j_{k-1}')}
{\bf 1}_{\{j_k'=j_{k-2}'\}}
\int\limits_t^T
\phi_{j_{k-1}'}(t_{k-2})
\int\limits_t^{t_{k-2}}
\phi_{j_{k-3}'}(t_{k-3})\ldots \int\limits_t^{t_2}\phi_{j_1'}(t_1)\times
$$

\vspace{1mm}
$$
\times d{\bf w}_{t_1}^{(1)}\ldots 
d{\bf w}_{t_{k-3}}^{(1)}d{\bf w}_{t_{k-2}}^{(1)}+\ldots
$$

\vspace{-2mm}
$$
\ldots+\sum\limits_{(j_2',\ldots,j_{k-1}')}
{\bf 1}_{\{j_k'=j_{1}'\}}
\int\limits_t^T
\phi_{j_{k-2}'}(t_{k-2})\ldots \int\limits_t^{t_3}\phi_{j_2'}(t_2)
\int\limits_t^{t_2}
\phi_{j_{k-1}'}(t_{1})\times
$$

$$
\Biggl.\times
d{\bf w}_{t_1}^{(1)}d{\bf w}_{t_2}^{(1)}\ldots 
d{\bf w}_{t_{k-2}}^{(1)}\Biggr)\times
$$

\vspace{-1mm}
$$
\times\sum\limits_{(j_1,\ldots,j_q)}\int\limits_t^T \phi_{j_q}(t_q)\ldots
\int\limits_t^{t_2}\phi_{j_1}(t_1)
d{\bf w}_{t_1}^{(i_1)}\ldots d{\bf w}_{t_q}^{(i_q)}=
$$

\vspace{-2mm}
$$
=\Biggl(
{\bf 1}_{\{j_k'=j_{k-1}'\}}\sum\limits_{(j_1',\ldots,j_{k-2}')}
\int\limits_t^T
\phi_{j_{k-2}'}(t_{k-2})\ldots \int\limits_t^{t_2}\phi_{j_1'}(t_1)d{\bf w}_{t_1}^{(1)}\ldots 
d{\bf w}_{t_{k-2}}^{(1)}+
\Biggr.
$$
$$
+{\bf 1}_{\{j_k'=j_{k-2}'\}}\sum\limits_{(j_1',\ldots,j_{k-3}', j_{k-1}')}
\int\limits_t^T
\phi_{j_{k-1}'}(t_{k-2})
\int\limits_t^{t_{k-2}}
\phi_{j_{k-3}'}(t_{k-3})\ldots \int\limits_t^{t_2}\phi_{j_1'}(t_1)\times
$$

\vspace{1mm}
$$
\times d{\bf w}_{t_1}^{(1)}\ldots 
d{\bf w}_{t_{k-3}}^{(1)}d{\bf w}_{t_{k-2}}^{(1)}+\ldots
$$

\vspace{-2mm}
$$
\ldots +{\bf 1}_{\{j_k'=j_{1}'\}}\sum\limits_{(j_2',\ldots,j_{k-1}')}
\int\limits_t^T
\phi_{j_{k-2}'}(t_{k-2})\ldots \int\limits_t^{t_3}\phi_{j_2'}(t_2)
\int\limits_t^{t_2}
\phi_{j_{k-1}'}(t_{1})\times
$$

$$
\Biggl.\times
d{\bf w}_{t_1}^{(1)}d{\bf w}_{t_2}^{(1)}\ldots 
d{\bf w}_{t_{k-2}}^{(1)}\Biggr)\times
$$

\vspace{-2mm}
$$
\times\sum\limits_{(j_1,\ldots,j_q)}\int\limits_t^T \phi_{j_q}(t_q)\ldots
\int\limits_t^{t_2}\phi_{j_1}(t_1)
d{\bf w}_{t_1}^{(i_1)}\ldots d{\bf w}_{t_q}^{(i_q)}=
$$

\vspace{-2mm}
$$
={\bf 1}_{\{j_k'=j_{k-1}'\}}\sum\limits_{(j_1,\ldots,j_q,j_1',\ldots,j_{k-2}')}
\int\limits_t^T \phi_{j_q}(t_q)\ldots \int\limits_t^{t_2}\phi_{j_1}(t_1)
\int\limits_t^{t_1}\phi_{j_{k-2}'}(t_{k-2}')\ldots \int\limits_t^{t_2'}
\phi_{j_1'}(t_1')\times
$$

\vspace{1mm}
$$
\times d{\bf w}_{t_1'}^{(1)}\ldots d{\bf w}_{t_{k-2}'}^{(1)}
d{\bf w}_{t_1}^{(i_1)}\ldots d{\bf w}_{t_q}^{(i_q)}+
$$

\vspace{-2mm}
$$
+{\bf 1}_{\{j_k'=j_{k-2}'\}}\sum\limits_{(j_1,\ldots,j_q,j_1',\ldots,j_{k-3}',j_{k-1}')}
\int\limits_t^T \phi_{j_q}(t_q)\ldots \int\limits_t^{t_2}\phi_{j_1}(t_1)
\int\limits_t^{t_1}\phi_{j_{k-1}'}(t_{k-2}')\times
$$

\vspace{-2mm}
$$
\times\int\limits_t^{t_{k-2}'}
\phi_{j_{k-3}'}(t_{k-3}')
\ldots \int\limits_t^{t_2'}
\phi_{j_1'}(t_1')d{\bf w}_{t_1'}^{(1)}\ldots d{\bf w}_{t_{k-3}'}^{(1)}d{\bf w}_{t_{k-2}'}^{(1)}
d{\bf w}_{t_1}^{(i_1)}\ldots d{\bf w}_{t_q}^{(i_q)}+
$$

\vspace{-5mm}
$$
\ldots
$$

\vspace{-5mm}
$$
+{\bf 1}_{\{j_k'=j_{1}'\}}\sum\limits_{(j_1,\ldots,j_q,j_2',\ldots,j_{k-1}')}
\int\limits_t^T \phi_{j_q}(t_q)\ldots \int\limits_t^{t_2}\phi_{j_1}(t_1)\times
$$

\vspace{-2mm}
$$
\times
\int\limits_t^{t_1}\phi_{j_{k-2}'}(t_{k-2}')\ldots \int\limits_t^{t_{3}'}
\phi_{j_{2}'}(t_{2}')
\int\limits_t^{t_2'}
\phi_{j_{k-1}'}(t_1')d{\bf w}_{t_1'}^{(1)}d{\bf w}_{t_2'}^{(1)}\ldots d{\bf w}_{t_{k-2}'}^{(1)}
d{\bf w}_{t_1}^{(i_1)}\ldots d{\bf w}_{t_q}^{(i_q)}\stackrel{\sf def}{=}
$$
\begin{equation}
\label{new1900}
\stackrel{\sf def}{=}S_4(T).
\end{equation}

\vspace{3mm}

By analogy with (\ref{newx1020ss}) we obtain w.~p.~1
$$
\int\limits_t^T \phi_l(\tau)\phi_{j_r}(\tau)d\tau
\int\limits_t^{T} \phi_{j_{r-1}}(t_{r-1})\ldots \int\limits_t^{t_2}
\phi_{j_1}(t_1)d{\bf w}_{t_1}^{(i_1)}\ldots d{\bf w}_{t_{r-1}}^{(i_{r-1})}=
$$

\vspace{-2mm}
$$
=\int\limits_t^T \phi_l(\tau)\phi_{j_r}(\tau)
\int\limits_t^{\tau} \phi_{j_{r-1}}(t_{r-1})\ldots \int\limits_t^{t_2}
\phi_{j_1}(t_1)d{\bf w}_{t_1}^{(i_1)}\ldots d{\bf w}_{t_{r-1}}^{(i_{r-1})}d\tau+\ldots
$$

\vspace{-2mm}
\begin{equation}
\label{new5001}
\ldots +
\int\limits_t^{T} \phi_{j_{r-1}}(t_{r-1})\ldots \int\limits_t^{t_2}
\phi_{j_1}(t_1)\int\limits_t^{t_1} \phi_{l}(\tau)\phi_{j_r}(\tau)      
d\tau d{\bf w}_{t_1}^{(i_1)}\ldots d{\bf w}_{t_{r-1}}^{(i_{r-1})},
\end{equation}

\vspace{3mm}
\noindent
where $i_1,\ldots,i_{r-1}=0,1,\ldots,m.$

Using iteratively the It\^{o} formula, as well as (\ref{new5001})
and combinatorial reasoning, we get w.~p.~1 (see Remark~2 below for details)

\vspace{-2mm}
$$
\int\limits_t^{T}\phi_{j_k'}(\theta)d{\bf w}_{\theta}^{(1)}\times
$$

\vspace{-4mm}
$$
\times
\sum\limits_{(j_1,\ldots,j_q,j_1',\ldots,j_{k-1}')}
\int\limits_t^T \phi_{j_q}(t_q)\ldots \int\limits_t^{t_2}\phi_{j_1}(t_1)
\int\limits_t^{t_1}\phi_{j_{k-1}'}(t_{k-1}')\ldots \int\limits_t^{t_2'}
\phi_{j_1'}(t_1')\times
$$

\vspace{2mm}
$$
\times d{\bf w}_{t_1'}^{(1)}\ldots 
d{\bf w}_{t_{k-1}'}^{(1)}d{\bf w}_{t_1}^{(i_1)}\ldots d{\bf w}_{t_q}^{(i_q)}=
$$

\vspace{2mm}
$$
=\sum\limits_{(j_1,\ldots,j_q,j_1',\ldots,j_{k}')}
\int\limits_t^T \phi_{j_q}(t_q)\ldots \int\limits_t^{t_2}\phi_{j_1}(t_1)
\int\limits_t^{t_1}\phi_{j_{k}'}(t_{k}')\ldots \int\limits_t^{t_2'}
\phi_{j_1'}(t_1')\times
$$

\vspace{2mm}
$$
\times d{\bf w}_{t_1'}^{(1)}\ldots 
d{\bf w}_{t_{k}'}^{(1)}d{\bf w}_{t_1}^{(i_1)}\ldots d{\bf w}_{t_q}^{(i_q)}+
$$
$$
+\sum\limits_{(j_1,\ldots,j_q,j_1',\ldots,j_{k-1}')}
\Biggl(\int\limits_t^T \phi_{j_q}(t_q)\ldots \int\limits_t^{t_2}\phi_{j_1}(t_1)
\int\limits_t^{t_1}\phi_{j_{k}'}(\theta)\phi_{j_{k-1}'}(\theta)\int\limits_t^{\theta}
\phi_{j_{k-2}'}(t_{k-2}')\ldots\Biggr.
$$
$$
\Biggl.\ldots \int\limits_t^{t_2'}\phi_{j_1'}(t_1')
d{\bf w}_{t_1'}^{(1)}\ldots 
d{\bf w}_{t_{k-2}'}^{(1)}d{\bf w}_{\theta}^{(0)}d{\bf w}_{t_1}^{(i_1)}\ldots d{\bf w}_{t_q}^{(i_q)}+
$$

$$
+
\int\limits_t^T \phi_{j_q}(t_q)\ldots \int\limits_t^{t_2}
\phi_{j_1}(t_1)\int\limits_t^{t_1}\phi_{j_{k-1}'}(t_{k-1}')
\int\limits_t^{t_{k-1}'}\phi_{j_{k}'}(\theta)\phi_{j_{k-2}'}(\theta)\int\limits_t^{\theta}
\phi_{j_{k-3}'}(t_{k-3}')\ldots
$$
$$
\ldots
\int\limits_t^{t_2'}
\phi_{j_{1}'}(t_{1}')d{\bf w}_{t_1'}^{(1)}\ldots 
d{\bf w}_{t_{k-3}'}^{(1)}d{\bf w}_{\theta}^{(0)}d{\bf w}_{t_{k-1}'}^{(1)}
d{\bf w}_{t_1}^{(i_1)}\ldots d{\bf w}_{t_q}^{(i_q)}+\ldots
$$

$$
\ldots +
\int\limits_t^T \phi_{j_q}(t_q)\ldots \int\limits_t^{t_2}
\phi_{j_1}(t_1)\int\limits_t^{t_1}\phi_{j_{k-1}'}(t_{k-1}')\ldots
\int\limits_t^{t_{3}'}\phi_{j_{2}'}(t_2')\int\limits_t^{t_{2}'}
\phi_{j_{k}'}(\theta)\phi_{j_{1}'}(\theta)d{\bf w}_{\theta}^{(0)}\times
$$
$$
\Biggl.\times
d{\bf w}_{t_2'}^{(1)}\ldots 
d{\bf w}_{t_{k-1}'}^{(1)}
d{\bf w}_{t_1}^{(i_1)}\ldots d{\bf w}_{t_q}^{(i_q)}\Biggr)=
$$

\vspace{2mm}
$$
=\sum\limits_{(j_1,\ldots,j_q,j_1',\ldots,j_{k}')}
\int\limits_t^T \phi_{j_q}(t_q)\ldots \int\limits_t^{t_2}\phi_{j_1}(t_1)
\int\limits_t^{t_1}\phi_{j_{k}'}(t_{k}')\ldots \int\limits_t^{t_2'}
\phi_{j_1'}(t_1')\times
$$

$$
\times d{\bf w}_{t_1'}^{(1)}\ldots 
d{\bf w}_{t_{k}'}^{(1)}d{\bf w}_{t_1}^{(i_1)}\ldots d{\bf w}_{t_q}^{(i_q)}+
$$

$$
+\sum\limits_{(j_1,\ldots,j_q,j_1',\ldots,j_{k-2}')}
\Biggl\{\int\limits_t^T 
\phi_{j_{k}'}(\theta)\phi_{j_{k-1}'}(\theta)\int\limits_t^{\theta}
\phi_{j_q}(t_q)\ldots \int\limits_t^{t_2}\phi_{j_1}(t_1)
\int\limits_t^{t_1}
\phi_{j_{k-2}'}(t_{k-2}')\ldots\Biggr.
$$
$$
\Biggl.\ldots \int\limits_t^{t_2'}\phi_{j_1'}(t_1')
d{\bf w}_{t_1'}^{(1)}\ldots 
d{\bf w}_{t_{k-2}'}^{(1)}d{\bf w}_{t_1}^{(i_1)}\ldots d{\bf w}_{t_q}^{(i_q)}d{\bf w}_{\theta}^{(0)}+\ldots
$$
$$
\ldots+\int\limits_t^T 
\phi_{j_q}(t_q)\ldots \int\limits_t^{t_2}\phi_{j_1}(t_1)
\int\limits_t^{t_1}
\phi_{j_{k-2}'}(t_{k-2}')\ldots \int\limits_t^{t_2'}\phi_{j_1'}(t_1')
\int\limits_t^{t_1'} 
\phi_{j_{k}'}(\theta)\phi_{j_{k-1}'}(\theta)d{\bf w}_{\theta}^{(0)}\times
$$
$$
\Biggl.
\times d{\bf w}_{t_1'}^{(1)}\ldots 
d{\bf w}_{t_{k-2}'}^{(1)}d{\bf w}_{t_1}^{(i_1)}\ldots d{\bf w}_{t_q}^{(i_q)}\Biggr\}+
$$

$$
+\sum\limits_{(j_1,\ldots,j_q,j_1',\ldots,j_{k-3}',j_{k-1}')}
\Biggl\{\int\limits_t^T 
\phi_{j_{k}'}(\theta)\phi_{j_{k-2}'}(\theta)\int\limits_t^{\theta}
\phi_{j_q}(t_q)\ldots \int\limits_t^{t_2}\phi_{j_1}(t_1)
\int\limits_t^{t_1}
\phi_{j_{k-1}'}(t_{k-1}')\times\Biggr.
$$
$$
\times\int\limits_t^{t_{k-1}'}
\phi_{j_{k-3}'}(t_{k-3}')
\ldots \int\limits_t^{t_2'}\phi_{j_1'}(t_1')
d{\bf w}_{t_1'}^{(1)}\ldots 
d{\bf w}_{t_{k-3}'}^{(1)}d{\bf w}_{t_{k-1}'}^{(1)}
d{\bf w}_{t_1}^{(i_1)}\ldots d{\bf w}_{t_q}^{(i_q)}
d{\bf w}_{\theta}^{(0)}+\ldots
$$

$$
\ldots+\int\limits_t^T 
\phi_{j_q}(t_q)\ldots \int\limits_t^{t_2}\phi_{j_1}(t_1)
\int\limits_t^{t_1}
\phi_{j_{k-1}'}(t_{k-1}')\int\limits_t^{t_{k-1}'}
\phi_{j_{k-3}'}(t_{k-3}')\ldots
\int\limits_t^{t_2'}\phi_{j_1'}(t_1')
\times
$$
$$
\Biggl.\times
\int\limits_t^{t_1'} 
\phi_{j_{k}'}(\theta)\phi_{j_{k-2}'}(\theta)d{\bf w}_{\theta}^{(0)}
d{\bf w}_{t_1'}^{(1)}\ldots 
d{\bf w}_{t_{k-3}'}^{(1)}d{\bf w}_{t_{k-1}'}^{(1)}
d{\bf w}_{t_1}^{(i_1)}\ldots d{\bf w}_{t_q}^{(i_q)}\Biggr\}+\ldots
$$

$$
\ldots +\sum\limits_{(j_1,\ldots,j_q,j_2',\ldots,j_{k-1}')}
\Biggl\{\int\limits_t^T 
\phi_{j_{k}'}(\theta)\phi_{j_{1}'}(\theta)\int\limits_t^{\theta}
\phi_{j_q}(t_q)\ldots \int\limits_t^{t_2}\phi_{j_1}(t_1)
\int\limits_t^{t_1}
\phi_{j_{k-1}'}(t_{k-1}')\ldots\Biggr.
$$
$$
\ldots\int\limits_t^{t_{3}'}
\phi_{j_2'}(t_2')
d{\bf w}_{t_2'}^{(1)}\ldots 
d{\bf w}_{t_{k-1}'}^{(1)}
d{\bf w}_{t_1}^{(i_1)}\ldots d{\bf w}_{t_q}^{(i_q)}
d{\bf w}_{\theta}^{(0)}+\ldots
$$

$$
\ldots+\int\limits_t^T 
\phi_{j_q}(t_q)\ldots \int\limits_t^{t_2}\phi_{j_1}(t_1)
\int\limits_t^{t_1}
\phi_{j_{k-1}'}(t_{k-1}')\ldots
\int\limits_t^{t_3'}\phi_{j_2'}(t_2')\int\limits_t^{t_2'}
\phi_{j_{k}'}(\theta)\phi_{j_{1}'}(\theta)d{\bf w}_{\theta}^{(0)}\times
$$
$$
\Biggl.\times
d{\bf w}_{t_2'}^{(1)}\ldots 
d{\bf w}_{t_{k-1}'}^{(1)}
d{\bf w}_{t_1}^{(i_1)}\ldots d{\bf w}_{t_q}^{(i_q)}\Biggr\}=
$$
$$
=\sum\limits_{(j_1,\ldots,j_q,j_1',\ldots,j_{k}')}
\int\limits_t^T \phi_{j_q}(t_q)\ldots \int\limits_t^{t_2}\phi_{j_1}(t_1)
\int\limits_t^{t_1}\phi_{j_{k}'}(t_{k}')\ldots \int\limits_t^{t_2'}
\phi_{j_1'}(t_1')\times
$$

$$
\times d{\bf w}_{t_1'}^{(1)}\ldots 
d{\bf w}_{t_{k}'}^{(1)}d{\bf w}_{t_1}^{(i_1)}\ldots d{\bf w}_{t_q}^{(i_q)}+
$$

$$
+\int\limits_t^T 
\phi_{j_{k}'}(\theta)\phi_{j_{k-1}'}(\theta)d\theta\sum\limits_{(j_1,\ldots,j_q,j_1',\ldots,j_{k-2}')}
\int\limits_t^{T}
\phi_{j_q}(t_q)\ldots \int\limits_t^{t_2}\phi_{j_1}(t_1)
\int\limits_t^{t_1}
\phi_{j_{k-2}'}(t_{k-2}')\ldots
$$
$$
\ldots \int\limits_t^{t_2'}\phi_{j_1'}(t_1')
d{\bf w}_{t_1'}^{(1)}\ldots 
d{\bf w}_{t_{k-2}'}^{(1)}d{\bf w}_{t_1}^{(i_1)}\ldots d{\bf w}_{t_q}^{(i_q)}+
$$

$$
+\int\limits_t^T 
\phi_{j_{k}'}(\theta)\phi_{j_{k-2}'}(\theta)d\theta
\sum\limits_{(j_1,\ldots,j_q,j_1',\ldots,j_{k-3}',j_{k-1}')}
\int\limits_t^{T}
\phi_{j_q}(t_q)\ldots \int\limits_t^{t_2}\phi_{j_1}(t_1)
\int\limits_t^{t_1}
\phi_{j_{k-1}'}(t_{k-1}')\times\Biggr.
$$
$$
\times\int\limits_t^{t_{k-1}'}
\phi_{j_{k-3}'}(t_{k-3}')
\ldots \int\limits_t^{t_2'}\phi_{j_1'}(t_1')
d{\bf w}_{t_1'}^{(1)}\ldots 
d{\bf w}_{t_{k-3}'}^{(1)}d{\bf w}_{t_{k-1}'}^{(1)}
d{\bf w}_{t_1}^{(i_1)}\ldots d{\bf w}_{t_q}^{(i_q)}
+\ldots
$$

$$
\ldots +\int\limits_t^T 
\phi_{j_{k}'}(\theta)\phi_{j_{1}'}(\theta)d\theta\sum\limits_{(j_1,\ldots,j_q,j_2',\ldots,j_{k-1}')}
\int\limits_t^{T}
\phi_{j_q}(t_q)\ldots \int\limits_t^{t_2}\phi_{j_1}(t_1)
\int\limits_t^{t_1}
\phi_{j_{k-1}'}(t_{k-1}')\ldots
$$
$$
\ldots\int\limits_t^{t_{3}'}
\phi_{j_2'}(t_2')
d{\bf w}_{t_2'}^{(1)}\ldots 
d{\bf w}_{t_{k-1}'}^{(1)}
d{\bf w}_{t_1}^{(i_1)}\ldots d{\bf w}_{t_q}^{(i_q)}=
$$

$$
=\sum\limits_{(j_1,\ldots,j_q,j_1',\ldots,j_{k}')}
\int\limits_t^T \phi_{j_q}(t_q)\ldots \int\limits_t^{t_2}\phi_{j_1}(t_1)
\int\limits_t^{t_1}\phi_{j_{k}'}(t_{k}')\ldots \int\limits_t^{t_2'}
\phi_{j_1'}(t_1')\times
$$

\vspace{1mm}
\begin{equation}
\label{new1901}
\times d{\bf w}_{t_1'}^{(1)}\ldots 
d{\bf w}_{t_{k}'}^{(1)}d{\bf w}_{t_1}^{(i_1)}\ldots d{\bf w}_{t_q}^{(i_q)}+S_4(T).
\end{equation}

\vspace{5mm}

From (\ref{newi9000}), (\ref{new1900}), and (\ref{new1901}) we conclude that
the equality (\ref{new1600}) is proved for $n=k.$
The equality (\ref{new1600}) is proved.

\vspace{2mm}

{\bf Remark~2.}\ {\it It should be noted that the sums with respect to 
permutations 
$$
\sum\limits_{(j_1,\ldots,j_q,j_1',\ldots,j_{k-1}')}
$$

\noindent
in {\rm (\ref{new1901}),} containing the expressions 
$\phi_{j_{k}'}(\theta)\phi_{j_{k-1}'}(\theta),\ldots,
\phi_{j_{k}'}(\theta)\phi_{j_{1}'}(\theta),$
should be understood in a special way.
Let us explain this rule on the following sum

\vspace{-6mm}
$$
\sum\limits_{(j_1,\ldots,j_q,j_1',\ldots,j_{k-1}')}
\int\limits_t^T \phi_{j_q}(t_q)\ldots \int\limits_t^{t_2}\phi_{j_1}(t_1)
\int\limits_t^{t_1}\phi_{j_{k}'}(\theta)\phi_{j_{k-1}'}(\theta)\int\limits_t^{\theta}
\phi_{j_{k-2}'}(t_{k-2}')\ldots
$$
\begin{equation}
\label{new777100}
\Biggl.\ldots \int\limits_t^{t_2'}\phi_{j_1'}(t_1')
d{\bf w}_{t_1'}^{(1)}\ldots 
d{\bf w}_{t_{k-2}'}^{(1)}d{\bf w}_{\theta}^{(0)}d{\bf w}_{t_1}^{(i_1)}\ldots d{\bf w}_{t_q}^{(i_q)}.
\end{equation}

More precisely, permutations $\left(j_1,\ldots,j_q,j_1',\ldots,j_{k-1}'\right)$ 
when summing in {\rm (\ref{new777100})}
are performed in such a way that if
$j_r^{*}$ swapped with $j_d^{*}$ in the  
permutation $\left(j_{q+k-1}^{*},\ldots,j_1^{*}\right)=
\left(j_q,\ldots,j_1,j_{k-1}',j_{k-2}',\ldots,j_{1}'\right),$ 
then $i_r^{*}$ swapped with $i_d^{*}$ in 
the permutation 
$$
\left(i_{q+k-1}^{*},\ldots,i_1^{*}\right)=
\bigl(i_q,\ldots,i_1,0,\underbrace{1, \ldots ,1}_{k-2}\bigr).
$$

\noindent
Moreover, 
$\bar \phi_{j_r^{*}}$ swapped with $\bar \phi_{j_d^{*}}$
in the permutation 
$$
\bigl(\bar \phi_{j_{q+k-1}^{*}},\ldots,\bar\phi_{j_1^{*}}\bigr)=
\bigl(\phi_{j_q},\ldots,\phi_{j_1},\hspace{1.5mm} \phi_{j_{k}'}\hspace{-0.5mm}\cdot\hspace{-0.5mm}
\phi_{j_{k-1}'},\hspace{1.5mm}
\phi_{j_{k-2}'},\ldots, \phi_{j_{1}'}\bigr).
$$

\vspace{1mm}
\noindent
A similar rule should be applied to all other sums with respect to permutations
$$
\sum\limits_{(j_1,\ldots,j_q,j_1',\ldots,j_{k-1}')}
$$

\noindent
in {\rm (\ref{new1901})} that contain the expressions
$\phi_{j_{k}'}(\theta)\phi_{j_{k-2}'}(\theta),\ldots,
\phi_{j_{k}'}(\theta)\phi_{j_{1}'}(\theta).$}

\vspace{2mm}

Let us prove the equality (\ref{new1600a}). Consider the case  $n=1.$
By analogy with (\ref{new1700}) and (\ref{new1700a}) we obtain 
$$
\int\limits_t^{s}\phi_{j_1'}(\theta)d{\bf w}_{\theta}^{(0)}
\sum\limits_{(j_1,\ldots,j_q)}\int\limits_t^s \phi_{j_q}(t_q)\ldots
\int\limits_t^{t_2}\phi_{j_1}(t_1)
d{\bf w}_{t_1}^{(i_1)}\ldots {\bf w}_{t_q}^{(i_q)}=
\sum\limits_{(j_1,\ldots,j_q, j_1')}
J_{(j_q\ldots j_1 j_1')s,t}^{(i_q\ldots i_1 0)}
$$
w.~p.~1, where $J_{(j_r\ldots j_1)s,t}^{(i_r\ldots i_1)}$
is defined by (\ref{new100001}).
The equality (\ref{new1600a}) is proved for the case $n=1.$

Let us assume that the equality (\ref{new1600a}) is true for $n=2, 3, \ldots, k-1$, and prove
its validity for $n=k.$

In complete analogy with (\ref{new1301}) we get
$$
\int\limits_t^{s}\phi_{j_k'}(\theta)d\theta
\int\limits_t^s \phi_{j_{k-1}'}(t_{k-1})\ldots
\int\limits_t^{t_2}\phi_{j_1'}(t_1)
dt_1\ldots dt_{k-1}=
$$
\begin{equation}
\label{new5000}
=
J_{(j_k'j_{k-1}'\ldots j_1')s,t}^{(0\ldots  0)}
+J_{(j_{k-1}'j_k' j_{k-2}'\ldots j_1')s,t}^{(0\ldots 0)}+
\ldots + J_{(j_{k-1}'\ldots j_1' j_k')s,t}^{(0\ldots 0)}.
\end{equation}

\vspace{5mm}

Applying (\ref{new5000}), we have
$$
\sum\limits_{(j_1',\ldots,j_k')}
\int\limits_t^T \phi_{j_k'}(t_k')\ldots \int\limits_t^{t_2'}\phi_{j_1'}(t_1')
d{\bf w}_{t_1'}^{(0)}\ldots d{\bf w}_{t_k'}^{(0)}=
$$

\vspace{-2mm}
$$
=\sum\limits_{(j_1',\ldots,j_{k-1}')}
\left(J_{(j_k'j_{k-1}'\ldots j_1')s,t}^{(0\ldots  0)}
+J_{(j_{k-1}'j_k' j_{k-2}'\ldots j_1')s,t}^{(0\ldots 0)}+
\ldots + J_{(j_{k-1}'\ldots j_1' j_k')s,t}^{(0\ldots 0)}\right)=
$$
\begin{equation}
\label{new3005}
~~~~=\int\limits_t^{T}\phi_{j_k'}(\theta)d\theta\sum\limits_{(j_1',\ldots,j_{k-1}')}
\int\limits_t^T \phi_{j_{k-1}'}(t_{k-1})\ldots \int\limits_t^{t_2'}\phi_{j_1'}(t_1)
d{\bf w}_{t_1}^{(0)}\ldots d{\bf w}_{t_{k-1}}^{(0)}.
\end{equation}

\vspace{3mm}

Using (\ref{new3005}) and the induction hypothesis, we obtain w.~p.~1
$$
\sum\limits_{(j_1',\ldots,j_k')}
\int\limits_t^T 
\phi_{j_k'}(t_k)\ldots \int\limits_t^{t_2}\phi_{j_1'}(t_1)d{\bf w}_{t_1}^{(0)}\ldots 
d{\bf w}_{t_k}^{(0)}\times
$$

\vspace{-2mm}
$$
\times\sum\limits_{(j_1,\ldots,j_q)}\int\limits_t^T \phi_{j_q}(t_q)\ldots
\int\limits_t^{t_2}\phi_{j_1}(t_1)
d{\bf w}_{t_1}^{(i_1)}\ldots d{\bf w}_{t_q}^{(i_q)}=
$$

\vspace{-3mm}
$$
=\int\limits_t^{T}\phi_{j_k'}(\theta)d\theta\sum\limits_{(j_1',\ldots,j_{k-1}')}
\int\limits_t^T \phi_{j_{k-1}'}(t_{k-1}')\ldots \int\limits_t^{t_2'}\phi_{j_1'}(t_1')
d{\bf w}_{t_1'}^{(0)}\ldots d{\bf w}_{t_{k-1}'}^{(0)}\times
$$
$$
\times\sum\limits_{(j_1,\ldots,j_q)}\int\limits_t^T \phi_{j_q}(t_q)\ldots
\int\limits_t^{t_2}\phi_{j_1}(t_1)
d{\bf w}_{t_1}^{(i_1)}\ldots d{\bf w}_{t_q}^{(i_q)}=
$$

$$
=\int\limits_t^{T}\phi_{j_k'}(\theta)d\theta
\sum\limits_{(j_1,\ldots,j_q,j_1',\ldots,j_{k-1}')}
\int\limits_t^T \phi_{j_q}(t_q)\ldots \int\limits_t^{t_2}\phi_{j_1}(t_1)\times
$$

\vspace{-1mm}
$$
\times
\int\limits_t^{t_1}\phi_{j_{k-1}'}(t_{k-1}')\ldots \int\limits_t^{t_2'}
\phi_{j_1'}(t_1')d{\bf w}_{t_1'}^{(0)}\ldots d{\bf w}_{t_{k-1}'}^{(0)}
d{\bf w}_{t_1}^{(i_1)}\ldots d{\bf w}_{t_q}^{(i_q)}=
$$

$$
=
\sum\limits_{(j_1,\ldots,j_q,j_1',\ldots,j_{k-1}')}\int\limits_t^{T}\phi_{j_k'}(\theta)d\theta
\int\limits_t^T \phi_{j_q}(t_q)\ldots \int\limits_t^{t_2}\phi_{j_1}(t_1)\times
$$

\begin{equation}
\label{new4000}
\times
\int\limits_t^{t_1}\phi_{j_{k-1}'}(t_{k-1}')\ldots \int\limits_t^{t_2'}
\phi_{j_1'}(t_1')d{\bf w}_{t_1'}^{(0)}\ldots d{\bf w}_{t_{k-1}'}^{(0)}
d{\bf w}_{t_1}^{(i_1)}\ldots d{\bf w}_{t_q}^{(i_q)}.
\end{equation}

\vspace{4mm}

An iterative application of the It\^{o} formula leads to the following equality

\vspace{-2mm}
$$
\int\limits_t^{T}\phi_{j_k'}(\theta)d\theta
\int\limits_t^T \phi_{j_q}(t_q)\ldots \int\limits_t^{t_2}\phi_{j_1}(t_1)\times
$$
$$
\times
\int\limits_t^{t_1}\phi_{j_{k-1}'}(t_{k-1}')\ldots \int\limits_t^{t_2'}
\phi_{j_1'}(t_1')d{\bf w}_{t_1'}^{(0)}\ldots d{\bf w}_{t_{k-1}'}^{(0)}
d{\bf w}_{t_1}^{(i_1)}\ldots d{\bf w}_{t_q}^{(i_q)}=
$$

$$
=J_{(j_k'j_q \ldots j_1 j_{k-1}'\ldots j_1')T,t}^{(0 i_q\ldots i_1 0\ldots 0)}+
J_{(j_q j_k'j_{q-1} \ldots j_1 j_{k-1}'\ldots j_1')T,t}^{(i_q 0 i_{q-1}\ldots i_1 0\ldots 0)}+\ldots
J_{(j_q \ldots j_1 j_k' j_{k-1}'\ldots j_1')T,t}^{(i_q\ldots i_1 0\ldots 0)}+
$$

\vspace{1mm}
\begin{equation}
\label{new4001}
+J_{(j_q \ldots j_1 j_{k-1}' j_k' j_{k-2}'\ldots j_1')T,t}^{(i_q\ldots i_1 0\ldots 0)}+\ldots
+J_{(j_q \ldots j_1 j_{k-1}' \ldots j_1' j_k')T,t}^{(i_q\ldots i_1 0\ldots 0)}
\end{equation}

\vspace{5mm}
\noindent
w.~p.~1.

Combining (\ref{new4000}) and (\ref{new4001}) we finally obtain w.~p.~1
$$
\sum\limits_{(j_1,\ldots,j_q)}\int\limits_t^T \phi_{j_q}(t_q)\ldots
\int\limits_t^{t_2}\phi_{j_1}(t_1)
d{\bf w}_{t_1}^{(i_1)}\ldots d{\bf w}_{t_q}^{(i_q)}\times
$$
$$
\times \sum\limits_{(j_1',\ldots,j_k')}
\int\limits_t^T \phi_{j_k'}(t_k')\ldots \int\limits_t^{t_2'}\phi_{j_1'}(t_1')
d{\bf w}_{t_1'}^{(0)}\ldots d{\bf w}_{t_k'}^{(0)}=
$$

\vspace{1mm}
$$
=\sum\limits_{(j_1,\ldots,j_q,j_1',\ldots,j_k')}
\int\limits_t^T \phi_{j_q}(t_q)\ldots \int\limits_t^{t_2}\phi_{j_1}(t_1)
\int\limits_t^{t_1}\phi_{j_k'}(t_k')\ldots \int\limits_t^{t_2'}
\phi_{j_1'}(t_1')\times
$$

\vspace{3mm}
$$
\times d{\bf w}_{t_1'}^{(0)}\ldots d{\bf w}_{t_k'}^{(0)}d{\bf w}_{t_1}^{(i_1)}\ldots d{\bf w}_{t_q}^{(i_q)}.
$$

\vspace{5mm}

The equality (\ref{new1600a}) is proved for $n=k.$
The equality (\ref{new1600a}) is proved. Theorem~4 is proved.

Let us consider the following theorem.

{\bf Theorem~5.}\ {\it Suppose that
$\{\phi_j(x)\}_{j=0}^{\infty}$ is an arbitrary complete orthonormal system  
of functions in the space $L_2([t,T]).$
Then the following representation

\vspace{-1mm}
$$
J''[\phi_{j_1}\ldots\phi_{j_k}]^{(i_1 \ldots i_k)}_{T,t}=
\prod_{l=1}^k\zeta_{j_l}^{(i_l)}+\sum\limits_{r=1}^{[k/2]}
(-1)^r \times \Biggr.
$$

\vspace{-3mm}
\begin{equation}
\label{leto6000xxa}
\times
\sum_{\stackrel{(\{\{g_1, g_2\}, \ldots, 
\{g_{2r-1}, g_{2r}\}\}, \{q_1, \ldots, q_{k-2r}\})}
{{}_{\{g_1, g_2, \ldots, 
g_{2r-1}, g_{2r}, q_1, \ldots, q_{k-2r}\}=\{1,2, \ldots, k\}}}}
\prod\limits_{s=1}^r
{\bf 1}_{\{i_{g_{{}_{2s-1}}}=~i_{g_{{}_{2s}}}\ne 0\}}
\Biggl.{\bf 1}_{\{j_{g_{{}_{2s-1}}}=~j_{g_{{}_{2s}}}\}}
\prod_{l=1}^{k-2r}\zeta_{j_{q_l}}^{(i_{q_l})}
\end{equation}

\vspace{2mm}
\noindent 
is valid w.~p.~{\rm 1,} where $i_1,\ldots,i_k=0,1,\ldots,m,$
$J''[\phi_{j_1}\ldots\phi_{j_k}]^{(i_1 \ldots i_k)}_{T,t}$
is defined by {\rm (\ref{chain10100}),}
$[x]$ is an integer part of a real number $x,$
$\prod\limits_{\emptyset}
\stackrel{\sf def}{=}1,$ $\sum\limits_{\emptyset}
\stackrel{\sf def}{=}0;$ 
another notations are the same as in Theorems~{\rm 1, 2.}}

\vspace{1mm}

{\bf Remark~3.}\ {\it It should be noted that the formulas {\rm (\ref{new1010}),}
{\rm (\ref{new100000}),} {\rm (\ref{new1601}),} {\rm (\ref{new1601a})}
follow from {\rm (\ref{leto6000xxa}).}
It is only necessary to set the values
of the corresponding indicators of the form ${\bf 1}_A$ from the formula
{\rm (\ref{leto6000xxa})} equal to $0$ or $1.$}

{\bf Proof.}\ The proof of Theorem~5 is carried out
by induction using the following recurrence relation

\vspace{-1mm}
$$
J''[\phi_{j_1}\ldots\phi_{j_k}]^{(i_1 \ldots i_k)}_{T,t}=
J''[\phi_{j_k}]^{(i_k)}_{T,t}\cdot
J''[\phi_{j_1}\ldots\phi_{j_{k-1}}]^{(i_1 \ldots i_{k-1})}_{T,t}-
$$
\begin{equation}
\label{recur1}
-\sum\limits_{l=1}^{k-1}{\bf 1}_{\{i_l=i_k\ne 0\}}
{\bf 1}_{\{j_l=j_k\}}\cdot 
J''[\phi_{j_1}\ldots\phi_{j_{l-1}}\phi_{j_{l+1}}\ldots\phi_{j_{k-1}}]^{(i_1
\ldots  i_{l-1}i_{l+1}\ldots i_{k-1})}_{T,t}\ \ \ \hbox{w.~p.~1.}
\end{equation}

\vspace{1mm}

Let us prove the recurrence relation (\ref{recur1}).
Using iteratively the It\^{o} formula, the orthonormality of $\{\phi_j(x)\}_{j=0}^{\infty}$,
as well as (\ref{new5001}) and 
combinatorial reasoning, we obtain w.~p.~1 (see Remark~4 below for details)

$$
J''[\phi_{j_k}]^{(i_k)}_{T,t}\cdot
J''[\phi_{j_1}\ldots\phi_{j_{k-1}}]^{(i_1 \ldots i_{k-1})}_{T,t}=
$$

\vspace{-3mm}
$$
=\int\limits_t^{T}\phi_{j_k}(\theta)d{\bf w}_{\theta}^{(i_k)}
\sum\limits_{(j_1,\ldots,j_{k-1})}
\int\limits_t^T \phi_{j_{k-1}}(t_{k-1})\ldots \int\limits_t^{t_2}\phi_{j_1}(t_1)
d{\bf w}_{t_1}^{(i_1)}\ldots d{\bf w}_{t_{k-1}}^{(i_{k-1})}=
$$

\vspace{-2mm}
$$
=\sum\limits_{(j_1,\ldots,j_{k-1})}\int\limits_t^{T}\phi_{j_k}(\theta)d{\bf w}_{\theta}^{(i_k)}
\int\limits_t^T \phi_{j_{k-1}}(t_{k-1})\ldots \int\limits_t^{t_2}\phi_{j_1}(t_1)
d{\bf w}_{t_1}^{(i_1)}\ldots d{\bf w}_{t_{k-1}}^{(i_{k-1})}=
$$

\vspace{-2mm}
$$
=\sum\limits_{(j_1,\ldots,j_{k})}
\int\limits_t^T \phi_{j_{k}}(t_{k})\ldots \int\limits_t^{t_2}\phi_{j_1}(t_1)
d{\bf w}_{t_1}^{(i_1)}\ldots d{\bf w}_{t_{k}}^{(i_{k})}+
$$

\vspace{-2mm}
$$
+\sum\limits_{(j_1,\ldots,j_{k-1})}
\Biggl({\bf 1}_{\{i_k=i_{k-1}\ne 0\}}
\int\limits_t^T \phi_{j_{k}}(\theta)\phi_{j_{k-1}}(\theta)
\int\limits_t^{\theta}
\phi_{j_{k-2}}(t_{k-2})
\ldots \int\limits_t^{t_2}\phi_{j_1}(t_1)\times\Biggr.
$$

$$
\times
d{\bf w}_{t_1}^{(i_1)}\ldots d{\bf w}_{t_{k-2}}^{(i_{k-2})}d{\bf w}_{\theta}^{(0)}+
$$

\vspace{-2mm}
$$
+{\bf 1}_{\{i_k=i_{k-2}\ne 0\}}
\int\limits_t^T \phi_{j_{k-1}}(t_{k-1})
\int\limits_t^{t_{k-1}} \phi_{j_{k}}(\theta)\phi_{j_{k-2}}(\theta)
\int\limits_t^{\theta}
\phi_{j_{k-3}}(t_{k-3})
\ldots \int\limits_t^{t_2}\phi_{j_1}(t_1)\times\Biggr.
$$

$$
\times
d{\bf w}_{t_1}^{(i_1)}\ldots d{\bf w}_{t_{k-3}}^{(i_{k-3})}d{\bf w}_{\theta}^{(0)}
d{\bf w}_{t_{k-1}}^{(i_{k-1})}+ \ldots
$$

$$
\ldots +{\bf 1}_{\{i_k=i_{1}\ne 0\}}
\int\limits_t^T \phi_{j_{k-1}}(t_{k-1})\ldots 
\int\limits_t^{t_3}
\phi_{j_{2}}(t_{2})
\int\limits_t^{t_{2}} \phi_{j_{k}}(\theta)\phi_{j_{1}}(\theta)\times
$$
$$
\Biggl.\times
d{\bf w}_{\theta}^{(0)}d{\bf w}_{t_2}^{(i_2)}\ldots d{\bf w}_{t_{k-1}}^{(i_{k-1})}
\Biggr)=
$$

\vspace{-1mm}
$$
=\sum\limits_{(j_1,\ldots,j_{k})}
\int\limits_t^T \phi_{j_{k}}(t_{k})\ldots \int\limits_t^{t_2}\phi_{j_1}(t_1)
d{\bf w}_{t_1}^{(i_1)}\ldots d{\bf w}_{t_{k}}^{(i_{k})}+
$$

\vspace{-2mm}
$$
+\sum\limits_{(j_1,\ldots,j_{k-2})}
{\bf 1}_{\{i_k=i_{k-1}\ne 0\}}
\Biggl\{\int\limits_t^T \phi_{j_{k}}(\theta)\phi_{j_{k-1}}(\theta)
\int\limits_t^{\theta}
\phi_{j_{k-2}}(t_{k-2})
\ldots \int\limits_t^{t_2}\phi_{j_1}(t_1)\times\Biggr.
$$

$$
\times
d{\bf w}_{t_1}^{(i_1)}\ldots d{\bf w}_{t_{k-2}}^{(i_{k-2})}d{\bf w}_{\theta}^{(0)}+\ldots
$$

\vspace{-2mm}
$$
\Biggl.\ldots +\int\limits_t^T 
\phi_{j_{k-2}}(t_{k-2})
\ldots \int\limits_t^{t_2}\phi_{j_1}(t_1)
\int\limits_t^{t_1}
\phi_{j_{k}}(\theta)\phi_{j_{k-1}}(\theta)
d{\bf w}_{\theta}^{(0)}
d{\bf w}_{t_1}^{(i_1)}\ldots d{\bf w}_{t_{k-2}}^{(i_{k-2})}\Biggr\}+
$$

\vspace{-2mm}
$$
+\sum\limits_{(j_1,\ldots,j_{k-3},j_{k-1})}
{\bf 1}_{\{i_k=i_{k-2}\ne 0\}}
\Biggl\{\int\limits_t^T \phi_{j_{k}}(\theta)\phi_{j_{k-2}}(\theta)
\int\limits_t^{\theta}
\phi_{j_{k-1}}(t_{k-1})
\int\limits_t^{t_{k-1}}
\phi_{j_{k-3}}(t_{k-3})
\ldots \Biggr.
$$

$$
\ldots \int\limits_t^{t_2}\phi_{j_1}(t_1)
d{\bf w}_{t_1}^{(i_1)}\ldots d{\bf w}_{t_{k-3}}^{(i_{k-3})}d{\bf w}_{t_{k-1}}^{(i_{k-1})}
d{\bf w}_{\theta}^{(0)}+\ldots
$$

\vspace{-2mm}
$$
\ldots +\int\limits_t^T 
\phi_{j_{k-1}}(t_{k-1})
\int\limits_t^{t_{k-1}}
\phi_{j_{k-3}}(t_{k-3})
\ldots \int\limits_t^{t_2}\phi_{j_1}(t_1)
\int\limits_t^{t_1}
\phi_{j_{k}}(\theta)\phi_{j_{k-2}}(\theta)\times
$$

$$
\Biggl.\times d{\bf w}_{\theta}^{(0)}
d{\bf w}_{t_1}^{(i_1)}\ldots d{\bf w}_{t_{k-3}}^{(i_{k-3})}d{\bf w}_{t_{k-1}}^{(i_{k-1})}\Biggr\}+\ldots
$$

$$
\ldots +\sum\limits_{(j_2,\ldots,j_{k-1})}
{\bf 1}_{\{i_k=i_{1}\ne 0\}}
\Biggl\{\int\limits_t^T \phi_{j_{k}}(\theta)\phi_{j_{1}}(\theta)
\int\limits_t^{\theta}
\phi_{j_{k-1}}(t_{k-1})
\ldots \int\limits_t^{t_3}\phi_{j_2}(t_2)\times\Biggr.
$$

$$
\times
d{\bf w}_{t_2}^{(i_2)}\ldots d{\bf w}_{t_{k-1}}^{(i_{k-1})}d{\bf w}_{\theta}^{(0)}+\ldots
$$
$$
\Biggl.\ldots +\int\limits_t^T 
\phi_{j_{k-1}}(t_{k-1})
\ldots \int\limits_t^{t_3}\phi_{j_2}(t_2)
\int\limits_t^{t_2}
\phi_{j_{k}}(\theta)\phi_{j_{1}}(\theta)
d{\bf w}_{\theta}^{(0)}
d{\bf w}_{t_2}^{(i_2)}\ldots d{\bf w}_{t_{k-1}}^{(i_{k-1})}\Biggr\}=
$$

\vspace{1mm}
$$
=\sum\limits_{(j_1,\ldots,j_{k})}
\int\limits_t^T \phi_{j_{k}}(t_{k})\ldots \int\limits_t^{t_2}\phi_{j_1}(t_1)
d{\bf w}_{t_1}^{(i_1)}\ldots d{\bf w}_{t_{k}}^{(i_{k})}+
$$

\vspace{-2mm}
$$
+\int\limits_t^T \phi_{j_{k}}(\theta)\phi_{j_{k-1}}(\theta)d\theta
\sum\limits_{(j_1,\ldots,j_{k-2})}
{\bf 1}_{\{i_k=i_{k-1}\ne 0\}}
\int\limits_t^{T}
\phi_{j_{k-2}}(t_{k-2})
\ldots \int\limits_t^{t_2}\phi_{j_1}(t_1)\times\Biggr.
$$

$$
\times
d{\bf w}_{t_1}^{(i_1)}\ldots d{\bf w}_{t_{k-2}}^{(i_{k-2})}+
$$

\vspace{-4mm}
$$
+\int\limits_t^T \phi_{j_{k}}(\theta)\phi_{j_{k-2}}(\theta)d\theta
\sum\limits_{(j_1,\ldots,j_{k-3},j_{k-1})}
{\bf 1}_{\{i_k=i_{k-2}\ne 0\}}
\int\limits_t^{T}
\phi_{j_{k-1}}(t_{k-1})
\int\limits_t^{t_{k-1}}
\phi_{j_{k-3}}(t_{k-3})
\ldots 
$$

$$
\ldots \int\limits_t^{t_2}\phi_{j_1}(t_1)
d{\bf w}_{t_1}^{(i_1)}\ldots d{\bf w}_{t_{k-3}}^{(i_{k-3})}d{\bf w}_{t_{k-1}}^{(i_{k-1})}
+\ldots
$$

\vspace{-2mm}

$$
\ldots +\int\limits_t^T \phi_{j_{k}}(\theta)\phi_{j_{1}}(\theta)d\theta
\sum\limits_{(j_2,\ldots,j_{k-1})}
{\bf 1}_{\{i_k=i_{1}\ne 0\}}
\int\limits_t^{T}
\phi_{j_{k-1}}(t_{k-1})
\ldots \int\limits_t^{t_3}\phi_{j_2}(t_2)\times
$$

$$
\times
d{\bf w}_{t_2}^{(i_2)}\ldots d{\bf w}_{t_{k-1}}^{(i_{k-1})}=
$$

\vspace{-3mm}

$$
=J''[\phi_{j_1}\ldots\phi_{j_k}]^{(i_1 \ldots i_k)}_{T,t}+
{\bf 1}_{\{i_k=i_{k-1}\ne 0\}}{\bf 1}_{\{j_k=j_{k-1}\}}
\cdot J''[\phi_{j_1}\ldots \phi_{j_{k-2}}]_{T,t}^{(i_1\ldots i_{k-2})}+
$$

\vspace{-1mm}
$$
+
{\bf 1}_{\{i_k=i_{k-2}\ne 0\}}{\bf 1}_{\{j_k=j_{k-2}\}}
\cdot J''[\phi_{j_1}\ldots \phi_{j_{k-3}}\phi_{j_{k-1}}]_{T,t}^{(i_1\ldots i_{k-3}i_{k-1})}+\ldots
$$

\vspace{-1mm}
$$
\ldots +
{\bf 1}_{\{i_k=i_{1}\ne 0\}}{\bf 1}_{\{j_k=j_{1}\}}
\cdot J''[\phi_{j_2}\ldots \phi_{j_{k-1}}]_{T,t}^{(i_2\ldots i_{k-1})}=
$$

\vspace{-1mm}
$$
=J''[\phi_{j_1}\ldots\phi_{j_k}]^{(i_1 \ldots i_k)}_{T,t}+
$$
\begin{equation}
\label{new00002}
~~~~~~+\sum\limits_{l=1}^{k-1}{\bf 1}_{\{i_l=i_k\ne 0\}}
{\bf 1}_{\{j_l=j_k\}}\cdot 
J''[\phi_{j_1}\ldots\phi_{j_{l-1}}\phi_{j_{l+1}}\ldots\phi_{j_{k-1}}]^{(i_1
\ldots  i_{l-1}i_{l+1}\ldots i_{k-1})}_{T,t}.
\end{equation}

\vspace{2mm}

The equality (\ref{recur1}) is proved. Theorem~5 is proved.

\vspace{2mm}

{\bf Remark~4.}\ {\it It should be noted that the sums with respect to 
permutations 
$$
\sum\limits_{(j_1,\ldots,j_{k-1})}
$$

\noindent
in {\rm (\ref{new00002}),} containing the expressions 
$$
{\bf 1}_{\{i_k=i_{k-1}\ne 0\}}\phi_{j_{k}}(\theta)\phi_{j_{k-1}}(\theta),\ldots,
{\bf 1}_{\{i_k=i_{1}\ne 0\}}\phi_{j_{k}}(\theta)\phi_{j_{1}}(\theta),
$$

\noindent
should be understood in a special way.
Let us explain this rule on following sum
$$
\sum\limits_{(j_1,\ldots,j_{k-1})}
{\bf 1}_{\{i_k=i_{k-1}\ne 0\}}
\int\limits_t^T \phi_{j_{k}}(\theta)\phi_{j_{k-1}}(\theta)
\int\limits_t^{\theta}
\phi_{j_{k-2}}(t_{k-2})
\ldots \int\limits_t^{t_2}\phi_{j_1}(t_1)\times
$$

\begin{equation}
\label{new00003}
\times
d{\bf w}_{t_1}^{(i_1)}\ldots d{\bf w}_{t_{k-2}}^{(i_{k-2})}d{\bf w}_{\theta}^{(0)}.
\end{equation}

More precisely, permutations $(j_1,\ldots,j_{k-1})$ 
when summing in {\rm (\ref{new00003})}
are performed in such a way that if
$j_r$ swapped with $j_d$ in the  
permutation $(j_1,\ldots,j_{k-1}),$ 
then $i_r$ swapped with $i_d$ in 
the permutation $(i_1,\ldots,i_{k-2} i_{k-1})$\ {\rm (}note that $i_{k-1}=0${\rm )}.
Moreover, 
$\bar \phi_{j_r}$ swapped with $\bar \phi_{j_d}$
in the permutation 
$$
\bigl(\bar \phi_{j_{1}},\ldots,\bar\phi_{j_{k-1}}\bigr)=
\bigl(\phi_{j_1},\ldots,\phi_{j_{k-2}},\hspace{1.5mm} 
{\bf 1}_{\{i_k=i_{k-1}\ne 0\}}\cdot \phi_{j_k}\cdot \phi_{j_{k-1}}\bigr),
$$

\noindent
where $\bar\phi_{j_{k-1}}(\tau)=
{\bf 1}_{\{i_k=i_{k-1}\ne 0\}}\phi_{j_k}(\tau)\phi_{j_{k-1}}(\tau).$

A similar rule should be applied to all other sums with respect to permutations
$$
\sum\limits_{(j_1,\ldots,j_{k-1})}
$$

\noindent
in {\rm (\ref{new00002})} that contain the expressions
$$
{\bf 1}_{\{i_k=i_{k-2}\ne 0\}}\phi_{j_{k}}(\theta)\phi_{j_{k-2}}(\theta),\ldots,
{\bf 1}_{\{i_k=i_{1}\ne 0\}}\phi_{j_{k}}(\theta)\phi_{j_{1}}(\theta).
$$
}

\vspace{-6mm}

\section{Main Results}

\subsection{Generalizations of Theorem 2 to the Case of an Arbitrary 
Complete Ortho\-nor\-mal Systems of Functions in the Space $L_2([t, T])$
and $\psi_1(\tau),$ $\ldots,\psi_k(\tau)\in L_2([t, T])$}

Theorems~3--5 imply the following two theorems on expansion
of iterated It\^{o} stochastic integrals (\ref{ito}).

{\bf Theorem 6}\ \cite{2018a}, \cite{arxiv-1}.
{\it Suppose that
the condition {\rm ($\star$)} is fulfilled
for the multi-index $(i_1 \ldots i_k)$ {\rm (}see Sect.~{\rm 2.2)} 
and the condition {\rm (\ref{ziko999})} is also 
fulfilled.
Furthermore$,$ let 
$\psi_1(\tau),\ldots,\psi_k(\tau)\in L_2([t, T])$ and
$\{\phi_j(x)\}_{j=0}^{\infty}$ is an arbitrary complete orthonormal system  
of functions in the space $L_2([t,T]).$
Then the following expansion

\vspace{-7mm}
$$
J[\psi^{(k)}]_{T,t}^{(i_1\ldots i_k)}=
\hbox{\vtop{\offinterlineskip\halign{
\hfil#\hfil\cr
{\rm l.i.m.}\cr
$\stackrel{}{{}_{p_1,\ldots,p_k\to \infty}}$\cr
}} }
\sum\limits_{j_1=0}^{p_1}\ldots
\sum\limits_{j_k=0}^{p_k}
C_{j_k\ldots j_1}\times
$$

\vspace{-2mm}
\begin{equation}
\label{new9999}
\times
\prod_{l=1}^k\left({\bf 1}_{\{m_l=0\}}+{\bf 1}_{\{m_l>0\}}\left\{
\begin{matrix}
H_{n_{1,l}}\left(\zeta_{j_{h_{1,l}}}^{(i_l)}\right)\ldots 
H_{n_{d_l,l}}\left(\zeta_{j_{h_{d_l,l}}}^{(i_l)}\right),\ 
&\hbox{\rm if}\ \ \ 
i_l\ne 0\cr\cr
\left(\zeta_{j_{h_{1,l}}}^{(0)}\right)^{n_{1,l}}\ldots
\left(\zeta_{j_{h_{d_l,l}}}^{(0)}\right)^{n_{d_l,l}},\  &\hbox{\rm if}\ \ \ 
i_l=0
\end{matrix}\right.\ \right)
\end{equation}

\vspace{1mm}
\noindent
con\-verg\-ing in the mean-square sense is valid$,$
where 
$H_n(x)$ is the Hermite polynomial {\rm (\ref{ziko500}),}
${\bf 1}_A$ is the indicator of the set $A,$
$i_1,\ldots,i_k=0, 1,\ldots,m;$\ \
$n_{1,l}+n_{2,l}+\ldots+n_{d_l,l}=m_l;$\ \ $n_{1,l}, n_{2,l}, \ldots, n_{d_l,l}=1,\ldots, m_l;$\ \ 
$d_l=1,\ldots,m_l;$\ \ $l=1,\ldots,k;$\ \ $m_1+\ldots+m_k=k;$\ \  
the numbers $m_1,\ldots,m_k,$\ $g_1,\ldots,g_k$
depend on $(i_1,\ldots,i_k)$ and 
the numbers $n_{1,l},\ldots,n_{d_l,l},$\ $h_{1,l},\ldots,h_{d_l,l},$\ $d_l$
depend on $\{j_1,\ldots,j_k\};$ moreover$,$ $\left\{j_{g_1},\ldots,j_{g_k}\right\}
=\{j_1,\ldots,j_k\};$
$$
\zeta_{j}^{(i)}=
\int\limits_t^T \phi_{j}(\tau) d{\bf w}_{\tau}^{(i)}\ \ \ (i=0,1,\ldots,m;\ \ j=0,1,2,\ldots)
$$
are independent standard Gaussian random variables
for various
$i$ or $j$ {\rm(}in the case when $i\ne 0${\rm)}
and $d{\bf w}_{\tau}^{(0)}=d\tau;$ another notations as in Theorems~{\rm 1, 2.}}

\newpage
\noindent
\par
{\bf Theorem 7}\ \cite{2018a}, \cite{arxiv-1}.
{\it Suppose that
$\psi_1(\tau),\ldots,\psi_k(\tau)\in L_2([t, T])$ and
$\{\phi_j(x)\}_{j=0}^{\infty}$ is an arbitrary complete orthonormal system  
of functions in the space $L_2([t,T]).$
Then the following expansion
$$
J[\psi^{(k)}]_{T,t}^{(i_1\ldots i_k)}=
\hbox{\vtop{\offinterlineskip\halign{
\hfil#\hfil\cr
{\rm l.i.m.}\cr
$\stackrel{}{{}_{p_1,\ldots,p_k\to \infty}}$\cr
}} }
\sum\limits_{j_1=0}^{p_1}\ldots
\sum\limits_{j_k=0}^{p_k}
C_{j_k\ldots j_1}\Biggl(
\prod_{l=1}^k\zeta_{j_l}^{(i_l)}+\sum\limits_{r=1}^{[k/2]}
(-1)^r \times
\Biggr.
$$

\vspace{-3mm}
\begin{equation}
\label{razzar1}
\times
\sum_{\stackrel{(\{\{g_1, g_2\}, \ldots, 
\{g_{2r-1}, g_{2r}\}\}, \{q_1, \ldots, q_{k-2r}\})}
{{}_{\{g_1, g_2, \ldots, 
g_{2r-1}, g_{2r}, q_1, \ldots, q_{k-2r}\}=\{1, 2, \ldots, k\}}}}
\prod\limits_{s=1}^r
{\bf 1}_{\{i_{g_{{}_{2s-1}}}=~i_{g_{{}_{2s}}}\ne 0\}}
\Biggl.{\bf 1}_{\{j_{g_{{}_{2s-1}}}=~j_{g_{{}_{2s}}}\}}
\prod_{l=1}^{k-2r}\zeta_{j_{q_l}}^{(i_{q_l})}\Biggr)
\end{equation}

\vspace{1mm}
\noindent
con\-verg\-ing in the mean-square sense is valid$,$
where 
$[x]$ is an integer part of a real number $x;$
another notations are the same as in Theorems~{\rm 1, 2, 5.}}

\subsection{Modifications of Theorems 6, 7 for the Case of an Arbitrary 
Complete Ortho\-nor\-mal Systems of Functions in the Space $L_2([t, T])$
and $\Phi(t_1,\ldots,t_k)\in L_2([t, T]).$}

Replacing the function $K(t_1,\ldots,t_k)$ of the form (\ref{ppp}) in Theorems~6, 7 by 
the function $\Phi(t_1,\ldots,t_k)\in L_2([t,T]),$ we get the following two theorems.

\vspace{2mm}

{\bf Theorem 8}\ \cite{2018a}, \cite{arxiv-1}.
{\it Suppose that
the condition {\rm ($\star$)} is fulfilled
for the multi-index $(i_1 \ldots i_k)$ {\rm (}see Sect.~{\rm 2.2)} 
and the condition {\rm (\ref{ziko999})} is also 
fulfilled.
Furthermore$,$ let 
$\Phi(t_1,\ldots,t_k)\in L_2([t, T])$ and
$\{\phi_j(x)\}_{j=0}^{\infty}$ is an arbitrary complete orthonormal system  
of functions in the space $L_2([t,T]).$
Then the following expansion
$$
J''[\Phi]_{T,t}^{(i_1\ldots i_k)}=
\hbox{\vtop{\offinterlineskip\halign{
\hfil#\hfil\cr
{\rm l.i.m.}\cr
$\stackrel{}{{}_{p_1,\ldots,p_k\to \infty}}$\cr
}} }
\sum\limits_{j_1=0}^{p_1}\ldots
\sum\limits_{j_k=0}^{p_k}
C_{j_k\ldots j_1}\times
$$

\vspace{-2mm}
\begin{equation}
\label{new9999x}
\times
\prod_{l=1}^k\left({\bf 1}_{\{m_l=0\}}+{\bf 1}_{\{m_l>0\}}\left\{
\begin{matrix}
H_{n_{1,l}}\left(\zeta_{j_{h_{1,l}}}^{(i_l)}\right)\ldots 
H_{n_{d_l,l}}\left(\zeta_{j_{h_{d_l,l}}}^{(i_l)}\right),\ 
&\hbox{\rm if}\ \ \ 
i_l\ne 0\cr\cr
\left(\zeta_{j_{h_{1,l}}}^{(0)}\right)^{n_{1,l}}\ldots
\left(\zeta_{j_{h_{d_l,l}}}^{(0)}\right)^{n_{d_l,l}},\  &\hbox{\rm if}\ \ \ 
i_l=0
\end{matrix}\right.\ \right)
\end{equation}

\vspace{2mm}
\noindent
con\-verg\-ing in the mean-square sense is valid$,$
where the sum of iterated It\^{o}
stochastic integrals $J''[\Phi]_{T,t}^{(i_1\ldots i_k)}$ is defined by {\rm (\ref{chain10100}),}

\vspace{-2mm}
\begin{equation}
\label{ppppaxx}
C_{j_k\ldots j_1}=\int\limits_{[t,T]^k}
\Phi(t_1,\ldots,t_k)\prod_{l=1}^{k}\phi_{j_l}(t_l)dt_1\ldots dt_k
\end{equation}

\vspace{2mm}
\noindent
is the Fourier coefficient$,$ $H_n(x)$ is the Hermite polynomial {\rm (\ref{ziko500}),}
${\bf 1}_A$ is the indicator of the set $A,$
$i_1,\ldots,i_k=0, 1,\ldots,m;$\ \
$n_{1,l}+n_{2,l}+\ldots+n_{d_l,l}=m_l;$\ \ $n_{1,l}, n_{2,l}, \ldots, n_{d_l,l}=1,\ldots, m_l;$\ \ 
$d_l=1,\ldots,m_l;$\ \ $l=1,\ldots,k;$\ \ $m_1+\ldots+m_k=k;$\ \  
the numbers $m_1,\ldots,m_k,$\ $g_1,\ldots,g_k$
depend on $(i_1,\ldots,i_k)$ and 
the numbers $n_{1,l},\ldots,n_{d_l,l},$\ $h_{1,l},\ldots,h_{d_l,l},$\ $d_l$
depend on $\{j_1,\ldots,j_k\};$ moreover$,$ $\left\{j_{g_1},\ldots,j_{g_k}\right\}
=\{j_1,\ldots,j_k\};$

\vspace{-2mm}
$$
\zeta_{j}^{(i)}=
\int\limits_t^T \phi_{j}(\tau) d{\bf w}_{\tau}^{(i)}\ \ \ (i=0,1,\ldots,m;\ \ j=0,1,2,\ldots)
$$

\vspace{2mm}
\noindent
are independent standard Gaussian random variables
for various
$i$ or $j$ {\rm(}in the case when $i\ne 0${\rm)}
and $d{\bf w}_{\tau}^{(0)}=d\tau.$
}

\vspace{2mm}

{\bf Theorem 9}\ \cite{2018a}, \cite{arxiv-1}.\ {\it Suppose that
$\Phi(t_1,\ldots,t_k)\in L_2([t, T]^k)$ 
and
$\{\phi_j(x)\}_{j=0}^{\infty}$ is an arbitrary complete orthonormal system  
of functions in the 
space $L_2([t,T]).$
Then the following expansion

\vspace{-5mm}
$$
J''[\Phi]_{T,t}^{(i_1\ldots i_k)}=
\hbox{\vtop{\offinterlineskip\halign{
\hfil#\hfil\cr
{\rm l.i.m.}\cr
$\stackrel{}{{}_{p_1,\ldots,p_k\to \infty}}$\cr
}} }
\sum\limits_{j_1=0}^{p_1}\ldots
\sum\limits_{j_k=0}^{p_k}
C_{j_k\ldots j_1}\Biggl(
\prod_{l=1}^k\zeta_{j_l}^{(i_l)}+\sum\limits_{r=1}^{[k/2]}
(-1)^r \times
\Biggr.
$$

\vspace{-1mm}
$$
\times
\sum_{\stackrel{(\{\{g_1, g_2\}, \ldots, 
\{g_{2r-1}, g_{2r}\}\}, \{q_1, \ldots, q_{k-2r}\})}
{{}_{\{g_1, g_2, \ldots, 
g_{2r-1}, g_{2r}, q_1, \ldots, q_{k-2r}\}=\{1, 2, \ldots, k\}}}}
\prod\limits_{s=1}^r
{\bf 1}_{\{i_{g_{{}_{2s-1}}}=~i_{g_{{}_{2s}}}\ne 0\}}
\Biggl.{\bf 1}_{\{j_{g_{{}_{2s-1}}}=~j_{g_{{}_{2s}}}\}}
\prod_{l=1}^{k-2r}\zeta_{j_{q_l}}^{(i_{q_l})}\Biggr)
$$

\vspace{3mm}
\noindent
con\-verg\-ing in the mean-square sense is valid$,$
where the sum of iterated It\^{o}
stochastic integrals $J''[\Phi]_{T,t}^{(i_1\ldots i_k)}$ is defined by {\rm (\ref{chain10100}),}
the Fourier coefficient
$C_{j_k\ldots j_1}$ has the form 
{\rm (\ref{ppppaxx});}
another
notations are the same as in Theorems~{\rm 1, 2, 5}.}

\section{Comparison with Other Results and Conclusions}

Before starting this section, we recall that
the sum of iterated It\^{o} stochastic integrals (\ref{chain10100}),
which plays a central role in the proofs of Theorems~6--9, is equal w.~p.~1
to the multiple Wiener stochastic integral with respect
to the components of a multidimensional Wiener process
(see the proof in \cite{2018a}, Sect.~1.11).

It should be noted that an analogue of Theorem 8 (more precisely, the expansion
like (\ref{new9999x}) for the case $i_1,\ldots,i_k=1,\ldots,m$) was obtained
in \cite{Rybakov1000}. The mentioned expansion is formulated in 
\cite{Rybakov1000} using the multiple Wiener stochastic integral and the Wick product.
Also note that the proof in \cite{Rybakov1000} is 
different from the proof given in this article.         
Let us describe these differences.

In \cite{Rybakov1000}, the author interprets the multiple Wiener
stochastic integral from a finite-dimensional kernel
$K_{p,\ldots,p}(t_1,\ldots,t_k)$ of the form (\ref{chain30001x})
as a linear operator and proves that this operator is bounded.
We note that the proof from \cite{Rybakov1000}
is essentially based on Theorem~3.1 from \cite{ito1951}.

In our proof of Theorems~6--9 we use
the sum of iterated It\^{o} stochastic integrals (\ref{chain10100})
several times and do not explicitly use
the multiple Wiener stochastic integral.
Moreover, our proof of Theorems~6--9 is based on the It\^{o} formula 
and does not use Theorem~3.1 from \cite{ito1951}.
The methodology of our proof is a direct development
of the approach we used to prove Theorem~5.1 in \cite{2006} (2006).

Note that the results of \cite{Rybakov1000}, as well as 
the results of this article, are based on our idea 
\cite{1} (2006) on the expansion of the kernel (\ref{ppp}) (or $\Phi(t_1,\ldots,t_k)\in L_2([t,T]^k)$)
into a generalized multiple Fourier series 
(see \cite{1}, Chapter~5, Theorem~5.1, pp.~235-245 
or \cite{2018a}, Chapter~1 for details). 

We also note a number of works \cite{ito1951}, \cite{Kuo}-\cite{major2} in which the properties
of multiple Wiener stochastic integrals were studied using
measure theory, in particular, the formulas for the product
of such integrals were obtained.

First of all, let us compare Theorem~5 with Proposition~5.1 from \cite{fox}.
An analogue of the right-hand side of (\ref{leto6000xxa})
for nonrandom $x_1,\ldots,x_k$
is constructed in \cite{fox} using diagrams (see the formula (5.1) in \cite{fox}).
This means that the application of the formula (5.1) from \cite{fox},
unlike the formula (\ref{leto6000xxa}), is difficult when
performing algebraic transformations.

Further, we note that the formula (5.1) from \cite{fox}
was applied to the representation of the multiple Wiener stochastic integral
somewhat differently than the formula (\ref{leto6000xxa}).
Namely, using Proposition~5.1 \cite{fox}.
Let us expain this difference in more detail.

Proposition~5.1 from \cite{fox} in our degree of generality 
and in our notations can be written as

\vspace{-2mm}
$$
J''\left[\phi_{j_1}\ldots \phi_{j_k}\right]_{T,t}^{(i_1\ldots i_k)}=
$$

\vspace{-2mm}
$$
=
J''\biggl[\underbrace{\phi_{j_1}
\ldots \phi_{j_1}}_{m_1}
\underbrace{\phi_{j_2}
\ldots \phi_{j_2}}_{m_2}\ldots 
\underbrace{\phi_{j_p}\ldots
\phi_{j_p}}_{m_p}\biggr]_{T,t}^
{(\overbrace{\hspace{0.5mm}{}_{i_1 \ldots i_{m_1}}}^{m_1}
\overbrace{\hspace{0.3mm}{}_{i_{m_1+1} \ldots i_{m_2}}}^{m_2}
\ldots \overbrace{\hspace{0.3mm}{}_{i_{m_1+\ldots +m_{p-1}+1} \ldots i_k}}^{m_p})}=
$$

\vspace{-2mm}
\begin{equation}
\label{new54321}
=J''\left[\phi_{j_1}
\hspace{-0.3mm}\ldots \hspace{-0.3mm}\phi_{j_1}\right]
_{T,t}^
{(\overbrace{\hspace{0.5mm}{}_{i_1 \ldots i_{m_1}}}^{m_1})}\hspace{-0.3mm}\cdot
J''\left[\phi_{j_2}\hspace{-0.3mm}
\ldots \hspace{-0.3mm}\phi_{j_2}\right]_{T,t}^
{(\overbrace{\hspace{0.3mm}{}_{i_{m_1+1} \ldots i_{m_2}}}^{m_2})}\hspace{-0.3mm}\cdot \ldots
\cdot J''\left[\phi_{j_p}
\hspace{-0.3mm}\ldots \hspace{-0.3mm}\phi_{j_p}\right]_{T,t}^
{(\overbrace{\hspace{0.3mm}{}_{i_{m_1+\ldots +m_{p-1}+1} \ldots i_k}}^{m_p})}
\end{equation}

\vspace{5mm}
\noindent
w.~p.~1, where

\vspace{-2mm}
$$
J''\left[\phi_{j_1}
\hspace{-0.3mm}\ldots \hspace{-0.3mm}\phi_{j_1}\right]
_{T,t}^
{(\overbrace{\hspace{0.5mm}{}_{i_1 \ldots i_{m_1}}}^{m_1})}\hspace{-1.5mm},
J''\left[\phi_{j_2}\hspace{-0.3mm}
\ldots \hspace{-0.3mm}\phi_{j_2}\right]_{T,t}^
{(\overbrace{\hspace{0.3mm}{}_{i_{m_1+1} \ldots i_{m_2}}}^{m_2})}\hspace{-1mm},\ldots,
J''\left[\phi_{j_p}
\hspace{-0.3mm}\ldots \hspace{-0.3mm}\phi_{j_p}\right]_{T,t}^
{(\overbrace{\hspace{0.3mm}{}_{i_{m_1+\ldots +m_{p-1}+1} \ldots i_k}}^{m_p})}
$$

\vspace{3.5mm}
\noindent
are defined by the right-hand side of the formula (5.1) from \cite{fox},
$m_1+\ldots +m_p=k,$ $m_1,\ldots, m_p>0,$ $j_q\ne j_d$ $(q\ne d,\ q,d=1,\ldots,p),$
$i_1,\ldots,i_k=1,\ldots,m.$

This actually means that in \cite{fox} an analogue of the formula
(\ref{leto6000xxa}) is constructed for the special case
$j_1=\ldots=j_k$. Moreover, the specified analogue 
is based on the formula (5.1) \cite{fox} obtained using diagrams.

Comparing the formulas (\ref{leto6000xxa}) and (\ref{new54321}) (or (5.1) from \cite{fox}), it is easy
to understand that the transition from 
(\ref{leto6000xxa}) and (\ref{new54321}) is obvious.
It is only necessary to set the values
of the corresponding indicators of the form ${\bf 1}_A$ from the formula
(\ref{leto6000xxa}) equal to $0$ or $1.$
The reverse transition from the formula (\ref{new54321})
to the formula (\ref{leto6000xxa}) is not obvious.
Note that the formula 
(\ref{leto6000xxa}) (not the formula (\ref{new54321})) is convenient for the  numerical
integration of It\^{o} stochastic differential equations (see \cite{2018a}, 
Chapter~5 for details).

Let us turn to the comparison of Theorem~5 with another interesting work \cite{major2} (2019).
As it turned out, a version of Theorem~5 was obtained in terms 
of Wick polynomials and for the case of vector valued random measures 
in \cite{major2} (see Theorem~7.2, p.~69).
However, much earlier the formula (\ref{leto6000xxa}) (Theorem~5) is obtained
in our monograph \cite{2009} (2009) as part of the formula
(5.30) (see \cite{2009}, p.~220).
Moreover, particular cases of the formula (\ref{leto6000xxa}) were obtained
even earlier in our works \cite{2006} (2006) and \cite{2007} (2007).
More precisely, particular cases $k=1,\ldots,5$ of the formula (\ref{leto6000xxa})
were obtained in \cite{2006} (2006) as parts of the formulas 
on the pages 243-244 and partiular cases $k=1,\ldots,7$ of the formula (\ref{leto6000xxa})
were obtained in \cite{2007} (2007) as parts of the formulas 
on the pages 208-218.

We also note that we have found an explicit expression for the 
Wick polynomial of degree $k$ of the arguments $\zeta_{j_1}^{(i_1)},\ldots,\zeta_{j_k}^{(i_k)}$ 
(see the formula (\ref{leto6000xxa})),
which is very convenient for the numerical simulation of
iterated It\^{o} stochastic integrals (\ref{ito})  \cite{Mikh-1}.
Note that the representation of the Wick polynomial
of the arguments $\zeta_{j_1}^{(i_1)},\ldots,\zeta_{j_k}^{(i_k)}$ 
in terms of the product of Hermite polynomials
is less convenient for the numerical simulation of
iterated It\^{o} stochastic integrals (\ref{ito}).
For example, the expression for $J''[\phi_{j_1}\phi_{j_2}\phi_{j_3}\phi_{j_4}]^{(i_1 i_2 i_3 i_4)}_{T,t}$
in terms of the product of Hermite polynomials,
even under the condition $i_1=i_2=i_3=i_4$, already contains
15 different expressions (see \cite{2018a}, Sect.~1.10).
At the same time, all these 15 expressions are contained 
in one formula (\ref{leto6000xxa}) provided that $k=4$ and $i_1=i_2=i_3=i_4$.
It is very convenient, since in computer simulation
using the formula (\ref{leto6000xxa}), in addition to
modeling of random variables $\zeta_{j_1}^{(i_1)},\ldots,\zeta_{j_k}^{(i_k)}$,
it remains only to set  
the values
of the corresponding indicators of the form ${\bf 1}_A$ from the formula
(\ref{leto6000xxa}) equal to $0$ or $1.$

It should be noted that in \cite{major} (Theorem~6.1)
a diagram formula was obtained for the product
of two multiple Wiener stochastic integrals
with respect to vector valued random measures. 
The formula (\ref{new1600}) can be derived from the diagram formula \cite{major}.
Although the proof of the diagram formula \cite{major}
is much more complicated than our proof of the formula (\ref{new1600}).

To conclude this article, we say a few words about expansions 
(\ref{new9999}) and (\ref{razzar1}).
The transition from the expansion (\ref{razzar1}) to the expansion 
(\ref{new9999}) is obvious. It is only necessary to set the values
of the corresponding indicators of the form ${\bf 1}_A$ from the formula
(\ref{razzar1}) equal to $0$ or $1.$
The reverse transition from the formula (\ref{new9999})
to the formula (\ref{razzar1}) is also possible but not obvious.
However, Theorems~4 and 5 provide a transition from 
(\ref{new9999}) to (\ref{razzar1}) and vice versa.
Note that the expansion (\ref{new9999}) is interesting from the point of 
view of studying the structure of the expansion of iterated It\^{o}
stochastic integrals. On the other hand, 
the expansion (\ref{razzar1}) is exceptionally convenient 
for applications 
\cite{Mikh-1}, \cite{new-new-1}.
For example, in \cite{Mikh-1}, \cite{new-new-1}, approximations 
of iterated It\^{o} stochastic integrals of multiplicities 1 to 6 
in the Python programming 
language were successfully implemented using (\ref{razzar1}) $(k=1,\ldots,6)$ and 
Legendre polynomials.

\linespread{0.7}

\end{document}